\documentclass[aap,preprint]{imsart}

\RequirePackage[OT1]{fontenc}
\RequirePackage{amsthm,amsmath}
\RequirePackage[numbers]{natbib}
\RequirePackage[colorlinks,citecolor=blue,urlcolor=blue]{hyperref}
\usepackage{graphicx,psfrag,epsf}
\usepackage{enumerate}
\usepackage[english]{babel}
\usepackage{amsmath}
\usepackage{amssymb}
\usepackage{marvosym}
\usepackage{amsfonts}
\usepackage{array}
\usepackage{bm,latexsym,mathrsfs}

\startlocaldefs
\numberwithin{equation}{section}
\theoremstyle{plain}
\newtheorem{thm}{Theorem}[section]
\newtheorem{rem}{Remark}[section]
\newtheorem{cor}{Corollary}[section]

\newtheorem{prop}{Proposition}[section]
\newtheorem{lem}{Lemma}[section]

 \newcommand{\ind}{\mathbf{1}}
 \def\proof{\noindent {\bf Proof.}\ }
\def\endproof{{\mbox{}\nolinebreak\hfill\rule{2mm}{2mm}\par\medbreak} }
\endlocaldefs

\begin{document}

\begin{frontmatter}
\title{Central limit Theorem for an Adaptive Randomly Reinforced Urn Model}
\runtitle{CLT for an Adaptive Urn Model}

\begin{aug}
\author{\fnms{Andrea} \snm{Ghiglietti}\thanksref{m1}\ead[label=e1]{andrea.ghiglietti@unimi.it}},
\author{\fnms{Anand} \snm{N. Vidyashankar}\thanksref{m2,t2}\ead[label=e2]{avidyash@gmu.edu}}\\
\and
\author{\fnms{William} \snm{F. Rosenberger}\thanksref{m2}\ead[label=e3]{wrosenbe@gmu.edu}}

\runauthor{Ghiglietti \textit{et al.}}

\affiliation{Universit\`{a} degli Studi di Milano\thanksmark{m1} and George Mason University\thanksmark{m2}}

\address{Department of Mathematics ``F. Enriques''\\
Universit\`{a} degli Studi di Milano\\
via Saldini 50\\
20133 Milan, Italy\\
\printead{e1}}

\address{Department of Statistics\\
George Mason University\\
4400 University Drive MS4A7\\
Fairfax, VA 22030 USA\\
\printead{e2}\\
\printead{e3}}

\thankstext{t2}{Supported in part by NsF grant DMS 000-03-07057.}

\end{aug}

\begin{abstract}
:\ The generalized P\'olya urn (GPU) models and their variants have been investigated in several disciplines.
However, typical assumptions made with respect to the GPU
do not include urn models with diagonal replacement matrix,
which arise in several applications, specifically in clinical trials.
To facilitate mathematical analyses of models in these applications,
we introduce an adaptive randomly reinforced urn model that uses accruing statistical information
to adaptively skew the urn proportion toward specific targets.
We study several probabilistic aspects that are important in implementing the urn model in practice.
Specifically, we establish the law of large numbers and a central limit theorem
for the number of sampled balls.
To establish these results, we develop new techniques involving last exit times
and crossing time analyses of the proportion of balls in the urn.
To obtain precise estimates in these techniques, we establish results
on the harmonic moments of the total number of balls in the urn.
Finally, we describe our main results in the context an application to response-adaptive randomization in clinical trials.
Our simulation experiments in this context demonstrate the ease and scope of our model.
\end{abstract}

\begin{keyword}[class=MSC]
\kwd[Primary ]{60F15}
\kwd{60F05}\kwd{60E20}\kwd{60G99}
\kwd[; secondary ]{68Q87}
\kwd{97K50}
\end{keyword}

\begin{keyword}
\kwd{Clinical trials}
\kwd{crossing times}
\kwd{harmonic moments}
\kwd{last exit times}
\kwd{generalized P\'olya urn}
\kwd{target allocation}
\end{keyword}

\end{frontmatter}

\section{Introduction}   \label{section_introduction}


A generalized P\'olya urn (GPU) model~\cite{Athreya.Karlin.68} is characterized by the pair $(Y_{1,n},Y_{2,n})$ of random variables representing
the number of balls of two colors, red and white, for instance.
The process is described as follows: at time $n=0$, the process starts with $(y_{1,0},y_{2,0})$ balls.
A ball is drawn at random. If the color is red, the ball is returned to the urn along with
the random numbers $(D_{11,1},D_{12,1})$ of red and white balls;
otherwise, the ball is returned to the urn along with
the random numbers $(D_{21,1},D_{22,1})$ of red and white balls, respectively.
Let $Y_{1,1}=y_{1,0}+D_{11,1}$ and $Y_{2,1}=y_{2,0}+D_{12,1}$ denote the urn composition when the sampled ball is red;
similarly, let $Y_{1,1}=y_{1,0}+D_{21,1}$ and $Y_{2,1}=y_{2,0}+D_{22,1}$ denote the urn composition when the sampled ball is white.
The process is repeated yielding the collection $\{(Y_{1,n},Y_{2,n});n\geq1\}$.
The quantities $R_1=\{(D_{11,n},D_{12,n});n\geq1\}$ and $R_2=\{(D_{21,n},D_{22,n});n\geq1\}$ are collections of independent and identically distributed (i.i.d.)
non-negative integer valued random variables, and $R_1$ is assumed to be independent of $R_2$.
We refer to
\[D_{n}=\left[
\begin{array}{ll}
D_{11,n}\ &\ D_{12,n} \\
D_{21,n}\ &\ D_{22,n}\\
\end{array}
\right]\]
as a replacement matrix.

In this paper, we focus on an extension of the randomly reinforced urn (RRU) model,
a variant of the randomized P\'olya urn (RPU) models, whose replacement matrix is given by
\[D_{n}=\left[
\begin{array}{ll}
D_{11,n}\ &\ D_{12,n} \\
D_{21,n}\ &\ D_{22,n}\\
\end{array}
\right]\equiv\left[
\begin{array}{ll}
D_{1,n}\ &\ 0 \\
0\ &\ D_{2,n}\\
\end{array}
\right].
\]
where the random variables $D_{1,n}$ and $D_{2,n}$ are supported on $[0,\infty)$,
rather than on the set of non-negative integers.
Let $m_1:=\bm{E}[D_{1,n}]$ and $m_{2}:=\bm{E}[D_{2,n}]$.
For the RRU model, a law of large numbers was established in~\cite{Muliere.et.al.06}; i.e.
\begin{equation}\label{eq:RRU}
Z_n\ =\ \frac{Y_{1,n}}{Y_{1,n}+Y_{2,n}}\ \stackrel{a.s.}{\rightarrow}\
\begin{cases}
1\cdot\ind_{\{m_1>m_2\}}+0\cdot\ind_{\{m_1<m_2\}}\ &\text{if } m_1\neq m_2,\\
Z_{\infty}\ &\text{if } m_1= m_2,
\end{cases}\end{equation}
where $\stackrel{a.s.}{\rightarrow}$ stands for almost sure convergence and $Z_{\infty}$ is a random variable supported on $(0,1)$.
The properties of the distribution of $Z_{\infty}$ were studied in~\cite{Aletti.et.al.09,Aletti.et.al.12}.
Denoting $\{(N_{1,n},N_{2,n});n\geq1\}$ the number of balls of red and white colors sampled from the urn,
one can deduce from the above LLN that $N_{1,n}/n$ converges to the same limit as $Z_n$.

Notice that the limit of the RRU in~\eqref{eq:RRU} is always 1 or 0 when $m_1\neq m_2$.
However, in applications it is common to target a specific value $\rho\in(0,1)$.
This was achieved in~\cite{Aletti.et.al.13}, where
the modified randomly reinforced urn (MRRU) model was introduced.
The MRRU model is an RRU model with two fixed thresholds $0<\rho_2\leq \rho_1<1$,
such that if $Z_n<\rho_2$, no white balls are replaced in urn,
while if $Z_n>\rho_1$, no red balls are replaced in the urn.
These changes occur at random times and will in general depend on $m_1$ and $m_2$.
Thus, even if the sequences $\{D_{1,n};n\geq1\}$ and $\{D_{2,n};n\geq1\}$ are i.i.d.,
the replacements matrices of the MRRU model are not i.i.d.
Indeed the replacement matrix has the following representation:
\[D_{n}=\left[
\begin{array}{ll}
D_{1,n}\cdot\ind_{\{Z_{n-1}\leq\rho_1\}}\ &\ 0 \\
0\ &\ D_{2,n}\cdot\ind_{\{Z_{n-1}\geq\rho_2\}}\\
\end{array}
\right].
\]
The LLN for the MRRU when $m_1\neq m_2$ is established as
$$Z_n\stackrel{a.s.}{\rightarrow}\rho_1\cdot\ind_{\{m_1>m_2\}}+\rho_2\cdot\ind_{\{m_1<m_2\}}.$$
A second order result for $Z_n$, namely the asymptotic distribution of $Z_n$
after appropriate centering, was derived in~\cite{Ghiglietti.et.al.12}.
\textit{However, the validity of the CLT for $N_{1,n}/n$ in the MRRU model is not known}.

A critical issue in the MRRU model is that $\rho_1$ and $\rho_2$ are typically unknown in real applications.
In this paper, we use the accruing information concerning the balls in the urn to provide random thresholds which
converge a.s. to specified targets.
More specifically, our replacement matrix becomes
\begin{equation}\label{eq:replacement_matrix}
D_{n}=\left[
\begin{array}{ll}
D_{1,n}\cdot\ind_{\{Z_n\leq\hat{\rho}_{1,n}\}}\ &\ 0 \\
0\ &\ D_{2,n}\cdot\ind_{\{Z_n\geq\hat{\rho}_{2,n}\}}\\
\end{array}
\right],\end{equation}
where $\hat{\rho}_{1,n}$ and $\hat{\rho}_{2,n}$ represent the random thresholds.
We call this adaptive urn model an {\em adaptive randomly reinforced urn} (ARRU), to distinguish it
from the RRU and the MRRU.
In this paper, we investigate the asymptotic properties of the ARRU model when $m_1\neq m_2$.
Specifically, we establish the LLN for $Z_n$ and $N_{1,n}/n$,
and the CLT for $N_{1,n}/n$.
Before concluding this section,
we describe some recent works in the literature which are similar in spirit to the present work
but are quite different from our proposed model.

Let $H_n:=\bm{E}[D_n|\mathcal{F}_{n-1}]$, where $\mathcal{F}_{n-1}$ is the ``information'' up to the time $(n-1)$.
This is referred to as the generating matrix.
Asymptotic properties of the urn composition for homogeneous GPU, i.e. $H_n=H$ for all $n\geq1$, have been studied in~\cite{Athreya.Karlin.68}
under the assumption that $H$ is irreducible.
In~\cite{Smythe.et.al.96}, the extended P\'olya urn (EPU) is defined as a GPU such that all the rows of $H$ sum to the same
positive constant, i.e.
\begin{equation}\label{eq:H1c1}
H\bm{1}\ =\ c\bm{1}.
\end{equation}
Under the assumption that $H$ has simple eigenvalues, second-order asymptotic properties
on the proportion of sampled color extracted from the urn are obtained in~\cite{Smythe.et.al.96}.
In~\cite{Janson.04}, the limiting distribution of the proportion of sampled balls for homogeneous urn models are derived.
In~\cite{Bai.et.al.99}, weak consistency and asymptotic normality of the urn composition for non-homogeneous GPU are established.
However, in~\cite{Bai.et.al.99}, the sequence $\{H_n;n\geq1\}$ is deterministic and converges to a matrix $H$ satisfying~\eqref{eq:H1c1}.~\cite{Bai.et.al.02,Bai.et.al.05} extend~\cite{Bai.et.al.99} to random generating matrices
and establish almost sure convergence of the proportion of sampled balls.
They also investigate the second-order properties.
A key assumption in~\cite{Bai.et.al.02,Bai.et.al.05} is~\eqref{eq:H1c1}.
In~\cite{Cheung.et.al.06} the sequence of generating matrices is defined as function of adaptive estimators,
which guarantees the convergence of $H_n$ to a limiting matrix $H$ satisfying~\eqref{eq:H1c1}.
For ``immigrated'' urn models, theoretical results have been obtained in~\cite{Zhang.et.al.11} under the assumptions~\eqref{eq:H1c1},
or $H\bm{1}\ <\ 0$.
These extensions do not include the RRU model, where $H_n$ is diagonal, non-negative and~\eqref{eq:H1c1} is not satisfied.
For distributional results concerning large P\'olya urns, see~\cite{Chauvin.et.al.11}.
We now describe application to clinical trial literature (see~\cite{Durham.et.al.98}).
For applications to computer science, we refer the reader to~\cite{Mahmoud.08}.

\subsection{Applications to clinical trials}   \label{subsection_intro_clinical_trials}

Urn models have a long history of applications in clinical trials, by providing randomization procedures that target
certain objectives (for a review, see~\cite{Rosenberger.02}).
In this context, patients are sequentially allocated to treatments according to the sampled colors
and the associated responses are used to update the urn.
This is referred to as \emph{response-adaptive}, since the probability of assignment depends on information about the treatment performances.
For a literature review on response-adaptive designs in clinical trials see~\cite{Hu.et.al.06,Lachin.et.al.02}.
In an RRU model, responses to treatments are typically transformed by a utility function to obtain the reinforcement values,
so that the higher the reinforcement, the better the treatment.
This yields a more ethical allocation in clinical trials, because~\eqref{eq:RRU} shows that the RRU assigns more patients to the superior treatment.
However, response-adaptive designs usually aim at obtaining good inferential properties by targeting a certain proportion $\rho\in(0,1)$,
which is typically chosen to satisfy some optimality criteria (see~\cite{Rosenberger.et.al.01}).
For this reason, in~\cite{Aletti.et.al.13} the RRU was modified to asymptotically attain any target allocation proportion, $\rho\in(0,1)$.
This guarantees the MRRU design to have an asymptotic allocation within $(0,1)$
there by incorporating ethical constraints (\textit{viz.} assigning more subjects to the superior treatment).
The main issue is that $\rho_1$ and $\rho_2$ are typically functions of unknown parameters (see~\cite{Rosenberger.et.al.01}).
The ARRU model presented in this paper allows $\rho_1$ and $\rho_2$ to be  functions of such unknown parameters,
and adaptively updates by substituting sequential estimates for the parameters. The limiting results in this paper
demonstrate that such procedures target the unknown optimal allocation and provide an appropriate randomization procedure
for such trials in large samples.  We also demonstrate by simulation that the properties hold relatively well for moderate
sample sizes.

\subsection{Structure of the paper}   \label{subsection_organization_paper}

The paper is organized as follows.
In Section~\ref{section_the_model}, we present the notation and assumptions concerning the ARRU model and related main results.
Specifically, in Subsection~\ref{subsection_LLN}, we present the LLN;
in Subsection~\ref{subsection_CLT} we present the CLT under the assumption that
the thresholds are updated at exponentially changing times.
Subsection~\ref{subsection_choices_f} is devoted to the implications of the main results in the context of clinical trials.

In Section~\ref{section_preliminary_results}, we describe several fundamental results concerning the ARRU model that are needed in the proof of the CLT.
Specifically, we prove that the harmonic moments of the total number of balls in the ARRU are uniformly bounded.
Then, we use this to obtain a uniform $L_1$-bound for the distance between the urn proportion at successive update times and the adaptive thresholds.
In Section~\ref{section_proofs} the proofs of the main results are provided,
while Section~\ref{section_simulation} contains results of a simulation study.
Section~\ref{section_discussion} contains extensions to multi-color urn models.

Finally, some remarks concerning proofs are in order.
The LLN and CLT for $N_{1,n}/n$ are deduced using the asymptotic properties of $Z_n$.
For this reason, in several results of this paper we will provide a detailed probabilistic description of the sequence $\{Z_n;n\geq1\}$.

\section{Model assumptions, notation and main results}   \label{section_the_model}

We begin by describing our model precisely.
Let $\bm{\xi_1}=\left\{\xi_{1,n};n\geq1\right\}$ and $\bm{\xi_2}=\left\{\xi_{2,n};n\geq1\right\}$ be two sequences of i.i.d. random variables,
with probability distributions $\mu_1$ and $\mu_2$ respectively.
Without loss of generality (Wlog), assume that the support of $\xi_{1,n}$ and $\xi_{2,n}$ is the same.
We denote it by $S$.
Consider an urn containing $y_{1,0}>0$ red balls and $y_{2,0}>0$ white balls, and define $y_0=y_{1,0}+y_{2,0}$.
At time $n=1$, a ball is drawn at random from the urn and its color is observed.
Let the random variable $X_1$ be such that
\begin{equation*}
X_1\ =\
\begin{cases}
1\ &\text{if the extracted ball is red},\\
0\ &\text{if the extracted ball is white}.
\end{cases}\end{equation*}
We assume $X_1$ to be independent of the sequences $\bm{\xi_{1}}$ and $\bm{\xi_{2}}$.
Note that $X_1$ is a Bernoulli random variable with parameter $z_0 = y_{1,0}/y_0$.

Let $\hat{\rho}_{1,0}$ and $\hat{\rho}_{2,0}$ be two random variables such that $\hat{\rho}_{1,0},\hat{\rho}_{2,0}\in\left(0,1\right)$ and $\hat{\rho}_{1,0}\geq\hat{\rho}_{2,0}$ a.s.
Let $u:S\rightarrow \left[a, b\right]$, $0< a \leq b < \infty$.
If $X_{1}=1$ and $z_0\leq\hat{\rho}_{1,0}$, we return the extracted ball to the urn together with $D_{1,1}=u\left(\xi_{1,1}\right)$ new red balls.
While, if $X_{1}=0$ and $z_0\geq\hat{\rho}_{2,0}$, we return it to the urn together with $D_{2,1}=u\left(\xi_{2,1}\right)$ new white balls.
If $X_{1}=1$ and $z_0>\hat{\rho}_{1,0}$, or if $X_{1}=0$ and $z_0<\hat{\rho}_{2,0}$, the urn composition is not modified.
To ease notation, let denote $w_{1,0}=\ind_{\{z_0\leq\hat{\rho}_{1,0}\}}$ and $w_{2,0}=\ind_{\{z_0\geq\hat{\rho}_{2,0}\}}$.
Formally, the extracted ball is always replaced in the urn together with
$$X_{1}D_{1,1}w_{1,0}+\left(1-X_{1}\right)D_{2,1}w_{2,0}$$
new balls of the same color to the extracted one;
now, the urn composition becomes
\[\left\{
\begin{array}{l}
Y_{1,1}=y_{1,0}+X_1D_{1,1}w_{1,0}\\
\\
Y_{2,1}=y_{2,0}+(1-X_1)D_{2,1}w_{2,0}.
\end{array}
\right.\]
Set $Y_{1}=Y_{1,1}+Y_{2,1}$ and $Z_{1}=Y_{1,1}/Y_1$.
Now, by iterating the above procedure we define
$\hat{\rho}_{1,1}$ and $\hat{\rho}_{2,1}$ to be two random variables, measurable with respect to the $\sigma$-algebra $\mathcal{F}_1=\sigma\left(X_1,X_1\xi_{1,1}+(1-X_1)\xi_{2,1}\right)$, with $\hat{\rho}_{1,1},\hat{\rho}_{2,1}\in\left(0,1\right)$ and $\hat{\rho}_{1,1}\geq\hat{\rho}_{2,1}$ a.s.
Let $m_1 =\int u\left(y\right)\mu_1\left(dy\right)$ and $m_2 =\int u\left(y\right)\mu_2\left(dy\right)$ be the means of
$\{D_{1,n};n\geq1\}$ and $\{D_{2,n};n\geq1\}$, respectively.
We assume throughout the paper the following condition:
\begin{equation}\label{ass:different_mean}
m_1\neq m_2.
\end{equation}
The urn process is then repeated for all $n\geq1$.
Let $\hat{\rho}_{1,n}$ and $\hat{\rho}_{2,n}$ be two random variables, measurable with respect to the $\sigma$-algebra $$\mathcal{F}_n=\sigma\left(X_1,X_1\xi_{1,1}+\left(1-X_1\right)\xi_{2,1},...,X_n,X_1\xi_{1,n}+\left(1-X_n\right)\xi_{2,n}\right),$$
with $\hat{\rho}_{1,n},\hat{\rho}_{2,n}\in\left(0,1\right)$ and $\hat{\rho}_{1,n}\geq\hat{\rho}_{2,n}$ a.s.
We will refer to $\hat{\rho}_{j,n}$ $j=1,2$ as threshold parameters.
At time $n+1$, a ball is extracted and let $X_{n+1}=1$ if the ball is red and $X_{n+1}=0$ otherwise.
Then, the ball is returned to the urn together with
$$X_{n+1}D_{1,n+1}W_{1,n}+\left(1-X_{n+1}\right)D_{2,n+1}W_{2,n}$$
balls of the same color, where $D_{1,n+1}=u\left(\xi_{1,n+1}\right)$, $D_{2,n+1}=u\left(\xi_{2,n+1}\right)$, $W_{1,n}=\ind_{\{Z_n\leq\hat{\rho}_{1,n}\}}$,
$W_{2,n}=\ind_{\{Z_n\geq\hat{\rho}_{2,n}\}}$ and $Z_{n+1}=Y_{1,{n}}/Y_{n}$ for any $n\geq1$.
Formally,
\[\left\{
\begin{array}{l}
Y_{1,n+1}=y_{1,0}+\sum_{i=1}^{n+1}X_iD_{1,i}W_{1,i-1}\\
\\
Y_{2,n+1}=y_{2,0}+\sum_{i=1}^{n+1}\left(1-X_i\right)D_{2,i}W_{2,i-1}
\end{array}
\right.\]
and $Y_{n+1}=Y_{1,n+1}+Y_{2,n+1}$.
If $X_{n+1}=1$ and $Z_n>\hat{\rho}_{1,n}$, i.e. $W_{1,n}=0$, or if $X_{n+1}=0$ and $Z_n<\hat{\rho}_{2,n}$, i.e. $W_{2,n}=0$,
the urn composition does not change at time $n+1$.
Note that condition $\hat{\rho}_{1,n}\geq\hat{\rho}_{2,n}$ a.s., which implies $W_{1,n}+W_{2,n}\geq1$, ensures that the urn composition can change with positive probability for any $n\geq1$, since the replacement matrix~\eqref{eq:replacement_matrix} is never a zero matrix.
Since, conditionally to the $\sigma$-algebra $\mathcal{F}_n$, $X_{n+1}$ is assumed to be independent of $\bm{\xi_{1}},\bm{\xi_{2}}$,
$X_{n+1}$ is conditionally Bernoulli distributed with parameter $Z_n$.

We will denote by $N_{1,n}$ and $N_{2,n}$ the number of red and white sampled balls, respectively, after the first $n$ draws,
that is $N_{1,n}=\sum_{i=1}^nX_i$ and $N_{2,n}=\sum_{i=1}^n\left(1-X_i\right)$.
Let $\rho_1$ and $\rho_2$ be two constants such that $0<\rho_2\leq\rho_1<1$.
We will adopt the following notation:
\[\begin{aligned}\hat{\rho}_{n}&&:=&\hat{\rho}_{1,n}\ind_{\{m_1>m_2\}}\ +\ \hat{\rho}_{2,n}\ind_{\{m_1<m_2\}};\\
\rho&&:=&\rho_1\ind_{\{m_1>m_2\}}\ +\ \rho_2\ind_{\{m_1<m_2\}}.\end{aligned}\]

\subsection{Law of large numbers}   \label{subsection_LLN}

Our first result is concerned with the LLN.
\begin{thm}\label{thm:convergence_almost_sure}
Under assumptions~\eqref{ass:different_mean} and
\begin{equation}\label{eq:rho_convergence_almost_sure}
\lim_{n\rightarrow\infty}\hat{\rho}_n\ =\ \rho\ \ \ a.s.
\end{equation}
we have that
\begin{equation}\label{eq:Z_convergence_almost_sure}
\lim_{n\rightarrow\infty}Z_n\ =\ \rho\ \ \ a.s.
\end{equation}
\end{thm}
From Theorem~\ref{thm:convergence_almost_sure} we can obtain the convergence of sampled balls, namely $N_{1,n}/n$.
\begin{cor}\label{cor:convergence_almost_sure_NR}
Under assumptions~\eqref{ass:different_mean} and~\eqref{eq:rho_convergence_almost_sure},
\begin{equation}\label{eq:NRn_convergence_almost_sure}
\lim_{n\rightarrow\infty}\frac{N_{1,n}}{n}\ =\ \rho\ \ \ a.s.
\end{equation}
\end{cor}

\subsection{Central limit theorem}   \label{subsection_CLT}

We next study the limit distribution of proportion of sampled balls $\frac{N_{1,n}}{n}$.
By the description of the model, $\frac{N_{1,n}}{n}$ depends on the sequence $\hat{\rho}_{j,n}$, $j=1,2$.
However, frequent changes to $\hat{\rho}_{j,n}$ make the sequence $\frac{N_{1,n}}{n}$ more erratic.
To stabilize the behavior of $\left\{\frac{N_{1,n}}{n};n\geq1\right\}$, we fix a constant $q>1$ and introduce the sequence $\tilde{\rho}_{j,n}$ $j=1,2$ as
\begin{equation}\label{def:adaptive_thresholds}
\tilde{\rho}_{j,n}:=\hat{\rho}_{j,\left[q^i\right]},\ \ \ \ as\ \left[q^i\right]\leq n< \left[q^{i+1}\right],
\end{equation}
for $i\in\mathbb{N}$; that is, we adapt the threshold parameters to change ``slowly'' at exponential times $\left\{\left[q^{i}\right], i=1,2,..\right\}$.
An alternative definition of $\tilde{\rho}_{j,n}$ $j=1,2$, which is used in some proofs, is the following
\begin{equation}\label{def:k_n}
\left(\tilde{\rho}_{1,n},\tilde{\rho}_{2,n}\right)\ :=
\ \left(\hat{\rho}_{1,\left[q^{k_n}\right]},\hat{\rho}_{2,\left[q^{k_n}\right]}\right),\ \ \
k_n\ :=\ [log_{q}(n)],
\end{equation}
for any $n\geq1$.
We will denote by
$$\tilde{\rho}_{n}=\tilde{\rho}_{1,n}\ind_{\{m_1>m_2\}}\ +\ \tilde{\rho}_{2,n}\ind_{\{m_1<m_2\}}.$$
We now turn to the statement of the CLT.
In the following $\stackrel{d}{\rightarrow}$ represents the convergence in distribution.
\begin{thm}\label{thm:convergence_in_distribution}
Let $\tilde{\rho}_{1,n}$ and $\tilde{\rho}_{2,n}$ be as in~\eqref{def:adaptive_thresholds}.
Assume that for any $\epsilon>0$ and $j=1,2$, there exists $0<c_1<\infty$ such that
\begin{equation}\label{eq:rho_exponentially_fast}
\bm{P}\left(|\hat{\rho}_{j,n}-\rho_{j}|>\epsilon\right)\ \leq\  c_1\exp\left(-n\epsilon^2\right),
\end{equation}
for large $n$. Then, under assumption~\eqref{ass:different_mean}, we have that
\begin{equation}\label{eq:X_eta_convergence_distribution}
\sqrt{n}\left(\frac{N_{1,n}}{n}-\bar{\rho}_n\right)\ \stackrel{d}{\rightarrow}\ \mathcal{N}\left(0,\rho\left(1-\rho\right)\right).
\end{equation}
where $\bar{\rho}_n=\frac{\sum_{i=1}^{n}\tilde{\rho}_{i-1}}{n}$.
\end{thm}

\begin{rem}\label{rem:condition_easy_to_satisfy_2}
The result of Theorem~\ref{thm:convergence_in_distribution} continues to hold if~\eqref{eq:rho_exponentially_fast} is not satisfied,
but~\eqref{eq:rho_convergence_almost_sure} and the following conditions hold:
\begin{itemize}
\item[(c1)] $\limsup_{n\rightarrow\infty}\sqrt{n}\bm{E}\left[|\hat{\rho}_{n}-\rho|\right]\ <\ \infty.$
\item[(c2)] There exists $\epsilon\in\left(0,1/2\right)$ such that $\hat{\rho}_{j,n}\in[\epsilon,1-\epsilon]$ a.s. for any $n\geq1$, $j=1,2$.
\end{itemize}
\end{rem}

Theorem~\ref{thm:convergence_in_distribution} introduces an asymptotic bias for $N_{1,n}/n$ given by $(\bar{\rho}_n-\rho)$.
We show that this bias is exactly of order $O(n^{-1/2})$;
our next proposition makes this observation precise.
\begin{prop}\label{prop:rhomean_closed_rho}
Let $\tilde{\rho}_{1,n}$ and $\tilde{\rho}_{2,n}$ be as in~\eqref{def:adaptive_thresholds}.
Then, under assumption~\eqref{ass:different_mean} and~\eqref{eq:rho_exponentially_fast},
\begin{equation}\label{ass:closeness_barrhon_rho}
\limsup_{n\rightarrow\infty}\ n\cdot\bm{E}\left[|\bar{\rho}_n-\rho|^2\right]\ <\ \infty.
\end{equation}
\end{prop}

\begin{rem}\label{rem:condition_easy_to_satisfy_3}
The result of Proposition~\ref{prop:rhomean_closed_rho} continues to hold if~\eqref{eq:rho_exponentially_fast} is not satisfied,
but the following condition holds:
\begin{itemize}
\item[(c4)] $\limsup_{n\rightarrow\infty}n\bm{E}\left[|\hat{\rho}_{n}-\rho|^2\right]\ <\ \infty$.
\end{itemize}
\end{rem}

In the case when $\hat{\rho}_{1,n}=\rho_1$ and $\hat{\rho}_{2,n}=\rho_2$ for any $n\geq 0$,
Theorem~\ref{thm:convergence_in_distribution} provides a CLT for the allocation proportion of MRRU model.
This is summarized in the following corollary:
\begin{cor}\label{cor:MRRU_CLT}
In a MRRU, under assumption~\eqref{ass:different_mean}, we have that
\begin{equation*}
\sqrt{n}\left(\frac{N_{1,n}}{n}-\rho\right)\ \stackrel{d}{\rightarrow}\ \mathcal{N}\left(0,\rho\left(1-\rho\right)\right).
\end{equation*}
\end{cor}

\subsection{Application to clinical trials (revisited)}   \label{subsection_choices_f}

Consider two competing treatments $\mathcal{T}_1$ and $\mathcal{T}_2$.
The random variables $\xi_{1,n}$ and $\xi_{2,n}$ are interpreted as the potential responses to treatments $\mathcal{T}_1$ and $\mathcal{T}_2$, respectively,
given by subjects that sequentially enter the trial.
At all times $n\geq1$, a subject is allocated to a treatment according to the color of the sampled ball and a new response is collected.
Note that only one response is observable from every subject, that is $X_n\xi_{1,n}+ \left(1-X_n\right)\xi_{2,n}$.
The function $u$ transforms the responses into reinforcements $D_{1,n}$ and $D_{2,n}$ that update the urn.
Typically, $u$ is chosen such that $\mathcal{T}_1$ (or $\mathcal{T}_2$) is considered the superior treatment when $m_1>m_2$ ($m_1<m_2$).
We assume there exists a unique superior treatment, which is formally stated in assumption~\eqref{ass:different_mean}.

We now describe the role of the sequences $\{\hat{\rho}_{1,n};n\geq1\}$ and $\{\hat{\rho}_{2,n};n\geq1\}$ in clinical trails.
Assume the distributions $\mu_1$ and $\mu_2$ are parametric,
depending on the vectors $\bm{\theta}_1$ and $\bm{\theta}_2$ respectively,
with $\bm{\theta}=\left(\bm{\theta_1},\bm{\theta_2}\right)\in\Theta\subset R^d$, with $d\geq1$.
Let $\bm{\hat{\theta}}_n=\left(\bm{\hat{\theta}}_{1,n},\bm{\hat{\theta}}_{2,n}\right)$ be an estimator of $\bm{\theta}$ after the first $n$ allocations,
so that $\bm{\hat{\theta}}_n$ is measurable with respect to the $\sigma$-algebra $\mathcal{F}_n$.
We assume that the distributions $\mu_1$ and $\mu_2$ are parametrically independent,
in the sense that $\mu_1$ does not depend on $\bm{\theta}_2$ and $\mu_2$ does not depend on $\bm{\theta}_1$.
Hence, $\bm{\hat{\theta}}_{1,n}$ is computed with the $N_{1,n}$ observations $\{\xi_{1,i}:X_i=1,i\leq n\}$,
while $\bm{\hat{\theta}}_{2,n}$ is computed with the $N_{2,n}$ observations $\{\xi_{2,i}:X_i=0,i\leq n\}$.
Thus, $\{\hat{\rho}_{1,n};n\geq1\}$ and $\{\hat{\rho}_{2,n};n\geq1\}$ are defined as follows:
\begin{equation}\label{eq:f_clinical_trials}
\hat{\rho}_{1,n}:=f_{1}\left(\bm{\hat{\theta}}_{1,n}\right)\ \ \ \ \ \textit{and}\ \ \ \ \ \hat{\rho}_{2,n}:=f_{2}\left(\bm{\hat{\theta}}_{2,n}\right), \ \ \forall n\geq1,
\end{equation}
where $f_{1}:\Theta\rightarrow\left(0,1\right)$ and $f_{2}:\Theta\rightarrow\left(0,1\right)$ are two continuous functions such that
$$f_{1}\left(\bm{x}\right)\geq f_{2}\left(\bm{x}\right),\ \ \ \forall \bm{x}\in\Theta;$$
this implies $\hat{\rho}_{1,n}\geq\hat{\rho}_{2,n}$ a.s. for every $n\geq1$.
Moreover, set
\begin{equation*}
\rho_1:=f_{1}\left(\bm{\theta}\right)\ \ \ \ \ \ \ \ \textit{and}\ \ \ \ \ \ \ \ \rho_2:=f_{2}\left(\bm{\theta}\right).
\end{equation*}

The LLN presented in Theorem~\ref{thm:convergence_almost_sure}
suggests a direct interpretation for the functions $f_1$ and $f_2$ in a clinical trial context:
$f_1\left(\bm{\theta}\right)$ and $f_2\left(\bm{\theta}\right)$
represent the desired limiting allocations for the sequence $N_{1,n}/n$,
in case the superior treatment is $\mathcal{T}_1$ ($m_1>m_2$) or $\mathcal{T}_2$ ($m_1<m_2$), respectively.
This is a great improvement, since the design can target an arbitrary known function of all the parameters of the response distributions.

Ideally, $f_1$ and $f_2$ are chosen to obtain good statistical properties from the design.
Typically, in clinical trials, a design is constructed to satisfy certain optimality criteria related to its statistical performances (e.g., power; see~\cite{Rosenberger.et.al.01}).
Letting $\eta\left(\bm{\theta}\right)$ denote the limit proportion of subjects to be allocated to treatment $\mathcal{T}_1$,
this design can be obtained by the urn model described in Section~\ref{section_the_model} by choosing $f_1\left(\bm{\theta}\right)=f_2\left(\bm{\theta}\right)=\eta\left(\bm{\theta}\right)$.
However, in some experiments, ethical aspects are important and the main goal may be to assign fewer subjects to the inferior treatment;
in this case we choose $f_1\left(\bm{\theta}\right)\simeq1$ and $f_2\left(\bm{\theta}\right)\simeq0$.
Designs requiring both ethical and statistical goals can also be obtained from our design, by setting $f_1\left(\bm{\theta}\right)\geq\eta\left(\bm{\theta}\right)\geq f_2\left(\bm{\theta}\right)$.
For instance, we may take
\begin{equation}\label{eq:linear_combination}
f_{1}\left(\bm{\theta}\right)=p\cdot \eta\left(\bm{\theta}\right)+\left(1-p\right)\cdot 1,\ \ \ f_{2}\left(\bm{\theta}\right)=p\cdot \eta\left(\bm{\theta}\right)+\left(1-p\right)\cdot 0,\ \ \ p\in\left(0,1\right],
\end{equation}
where $p$ is a biasing term, which introduces a trade-off between the ethics and statistical properties.

Finally, it is worth emphasizing that conditions~\eqref{eq:rho_convergence_almost_sure} and~\eqref{eq:rho_exponentially_fast}
required in the LLN of Theorem~\ref{thm:convergence_almost_sure} and in the CLT of Theorem~\ref{thm:convergence_in_distribution}, respectively,
are straightforwardly satisfied when we take $\bm{\hat{\theta}}_{n}$ to be maximum likelihood estimators (MLEs) for $\bm{\theta}$.

Moreover, condition (c2) in Remark~\ref{rem:condition_easy_to_satisfy_2} is equivalent of the assumption that
the ranges of $f_1$ and $f_2$ are subsets of $\left[\epsilon,1-\epsilon\right]$, for some $\epsilon\in\left(0,1/2\right)$.

\section{Harmonic moments and related asymptotics}   \label{section_preliminary_results}

\subsection{Harmonic moments}   \label{subsection_harmonic_moments}

In this subsection, we show that the harmonic moments of the total number of balls in the urn are uniformly bounded.
This is a key result which is needed in several probabilistic estimates, and in particular in the proof of the CLT.
More specifically, as explained previously the results concerning the asymptotic behavior of $N_{1,n}/n$,
depend critically on the behavior of $(Z_n-\hat{\rho}_n)$.
In Subsection~\ref{subsection_asymptotic_difference} we provide bounds for $Y_n(Z_n-\hat{\rho}_n)$, by using comparison arguments with the MRRU model.
Now, to replace the random scaling $Y_n$ by the deterministic scaling $n$, one needs to investigate the behavior of ${n}/{Y_n}$.
Our next theorem provides a precise estimates of the $j^{th}$ moment of ${n}/{Y_n}$ for any $j\geq0$.


\begin{thm}\label{thm:Y_bounded_and_diverge}
Under assumption~\eqref{ass:different_mean} and~\eqref{eq:rho_exponentially_fast},
for any $j>0$, we have that
$$\sup_{n\rightarrow\infty}\ \bm{E}\left[\left(\frac{n}{Y_n}\right)^j\right]\ <\ \infty.$$
\end{thm}

In the proof of Theorem~\ref{thm:Y_bounded_and_diverge},
we need the following lemma that provides an upper bound on the increments of the urn process $Z_n$,
by imposing a condition on the total number of balls in the urn $Y_n$.
Hence, the proof of Theorem~\ref{thm:Y_bounded_and_diverge} is reported after the following result.
\begin{lem}\label{lem:Y_increments}
For any $\epsilon\in\left(0,1\right)$, we have that
\begin{equation}\label{eq:how_big_is_D}
\left\{\ Y_n > b \left(\frac{1-\epsilon}{\epsilon}\right)\ \right\}\ \ \ \subseteq\ \ \
\left\{\ |Z_{n+1}-Z_n|<\epsilon\ \right\}.
\end{equation}
\end{lem}

\proof
The difference $(Z_{n+1}-Z_{n})$ can be expresses as follows:
$$\frac{Y_{1,n}+X_{n+1}W_{1,n}D_{1,n+1}}
{Y_{n}+X_{n+1}W_{1,n}D_{1,n+1}+(1-X_{n+1})W_{2,n}D_{2,n+1}}-\frac{Y_{1,n}}{Y_{n}}$$
Consider $\{Z_{n+1}> Z_n\}$, since the case $\{Z_{n+1}< Z_n\}$ is analogous.
Note that $\{Z_{n+1}> Z_n\}$ implies that $\{X_{n+1}=1\}$ and $\{W_{1,n}=1\}$.
Then, since $D_{1,n+1}<b$ a.s., on the set $\{Z_{n+1}> Z_n\}$ we have
\[\begin{aligned}
Z_{n+1}-Z_{n}\ &&\leq&\ \frac{Y_{1,n}+D_{1,n+1}}{Y_{n}+D_{1,n+1}}-\frac{Y_{1,n}}{Y_{n}}\\
&&=&\ \frac{D_{1,n+1}}{D_{1,n+1}+Y_{n}}\left(1-Z_{n}\right)\ \leq\ \frac{b}{b+Y_{n}}\ <\ \epsilon,
\end{aligned}\]
where the last inequality follows from $\{Y_n > b(1-\epsilon)/\epsilon\}$ in~\eqref{eq:how_big_is_D}.
\endproof

\proof[Proof of Theorem~\ref{thm:Y_bounded_and_diverge}]
In this proof, when we have set of integers $\{[a_1],..,[b_1]\}$ with $a_1,b_1\notin \mathbb{N}$,
to ease notation we will just write $\{a_1,..,b_1\}$, omitting the symbol $[\cdot]$.
First, note that, since $D_{1,i},D_{2,i}\geq a$ a.s. for any $i\geq1$ and $Y_0>0$, we have that
\begin{equation}\label{eq:intermidiate_2}\begin{aligned}
Y_{n}\ &=&&\ Y_0\ +\ \sum_{i=1}^n\left(D_{1,i}X_iW_{1,i-1}+D_{2,i}\left(1-X_i\right)W_{2,i-1}\right)\\
&\geq&&\ Y_0\ +\ a\cdot \sum_{i=1}^n\left(X_iW_{1,i-1}+\left(1-X_i\right)W_{2,i-1}\right),\\
&\geq&&\ Y_0\ +\ a\cdot \sum_{i=n\beta}^n\left(X_iW_{1,i-1}+\left(1-X_i\right)W_{2,i-1}\right),
\end{aligned}\end{equation}
for any $\beta\in(0,1)$. To keep calculation transparent we choose $\beta=1/2$.
We recall that, by construction, we have that
$W_{1,i-1},W_{2,i-1}\in\{0;1\}$ and $W_{1,i-1}+W_{2,i-1}\geq 1$ for any $i\geq1$;
hence, the random variables $X_iW_{1,i-1}+\left(1-X_i\right)W_{2,i-1}$
are, conditionally to the $\sigma$-algebra $\mathcal{F}_{i-1}$, Bernoulli distributed with parameter greater than or equal to
$\min\{Z_{i-1};1-Z_{i-1}\}$.
Hence, the behavior of $Y_n$ is intrinsically related to the behavior of $Z_n$.

\noindent Thus, let us introduce the sets $A_{d,n}$ (down), $A_{c,n}$ (center) and $A_{u,n}$ (up) as follows:
\[\begin{aligned}
A_{d,n}\ &&:=&\ \left\{\bigcup_{n/2\leq i\leq n}\left\{Z_i<c\right\}\right\},\\
A_{c,n}\ &&:=&\ \left\{\bigcap_{n/2\leq i\leq n}\left\{Z_i\in\left[c,1-c\right]\right\}\right\},\\
A_{u,n}\ &&:=&\ \left\{\bigcup_{n/2\leq i\leq n}\left\{Z_i>1-c\right\}\right\},
\end{aligned}\]
where $c\in(0,1)$ will be appropriately fixed more ahead in the proof.
Then, we perform the following decomposition on the behavior of $\{Z_i;n/2\leq i\leq n\}$,
$$\bm{E}\left[\left(\frac{n}{Y_n}\right)^j\right]\ \leq\ \left(\frac{n}{Y_0}\right)^j\cdot \bm{P}(A_{d,n})\ +\
\bm{E}\left[\left(\frac{n}{Y_n}\right)^j\ind_{A_{c,n}}\right]\ +\ \left(\frac{n}{Y_0}\right)^j\cdot \bm{P}(A_{u,n}).$$
On the set $A_{c,n}$ the process $\{Z_i;n/2\leq i\leq n\}$ is bounded away from the extreme values $\{0;1\}$;
Hence we can use comparison arguments with a sequence of i.i.d. Bernoulli random variables with parameter $c$ to get the boundedness of $\bm{E}\left[\left(n/Y_n\right)^j\ind_{A_{c,n}}\right]$.
After that, we will focus on proving that $\bm{P}(A_{d,n})$ and $\bm{P}(A_{u,n})$ converge to zero exponentially fast.

First, note that on the set $A_{c,n}$ the random variables
$$X_iW_{1,i-1}+\left(1-X_i\right)W_{2,i-1}$$
are, conditionally to the $\sigma$-algebra $\mathcal{F}_{i-1}$, Bernoulli with parameter with parameter greater than or equal to $c$ for any $i=n/2,..,n$.
Hence, if we introduce $\{B_i;i\geq1\}$ a sequence of i.i.d. Bernoulli random variable with parameter $c$,
from~\eqref{eq:intermidiate_2} we have that
$$\bm{E}\left[\left(\frac{n}{Y_{n}}\right)^j\ind_{A_{c,n}}\right]\ \leq\ \frac{1}{a^j}\bm{E}\left[\left(\frac{n}{Y_0/a+\sum_{i=n/2}^nB_i}\right)^j\right].$$
We now show that
$$\limsup_{n\rightarrow\infty}\ \bm{E}\left[\left(\frac{n}{K_0+\sum_{i=1}^nB_i}\right)^j\right]\ <\ \infty,$$
with $K_0=Y_0/a$.
To this end, we apply Theorem 2.1 of~\cite{Etemadi.et.al.97}, with $n_0=1$, $p=j$, $Z_{i,n}=B_i+Y_0/n$ for $i\leq n$.
All the assumptions of the theorem are satisfied in our case.
In fact, at first we have $\bm{E}\left[\bar{Z}_{n_0}^{-p}\right]<\infty$ because $\bm{E}\left[\left(Y_0+B_1\right)^{-j}\right]\leq K_0^{-j}<\infty$.\\
Secondly, note that $Z_{i,n}$ are identically distributed for all $i\leq n$, since $B_i$ are i.i.d. Bernoulli of parameter $c$.
Finally, $\bar{Z}_{n}$ converges in distribution, since $\bar{Z}_{n}=\sum_{i=1}^nB_i/n+K_0\stackrel{a.s.}{\rightarrow}c+K_0$.
Hence, by Theorem 2.1 of~\cite{Etemadi.et.al.97}, it follows that $\bm{E}\left[\bar{Z}_{n}^{-p}\right]$ is uniformly integrable. As a consequence,
$$\limsup_{n\rightarrow\infty}\ \bm{E}\left[\left(\frac{n}{K_0+\sum_{i=1}^nB_i}\right)^j\right]\ =
\ \limsup_{n\rightarrow\infty}\ \bm{E}\left[\bar{Z}_{n}^{-p}\right]\ <\ \infty.$$

Now, we will prove that $\bm{P}(A_{d,n})$ and $\bm{P}(A_{u,n})$ converge to zero exponentially fast.
We will show that this occurs because $\hat{\rho}_{1,n}$ and $\hat{\rho}_{2,n}$
are bounded away from the extreme values $\{0;1\}$, with probability that converge to one exponentially fast.
Formally, fix $\epsilon>0$, such that $\rho_1+\epsilon<1$ and $\rho_2-\epsilon>0$,
and define $\alpha_n:=n^{\alpha}$, $\alpha\in(0,1)$, for any $n\geq1$.
Now, for any $n\geq 1$ define the following sets:
\[\begin{aligned}
A_{1,n}\ &&:=&\ \left\{\ \sup_{i\geq \alpha_n}\{\hat{\rho}_{1,i}\}>\rho_1+\epsilon\ \right\},\\
A_{2,n}\ &&:=&\ \left\{\ \inf_{i\geq \alpha_n}\{\hat{\rho}_{2,i}\}<\rho_2-\epsilon\ \right\},\\
A_{3,n}\ &&:=&\ \left\{\ \inf_{i\geq \alpha_n}\{\min\{1-\hat{\rho}_{1,i};\hat{\rho}_{2,i}\}\}
\geq\min\{1-\rho_1;\rho_2\}-\epsilon\ \right\},\\
\end{aligned}\]
where we recall that $\hat{\rho}_{1,i}$ and $\hat{\rho}_{2,i}$ are the adaptive thresholds.
Note that $A_{1,n}\cup A_{2,n}\cup A_{3,n}=\Omega$.
We have that
\[\begin{aligned}
\bm{P}(A_{d,n})\ &&\leq&\ \bm{P}\left(A_{1,n}\right)\ +\ \bm{P}\left(A_{2,n}\right)\ +\ \bm{P}\left(A_{3,n}\cap A_{d,n}\right),\\
\bm{P}(A_{u,n})\ &&\leq&\ \bm{P}\left(A_{1,n}\right)\ +\ \bm{P}\left(A_{2,n}\right)\ +\ \bm{P}\left(A_{3,n}\cap A_{u,n}\right).
\end{aligned}\]

First, we prove that $\bm{P}\left(A_{1,n}\right)$ and $\bm{P}\left(A_{2,n}\right)$ converge to zero exponentially fast.
Consider the term $\bm{P}\left(A_{1,n}\right)$.
From the definition of $A_{1,n}$, we obtain
$$\bm{P}\left(A_{1,n}\right)= \bm{P}\left(\bigcup_{i\geq \alpha_n}\left\{\hat{\rho}_{1,i}>\rho_1+\epsilon\right\}\right)\leq\sum_{i\geq \alpha_n}\bm{P}\left(\hat{\rho}_{1,i}>\rho_1+\epsilon\right).$$
From~\eqref{eq:rho_exponentially_fast}, for large $i$ we have that
$$\bm{P}\left(\hat{\rho}_{1,i}>\rho_1+\epsilon\right)\leq c_1\exp\left(-i\epsilon^2\right),$$
with $0<c_1<\infty$.
Hence, using the fact that $Y_n$ is increasing we have that
\[\begin{aligned}
\bm{P}\left(A_{1,n}\right)\ &&\leq&\ \sum_{i\geq \alpha_n}\bm{P}\left(\hat{\rho}_{1,i}>\rho_1+\epsilon\right)\\
&&\leq&\ c_1\sum_{i\geq \alpha_n}\exp\left(-i\epsilon^2\right)\\
&&=&\ c_1\exp\left(-\alpha_n\epsilon^2\right).
\end{aligned}\]
Similar arguments can be applied to prove $\bm{P}\left(A_{2,n}\right)\rightarrow0$ exponentially fast.

Finally, we show that $\bm{P}\left(A_{3,n}\cap A_{d,n}\right)$ and $\bm{P}\left(A_{3,n}\cap A_{u,n}\right)$ converge to zero exponentially fast.
Consider $\bm{P}\left(A_{3,n}\cap A_{d,n}\right)$, since the proof for\\
$\bm{P}\left(A_{3,n}\cap A_{u,n}\right)$ is analogous.
First, let introduce $\phi:=\min\{\rho_{2};1-\rho_{1}\}$, and rewrite $A_{3,n}$ as follows:
$$A_{3,n} =\ \left\{ \inf_{i\geq \alpha_n}\{\hat{\rho}_{2,n};1-\hat{\rho}_{1,n}\}\geq\phi-\epsilon\ \right\}.$$
Define the set $\tilde{A}_{d,n}$ as follows:
$$\tilde{A}_{d,n} :=\ \left\{\bigcap_{\alpha_n\leq i\leq n/2}\left\{Z_i<c\right\}\right\}.$$
We now set an appropriate value of $c$ such that
\begin{equation}\label{eq:subset}
\left\{A_{3,n}\cap A_{d,n}\right\}\ \subset\ \left\{A_{3,n}\cap \tilde{A}_{d,n}\right\},
\end{equation}
for any $n\geq1$.
To do that, we need to set $c$ such that $\{Z_{i}\geq c\}\subset\{Z_{i+1}\geq c\}$ for any $i\geq \alpha_n$.
First, note that on the set $A_{3,n}$, $\{\hat{\rho}_{2,i}\geq(\phi-\epsilon)\}$ for any $i\geq\alpha_n$.
Hence, for any $c< (\phi-\epsilon)$, if $\{c\leq Z_i \leq (\phi-\epsilon)\}$ we have $W_{2,i}=0$, that implies $Z_{i+1}\geq Z_i$ and so $Z_{i+1}\geq c$.
Alternatively, if $\{Z_i \geq (\phi-\epsilon)>c\}$, the set $\{Z_{i+1}\leq Z_i\}$ is possible,
and hence we have to bound the increments of $Z_n$ to guarantee that $Z_{i+1}\geq c$,
i.e. set $c$ such that
$$|Z_{i+1}-Z_{i}|<(\phi-\epsilon)-c,\ \ \forall i\geq0.$$
Using~\eqref{eq:how_big_is_D}, we obtain
\begin{equation}\label{def:p_0}
c\leq p_0:=\frac{Y_0}{Y_0+b}\cdot\left(\phi-\epsilon\right).
\end{equation}
This guarantees~\eqref{eq:subset} holds for any $n\geq1$.

We now show that $\bm{P}(A_{3,n}\cap \tilde{A}_{d,n})$ converges to zero exponentially fast.
To this end, first note that on the set $A_{3,n}$, we have $\hat{\rho}_{2,i}>\rho_2-\epsilon$ for any $i=\alpha_n,..,n/2$;
moreover, on the set $\tilde{A}_{d,n}$, we have $Z_i<p_0$ for any $i=\alpha_n,..,n/2$.
These considerations imply that $W_{2,i}=0$ and $W_{1,i}=1$ for any $i=\alpha_n,..,n/2$, on the set $A_{3,n}\cap \tilde{A}_{d,n}$.
Hence, we can write
\begin{equation}\label{eq:previous_calculations}
Z_{n/2}=\frac{Y_{1,\alpha_n}+\sum_{i=\alpha_n}^{n/2}X_iD_{1,i}}{Y_{\alpha_n}+\sum_{i=\alpha_n}^{n/2}X_iD_{1,i}}
\geq \frac{y_{1,0}+a\sum_{i=\alpha_n}^{n/2}X_i}{(y_0+\alpha_nb)+a\sum_{i=\alpha_n}^{n/2}X_i},\end{equation}
where the inequality is because $Y_{1,\alpha_n}\geq y_{1,0}$, $Y_{\alpha_n}\leq y_0+\alpha_nb$ and $D_{1,i}\geq a$ a.s. for any $i\geq1$.
Now, define for any $n\geq1$ the set $A_{4,n}$ as follows:
$$A_{4,n}\ :=\ \left\{\ \sum_{i=\alpha_n}^{n/2}X_i\ >\ \frac{p_0}{a(1-p_0)}\left(y_0+\alpha_nb\right)\ \right\},$$
and consider the set $A_{\phi,n}\cap \tilde{A}_{d,n}\cap A_{4,n}$.
On the set $A_{\phi,n}\cap \tilde{A}_{d,n}$ we can use the definition of $A_{4,n}$ in~\eqref{eq:previous_calculations},
obtaining
$$\left\{A_{3,n}\cap \tilde{A}_{d,n}\cap A_{4,n}\right\}\subset\left\{\left\{Z_{n/2}>p_0\right\}\cap \tilde{A}_{d,n}\right\}.$$
However, $\left\{Z_{n/2}>p_0\right\}\cap \tilde{A}_{d,n}=\emptyset$. Hence, $\bm{P}(A_{\phi,n}\cap \tilde{A}_{d,n}\cap A_{4,n})=0$ and
it is sufficient to show that $\bm{P}\left(A_{3,n}\cap \tilde{A}_{d,n}\cap A_{4,n}^C\right)$ converges to zero exponentially fast.

To this end, note that on the set $A_{3,n}\cap \tilde{A}_{d,n}$ we have $Z_{i+1}\geq Z_i$ for any $i=\alpha_n,..,n/2$,
since we previously showed that $W_{2,i}=0$ and $W_{1,i}=1$.
Hence, on the set $A_{3,n}\cap \tilde{A}_{d,n}$, $\{X_i,i=\alpha_n,..,n/2\}$ are conditionally Bernoulli with parameter
$p_i\geq Z_{\alpha_n}$ a.s.
Now, let denote with $\{\varrho_{i,n};i=1,..,n/2-\alpha_n\}$ a sequence of i.i.d. Bernoulli random variable with parameter $z_{0,n}$, defined as
$$z_{0,n}\ :=\ \frac{y_{1,0}}{y_0+\alpha_nb}\ \leq\ Z_{\alpha_n}\ \ \ a.s.;$$
it follows that $\bm{P}\left(A_{3,n}\cap\tilde{A}_{d,n}\cap A_{4,n}^C\right)$ is less than or equal than
\begin{equation}\label{eq:bernoulli_chernoff}
\bm{P}\left(\sum_{i=1}^{n/2-\alpha_n}\varrho_{i,n}\ \leq\ \frac{p_0}{a(1-p_0)}\left(y_0+\alpha_nb\right)\right).
\end{equation}
Finally, we use the following Chernoff's upper bound for i.i.d. random variables in $[0,1]$ (see~\cite{Dembo.et.al.98})
\begin{equation}\label{eq:bernoulli_chernoff_result}
\bm{P}\left( S_n\leq c_0\cdot\bm{E}[S_n]\right)\ \leq\ \exp\left(-\frac{(1-c_0)^2}{2}\cdot\bm{E}[S_n]\right),
\end{equation}
with $c_0\in (0,1)$ and $S_n=\sum_i^nX_i$.
In our case, we have that~\eqref{eq:bernoulli_chernoff} can be written as $\bm{P}\left( S_n\leq c_n\cdot\bm{E}[S_n]\right)$,
where $S_n=\sum_{i=1}^{n/2-\alpha_n}\varrho_{i,n}$ and
$$\bm{E}\left[S_n\right]=\left(\frac{n}{2}-\alpha_n\right)\frac{y_{1,0}}{(y_0+\alpha_nb)}\ \ \ \texttt{and}\ \ \ c_n=\frac{p_0}{a(1-p_0)}\frac{(y_0+\alpha_nb)^2}{y_{1,0}(n/2-\alpha_n)};$$
since $c_n\rightarrow0$, we can define an integer $n_0$ such that $c_n<c_0$ for any $n\geq n_0$, so that
$$\bm{P}\left( S_n\leq c_n\cdot\bm{E}[S_n]\right)\ \leq\ \bm{P}\left( S_n\leq c_0\cdot\bm{E}[S_n]\right).$$
Hence, by using~\eqref{eq:bernoulli_chernoff_result}, for any $n\geq n_0$ we have that
$$\bm{P}\left(A_{3,n}\cap A_{4,n}^C\right)\ \leq\ \exp\left(-\frac{(1-c_0)^2}{2}\cdot\bm{E}[S_n]\right),$$
which converges to zero exponentially fast since
$$\bm{E}[S_n]\ =\ \frac{y_{1,0}(n/2-\alpha_n)}{y_0+\alpha_nb}\ \sim\ \frac{n}{\alpha_n}\ =\ n^{1-\alpha}.$$
This concludes the proof.
\endproof

\begin{rem}\label{rem:Y_bounded_and_diverge}
The result of Theorem~\ref{thm:Y_bounded_and_diverge} can be also obtained relaxing assumption~\eqref{eq:rho_exponentially_fast}.
In that case, we need condition~\eqref{eq:rho_convergence_almost_sure} and (c2) to be satisfied.
Then, the proof is the same by setting
\[\begin{aligned}
A_{1,n}\ &&:=&\ \left\{\ \sup_{i\geq \alpha_n}\{\hat{\rho}_{1,i}\}>1-\epsilon\ \right\},\\
A_{2,n}\ &&:=&\ \left\{\ \inf_{i\geq \alpha_n}\{\hat{\rho}_{2,i}\}<\epsilon\ \right\},\\
A_{3,n}\ &&:=&\ \left\{\ \inf_{i\geq \alpha_n}\{\min\{1-\hat{\rho}_{1,i};\hat{\rho}_{2,i}\}\}
\geq\epsilon\ \right\},\\
\end{aligned}\]
where $0<\epsilon<1/2$ is such that $\hat{\rho}_{j,n}\in[\epsilon,1-\epsilon]$ for any $n\geq1$ and $j=1,2$.
Then, $\bm{P}(A_{1,n})=\bm{P}(A_{2,n})=0$ for any $n\geq1$.
\end{rem}

\subsection{A uniform bound}   \label{subsection_asymptotic_difference}

In this subsection, we provide a uniform bound for the scaled difference between $Z_t$ and $\tilde{\rho}_t$ (which is $\mathcal{F}_{t^{-}}$-measurable).
To make precise statements, we start by defining some notations.
Set $\Delta_{j,k}:=\emph{sign}(m_1-m_2)\left(\tilde{\rho}_{q^j+k}-Z_{q^j+k}\right)$ and $\tilde{T}_{j,k}:=Y_{q^j+k}\Delta_{j,k}$,
for any $j\geq1$ and any $k=1,..,d_j$, where $d_j:=q^{j+1}-q^j$.
Note that, since from~\eqref{def:adaptive_thresholds} $\tilde{\rho}_{1,q^j+k}=\hat{\rho}_{1,q^j}$ for any $k\in\{1,..,d_j\}$,
we can also write $\Delta_{j,k}=\emph{sign}(m_1-m_2)\left(\hat{\rho}_{q^j}-Z_{q^j+k}\right)$.
Let $\{\tau_j; j\geq1\}$ be a sequence of stopping times defined as follows:
\begin{equation}\label{def:tauj}
\begin{aligned}
\tau_j & :=
\begin{cases}
\inf\left\{\ k\geq1\ :\ \tilde{T}_{j,k}\in \left[-b,0\right]\ \right\}& \text{if }
\left\{\ k\geq1\ :\ \tilde{T}_{j,k}\in \left[-b,0\right]\ \right\}\neq \emptyset;
\\
\infty & \text{otherwise}.
\end{cases}
\end{aligned}\end{equation}
In Theorem~\ref{thm:how_to_use_branching_process} we provide a $L_1$-uniform bound for the scaled distance among
urn proportion $Z_{q^j+k}$ and the threshold $\tilde{\rho}_{q^j+k}$, on the set $\{\tau_j\leq k\}$.

\begin{thm}\label{thm:how_to_use_branching_process}
Let $\tilde{\rho}_{1,n}$ and $\tilde{\rho}_{2,n}$ be as in~\eqref{def:adaptive_thresholds}.
Then, under assumption~\eqref{ass:different_mean} and~\eqref{eq:rho_exponentially_fast},
there exists a constant $C>0$ such that
\begin{equation}\label{eq:apply_branching_process_result}
\sup_{j\geq1}\sup_{1\leq k\leq d_j}\bm{E}\left[q^{j}\cdot|\Delta_{j,k}|\ind_{\{\tau_j\leq k\}}\right]\ <\ C,
\end{equation}
where $d_j=q^{j+1}-q^j$.
\end{thm}

The proof uses comparison arguments with the MRRU model and related asymptotic results.
Hence, we first present the results concerning the MRRU model in Subsection~\ref{subsubsection_results_MRRU_model}.
The proof of Theorem~\ref{thm:how_to_use_branching_process} is reported in Subsection~\ref{subsubsection_proof_asymptotic_difference}.


\subsubsection{Estimates for the MRRU model}   \label{subsubsection_results_MRRU_model}

In this subsection, we present some probabilistic estimates concerning the MRRU model which are needed in the proof of Theorem~\ref{thm:how_to_use_branching_process}.
We recall that for the MRRU the threshold are fixed, i.e. $\hat{\rho}_{j,n}=\rho_j$ for any $n\geq1$, $j=1,2$.
Hence, in this subsection we consider $W_{1,n}=\ind_{\{Z_n\leq \rho_1\}}$ and $W_{2,n}=\ind_{\{Z_n\geq \rho_2\}}$.
We start by introducing some quantities related to the MRRU model.
Let $\{T_n;n\geq0\}$ be the process defined as
\begin{equation}\label{def:T_n}
T_n:=\emph{sign}(m_1-m_2)\cdot Y_n\left(\rho-Z_n\right),
\end{equation}
which is sometimes useful to represent it as follows:
$$T_n=\emph{sign}(m_1-m_2)\cdot \left(\rho Y_{2,n}-(1-\rho) Y_{1,n}\right).$$
Then, let $t_0$ be the following stopping time
\begin{equation}\label{def:t_0_branching}
t_0\ :=\ \inf\left\{\ k\geq0\ :\ T_k\in\left[-b,0\right]\ \right\}.
\end{equation}
Let
$$S_n\ :=\ \left\{0\leq k\leq n:T_{n-k}\in\left[-b,0\right]\right\},$$
and let $\{s_n;n\geq1\}$ be a sequence of random times defined as
\begin{equation}\label{eq:s_n_lemma}
\begin{aligned}
s_n & =
\begin{cases}
\inf\{S_n\} & \text{if }S_n \neq\emptyset;\\
\infty & \text{otherwise}.
\end{cases}
\end{aligned}\end{equation}
where we recall that $b$ is the maximum value of the urn reinforcements, i.e. $D_{1,n},D_{2,n}\leq b$ a.s. for any $n\geq1$.
Note that by definition $\{s_n=\infty\}=\{t_0>n\}$.
In Theorem~\ref{thm:branching_process} we provide the $L_2$-uniform bound for $Y_n(Z_n-\rho)$, on the set $\{t_0\leq n\}$.

\begin{thm}\label{thm:branching_process}
For an MRRU, under assumption~\eqref{ass:different_mean}, there exists a constant $C>0$ such that
\begin{equation}
\sup_{n\geq1}\bm{E}\left[\left(Y_n|\rho-Z_n|\right)^2\ |\ t_0\leq n \right]\ \leq\ C.
\end{equation}
\end{thm}

The proof uses the boundedness of the moments of the excursion times $s_n$, which is provided in Theorem~\ref{thm:branching_process_times}.
Hence, we first present Theorem~\ref{thm:branching_process_times} and then we report the proof of Theorem~\ref{thm:branching_process}.

\begin{thm}\label{thm:branching_process_times}
For an MRRU, under assumption~\eqref{ass:different_mean}, there exists a constant $C>0$ such that
$$\sup_{n\geq1}\{\bm{E}\left[s_n^2|t_0\leq n\right]\}\ \leq C.$$
\end{thm}

In the proof of Theorem~\ref{thm:branching_process_times},
we need to couple the MRRU model with a particular urn model $\{\tilde{Z}_n;n\geq1\}$.
The processes are coupled, in the sense that: (i)
the potential reinforcements are the same, i.e. $\tilde{D}_{1,n}=D_{1,n}$ and $\tilde{D}_{2,n}=D_{2,n}$ a.s.;
(ii) the drawing process is defined on the same probability space, i.e. $\tilde{U}_n=U_n$ a.s.
where $\{U_n;n\geq1\}$ and $\{\tilde{U}_n;n\geq1\}$ are i.i.d.
uniform random variables such that $X_{n+1}:=\ind_{\{U_{n+1}<Z_n\}}$ and $\tilde{X}_{n+1}:=\ind_{\{\tilde{U}_{n+1}<\tilde{Z}_n\}}$for any $n\geq1$,
respectively.

We now describe the urn model $\{\tilde{Z}_n;n\geq1\}$.
Fix a constant $\tilde{y}_0\in\left(0,Y_0\right]$ and $z_0=\rho_1$.
The process $\{\tilde{Z}_n;n\geq1\}$ evolves as follows:

if $s_{n-1}=0$, i.e. $Z_{n-1}\geq\rho_1$, then $\tilde{X}_n=\ind_{\{\tilde{U}_{n}<\rho_1\}}$ and
\begin{equation}\label{def:urn_coupled_induction_1}\left\{
\begin{array}{l}
\tilde{Y}_{1,n}\ =\ \rho_1\cdot\tilde{y}_0\ +\ \tilde{X}_{n}\tilde{D}_{1,n},\\
\\
\tilde{Y}_{2,n}\ =\ \left(1-\rho_1\right)\cdot \tilde{y}_0+\left(1-\tilde{X}_{n}\right)\tilde{D}_{2,n};
\end{array}
\right.\end{equation}
if $s_{n-1}=k\geq1$, i.e. $Z_{n-1}<\rho_1$, then $\tilde{X}_n=\ind_{\{\tilde{U}_{n}<\tilde{Z}_{n-1}\}}$ and
\begin{equation}\label{def:urn_coupled_induction_2}\left\{
\begin{array}{l}
\tilde{Y}_{1,n}\ =\ \tilde{Y}_{1,n-1}+\tilde{X}_{n}\tilde{D}_{1,n},\\
\\
\tilde{Y}_{2,n}\ =\ \tilde{Y}_{2,n-1}+\left(1-\tilde{X}_{n}\right)\tilde{D}_{2,n};
\end{array}
\right.\end{equation}
where $\tilde{Y}_{n}:=\tilde{Y}_{1,n}+\tilde{Y}_{2,n}$ and $\tilde{Z}_{n}:=\tilde{Y}_{1,n}/\tilde{Y}_{n}$.
The urn model is well defined since $s_{n-1}$ is $\mathcal{F}_{n-1}$-measurable.
It is worth noticing that $\tilde{Z}_n$ represents a Generalized P\'olya urn evaluated after exactly $(s_{n-1}+1)$ steps,
with initial composition $\rho_1\tilde{y}_0$ red and $\rho_1\left(1-\tilde{y}_0\right)$ white balls.\\

In the next lemma, we state an important relation among the MRRU model and the process $\{\tilde{Z}_n;n\geq1\}$,
needed in the proof of Theorem~\ref{thm:branching_process_times}.

\begin{lem}\label{lem:proof_by_induction}
Consider the urn model $\{\tilde{Z}_n;n\geq1\}$ defined in~\eqref{def:urn_coupled_induction_1} and~\eqref{def:urn_coupled_induction_2}
coupled with the MRRU process $\{Z_n;n\geq1\}$.
Let $\tilde{T}_n:=\emph{sign}(m_1-m_2)\cdot \tilde{Y}_n\left(\rho-\tilde{Z}_n\right)$ for any $n\geq1$.
Then, on the set $\{\exists j< n: T_j\leq0\}$, we have that
$$\left\{T_n>0\right\}\ \subset\ \left\{\tilde{T}_n\geq T_n\right\}.$$
\end{lem}

\proof
Wlog assume $m_1>m_2$, which implies $\rho=\rho_1$ and $T_n=Y_n\left(\rho_1-Z_n\right)$.
Sometimes, we will prefer the following expression of $\tilde{T}_n$
$$\tilde{T}_n=\rho_1 \tilde{Y}_{2,n}-(1-\rho_1) \tilde{Y}_{1,n}.$$
The proof will be by induction.
Note that, on the set $\{\exists j< n: T_j\leq0\}$, $s_n$ is almost surely finite.
On the set $\{s_n=0\}$, i.e. $\{T_n\leq0\}$, we can immediately show that $\{T_{n+1}>0\}$ implies $\{\tilde{T}_{n+1}\geq T_{n+1}\}$ and $\{\tilde{Z}_{n+1}\leq Z_{n+1}\}$.
In fact, from $\{T_n\leq0\}$ and $\{T_{n+1}>0\}$ we have $X_{n+1}=0$ and $W_{2,n}=1$, so that
$$T_{n+1}=T_n+\rho_1 D_{2,n+1} \leq \rho_1 D_{2,n+1}=\rho_1 \tilde{D}_{2,n+1}= \tilde{T}_{n+1}$$
and
$$Z_{n+1}=\frac{Z_nY_n}{Y_n+D_{2,n+1}}\geq\frac{\rho_1\tilde{y}_0}{\tilde{y}_0+\tilde{D}_{2,n+1}}=\tilde{Z}_{n+1}.$$
Now, consider the set $\{s_n\geq1\}$ and assume by induction hypothesis that
\begin{equation}\label{hp:ind_3}
\left\{\ \tilde{T}_i\geq T_i>0,\ \tilde{Z}_i\leq Z_i<\rho_1,\ \forall\ i=n-s_n+1,...,n\ \right\}\ \cap\ \left\{T_{n+1}>0\right\}.
\end{equation}
Then, we will show that $\tilde{T}_{n+1}\geq T_{n+1}$ and $\tilde{Z}_{n+1}\leq Z_{n+1}$.
Since $T_n= \rho_1 Y_{2,n}-(1-\rho_1) Y_{1,n}$, we note that
$$T_{n+1}=T_{n-s_n}+\sum_{i=n-s_n+1}^{n+1}\left[\rho_1\left(1-X_i\right)D_{2,i}W_{2,i-1}-\left(1-\rho_1\right)X_iD_{1,i}W_{1,i-1}\right],$$
where we recall that for the MRRU model $W_{1,i}=\ind_{\{Z_i\leq \rho_1\}}$ and $W_{2,i}=\ind_{\{Z_i\geq \rho_2\}}$.
Since $T_{n-s_n}\leq0$, $W_{2,i-1}\leq1$, by~\eqref{hp:ind_3} $W_{1,i}=1$ for any $i=n-s_n+1,...,n$, $X_{n-s_n+1}=0$,
and by construction $\tilde{D}_{1,i}=D_{1,i}$ and $\tilde{D}_{2,i}=D_{2,i}$, we have
$$T_{n+1}\ \leq \sum_{i=n-s_n+1}^{n+1}\left[\rho_1\left(1-X_i\right)\tilde{D}_{2,i}-\left(1-\rho_1\right)X_i\tilde{D}_{1,i}\right].$$
Moreover, by~\eqref{hp:ind_3} we have
$X_{i+1}=\ind_{\{U_{i+1}<Z_i\}}\geq\ind_{\{\tilde{U}_{i+1}<\tilde{Z}_i\}}=\tilde{X}_{i+1}$ for any $i=n-s_n+1,...,n$.
Hence, we can write
$$T_{n+1} \leq \rho_1\sum_{i=n-s_n+1}^{n+1}\left(1-\tilde{X}_i\right)\tilde{D}_{2,i}-\left(1-\rho_1\right)
\sum_{i=n-s_n+1}^{n+1}\tilde{X}_i\tilde{D}_{1,i}\ =\ \tilde{T}_{n+1}.$$
Similarly, we can prove that $\tilde{Z}_{n+1}\leq Z_{n+1}$. Note that
$$Z_{n+1}=\frac{Z_{n-s_n}Y_{n-s_n}+\sum_{i=n-s_n+1}^{n+1}X_iD_{1,i}W_{1,i-1}}
{Y_{n-s_n}+\sum_{i=n-s_n+1}^{n+1}X_iD_{1,i}W_{1,i-1}+\sum_{i=n-s_n+1}^{n+1}\left(1-X_i\right)D_{2,i}W_{2,i-1}}.$$
Now, since $Z_{n-s_n}\geq\rho_1$, $Y_{n-s_n}\geq \tilde{y}_0$ and $X_{i+1}\geq\tilde{X}_{i+1}$ for any $i=n-s_n+1,...,n$, it follows that
$$Z_{n+1}\geq\frac{\rho_1Y_{0}+\sum_{i=n-s_n+1}^{n+1}\tilde{X}_i\tilde{D}_{1,i}}
{\tilde{y}_0+\sum_{i=n-s_n+1}^{n+1}\tilde{X}_i\tilde{D}_{1,i}+
\sum_{i=n-s_n+1}^{n+1}\left(1-\tilde{X}_i\right)\tilde{D}_{2,i}}\ =\ \tilde{Z}_{n+1},$$
which concludes our proof by induction.
\endproof

\proof[Proof of Theorem~\ref{thm:branching_process_times}]
Wlog assume $m_1>m_2$, which implies $\rho=\rho_1$ and $T_n=Y_n\left(\rho_1-Z_n\right)$.
The structure of the proof is the following.
The aim is to show that $\bm{P}\left(s_n=k | t_0\leq n \right)$ converges to zero fast enough such that $\bm{E}\left[s_n^2 | t_0\leq n \right]$ is bounded.
To this end, we consider the urn model $\{\tilde{Z}_n;n\geq1\}$ defined in~\eqref{def:urn_coupled_induction_1} and~\eqref{def:urn_coupled_induction_2}
coupled with the MRRU model,
such that $\bm{P}\left(s_n=k | t_0\leq n \right)$ can be expressed in terms of $\{\tilde{Z}_n;n\geq1\}$.
After some calculations, this is provided by Lemma~\ref{lem:proof_by_induction}.
Moreover, we compare $\{\tilde{Z}_n;n\geq1\}$ with a Generalized P\'olya urn model,
whose moments are uniformly bounded.

First, for any $n\geq1$ note that
$$\bm{E}\left[s_n^2 | t_0\leq n \right]\ =\ \sum_{k=1}^{n}k^2 \bm{P}\left(s_n=k | t_0\leq n \right),$$
since $\bm{P}\left(s_n=\infty | t_0\leq n \right)\ =\ \bm{P}\left(t_0>n | t_0\leq n \right)=0$.
In fact, by definition $t_0\leq n-s_n$ a.s.

Before considering the urn model $\{\tilde{Z}_n;n\geq1\}$, we express $\bm{P}\left(s_n=k | t_0\leq n \right)$ in terms of $\{T_n;n\geq1\}$.
Note that in the MRRU, if $T_j\geq-b$ for some $j\geq0$, then $\bm{P}\left(T_n<-b\right)=0$ for any $n\geq j$.
In fact, when $T_n\geq0$ ($Z_n\leq\rho_1$) we have $T_{n+1}\geq-b$, because the reinforcements are bounded by $b$ and so $|T_{n+1}-T_n|<b$ a.s.;
while when $-b\leq T_n<0$ ($Z_n>\rho_1$) we have  $T_{n+1}\geq T_n\geq-b$, because $Z_n>\rho_1$ implies $W_{1,n}=0$
and so the urn is not reinforced by red balls, i.e. $T_{n+1}\geq T_n$.
As a consequence, since $T_{t_0}\geq-b$ by definition, on the set $\{n\geq t_0\}$, we have $\{T_n\notin[-b,0]\}\subset\{T_n>0\}$.
Hence, since $t_0\leq n-s_n$, we have for all $1\leq k\leq n$
\begin{equation}\label{eq:prob_intersections}\begin{aligned}
\bm{P}\left(s_n=k | t_0\leq n \right)\ &&=&\ \bm{P}\left(\cap_{i=0}^{k-1}\left\{T_{n-i}>0\right\}\cap\left\{T_{n-k}\leq0\right\}|t_0\leq n\right)\\
&&\leq&\ \bm{P}\left(\cap_{i=0}^{k-1}\left\{T_{n-i}>0\right\}|\{T_{n-k}\leq0\}\cap\{t_0\leq n\}\right)\\
&&\leq&\ \bm{P}\left(\cap_{i=0}^{k-1}\left\{T_{n-i}>0\right\}|T_{n-k}\leq0\right),
\end{aligned}\end{equation}
where the last inequality follows from $\{T_{n-k}\leq0\}\subseteq\{t_0\leq n\}$.
To deal with~\eqref{eq:prob_intersections}, we consider the urn model $\{\tilde{Z}_n;n\geq1\}$ defined in~\eqref{def:urn_coupled_induction_1} and~\eqref{def:urn_coupled_induction_2}.
From Lemma~\ref{lem:proof_by_induction}, we have that, on the set $\{\exists j<n: T_j\leq0\}$, the event $\{T_n>0\}$ implies $\{\tilde{T}_n\geq T_n\}$.
Hence, we have that
\begin{equation}\label{eq:intermidiate_1}\begin{aligned}
\bm{P}\left(\cap_{i=0}^{k-1}\left\{T_{n-i}>0\right\}|T_{n-k}\leq0\right)\ &&\leq&\ \bm{P}\left(\cap_{i=0}^{k-1}\left\{\tilde{T}_{n-i}>0\right\}|T_{n-k}\leq0\right)\\
&&=&\ \bm{P}\left(\ \cap_{i=1}^k \left\{Z^G_i<\rho_1\right\}\ \right),
\end{aligned}\end{equation}
by construction, where $\{Z^G_i;i\geq1\}$ is the proportion of red balls of a Generalized P\'olya urn,
starting with a proportion of $Z^G_0=\rho_1$ and an initial number of balls $Y^G_0=\tilde{y}_0$, and
the same reinforcements distributions as $D_{1,n}$ and $D_{2,n}$.

Now, let $s^G$ be the first time the process $Z^G_i$ is above $\rho_1$, i.e.
\[\label{eq:s_G}
\begin{aligned}
s^G & :=
\begin{cases}
\inf\left\{\ i\geq 1\ :\ Z^G_i\geq\rho_1\ \right\} & \text{if }\left\{\ i\geq 1\ :\ Z^G_i\geq\rho_1\ \right\}\neq\emptyset;
\\
\infty & \text{otherwise}.
\end{cases}
\end{aligned}\]
It can be shown using standard arguments that there exists $k_0\in\mathbb{N}$ such that for any $k\geq k_0$,
there exist $0<c_1,c_2<\infty$
$$\bm{P}\left(s^G=k\right)\ \leq\ c_1\exp\left(-c_2 k\right),$$
which implies that $\bm{E}\left[\exp(\gamma s^G)\right]<\infty$ for some $\gamma>0$.

Now, returning to~\eqref{eq:intermidiate_1}, we have that
$$\bm{P}\left(\ \cap_{i=1}^k \left\{Z^G_i<\rho_1\right\}\ \right)\ =\ \bm{P}\left(s^G>k\right)\ \leq\ \frac{\bm{E}\left[\left(s^G\right)^4\right]}{k^4}\ =\
\frac{C_4}{k^4}.$$
Thus we have for any $k\geq1$
$$\bm{P}\left(s_n=k | t_0\leq n \right)\ \leq\ \frac{C_4}{k^4},$$
and hence
\[\begin{aligned}
\bm{E}\left[s_n^2 | t_0\leq n \right]\ &&=&\ \sum_{k=1}^{n}k^2 \bm{P}\left(s_n=k | t_0\leq n \right)\\
&&\leq&\ C_4\cdot\sum_{k=1}^{n} \frac{1}{k^2}\ <\ C\ <\ \infty.
\end{aligned}\]
This concludes the proof.
\endproof

\proof[Proof of Theorem~\ref{thm:branching_process}]
Wlog assume $m_1>m_2$, which implies $\rho=\rho_1$.
Since $T_n=Y_n\left(\rho_1-Z_n\right)$ we want to prove
$$sup_{n\geq1}\ \bm{E}\left[T_n^2 | t_0\leq n \right]\ <\ \infty.$$
Let $s_n$ be the random time defined in~\eqref{eq:s_n_lemma}.
Then, since $|T_{i+1}-T_i|\leq b$ a.s. for any $i\geq1$ and from~\eqref{eq:s_n_lemma} $T_{n-s_n}\in[-b,0]$, we have
\[\begin{aligned}
&&&\bm{E}\left[T_n^2\ | t_0\leq n \right]\ =\ \sum_{l=0}^{n}\bm{E}\left[T_n^2\ | \{s_n=l\}\cap\{t_0\leq n\} \right]\bm{P}\left(s_n=l\ | t_0\leq n \right)\\
&&=&\ b^2\ +\ \sum_{l=1}^{n}\bm{E}\left[\left(\sum_{i=n-l}^{n-1}(T_{i+1}-T_i)+T_{n-l}\right)^2\ | \{s_n=l\} \right]\bm{P}\left(s_n=l\ | t_0\leq n \right)\\
&&\leq&\ \sum_{l=0}^{n}\left(l+1\right)^2b^2 \bm{P}\left(s_n=l\ | t_0\leq n \right).
\end{aligned}\]
Now, using $\left(l+1\right)^2\leq 4 l^2$, we have that
$$\sum_{l=0}^{n}\left(l+1\right)^2b^2 \bm{P}\left(s_n=l | t_0\leq n \right)\ \leq\ 4b^2\cdot\bm{E}\left[s_n^2 | t_0\leq n \right].$$
Finally, using Theorem~\ref{thm:branching_process_times} we have that the last quantity is uniformly bounded by a constant $C$ independent of $n$,
so the proof is concluded.
\endproof

\begin{rem}\label{rem:C_depend_on_D}
From the proof of Theorem~\ref{thm:branching_process_times}, we have that the constant $C$ is independent of the initial proportion $Z_0$.
Moreover, $C$ provides a uniform bound for any other MRRU with initial number of balls $\geq Y_0$.
\end{rem}

\subsubsection{Proof of Theorem~\ref{thm:how_to_use_branching_process}}   \label{subsubsection_proof_asymptotic_difference}

\proof
Wlog, assume $m_1>m_2$, which implies $\rho=\rho_1$.
First, fix $j\in\mathbb{N}$ and apply Cauchy-Schwarz, so obtaining
$$\left(\bm{E}\left[q^{j}\cdot|\Delta_{j,k}|\ind_{\{\tau_j\leq k\}}\right]\right)^2\
\leq\ \bm{E}\left[\left(\tilde{T}_{j,k}\right)^2\ind_{\{\tau_j\leq k\}}\right]
\bm{E}\left[\left(\frac{q^{j}}{Y_{q^{j}}}\right)^2\right].$$
Since $\bm{E}\left[\left(\frac{q^{j}}{Y_{q^{j}}}\right)^2\right]$ is
uniformly bounded by Theorem~\ref{thm:Y_bounded_and_diverge}, it remains to prove that
\begin{equation*}
\bm{E}\left[\left(\tilde{T}_{j,k}\right)^2\ind_{\{\tau_j\leq k\}}\right]\ <\ C,
\end{equation*}
for any $j\geq1$ and any $k=1,..,d_j$.
To this end, fix $j\in\mathbb{N}$ and note that since $\tilde{\rho}_{1,q^j+k}=\hat{\rho}_{1,q^j}$ for any $k\in\{1,..,d_j\}$,
the process $\left\{Z_{q^j+k};k=1,..,d_j\right\}$ can be considered as the urn proportion of the MRRU model,
with initial composition ($Y_{1,q^j}$,$Y_{2,q^j}$) and fixed threshold parameters $\hat{\rho}_{1,q^j}$ and $\hat{\rho}_{2,q^j}$.
Then, for each $j\in\mathbb{N}$ we can apply Theorem~\ref{thm:branching_process}, with $t_0$ defined in~\eqref{def:t_0_branching} equal to $\tau_j$, so obtaining
\begin{equation}
\bm{E}\left[\left(\tilde{T}_{j,k}\right)^2\ind_{\{\tau_j\leq k\}}\right]\ \leq\ C_j,
\end{equation}
where $C_j$ is a constant depending on the initial composition ($Y_{1,q^j}$,$Y_{2,q^j}$).
However, from Remark~\ref{rem:C_depend_on_D} we have that there exists a uniform bound $C>0$ such that $C_j\leq C$ for any $j\geq1$,
since all the processes $\{Z_{q^j+k},k=1,..,d_j\}$ $j\geq1$ can be considered as MRRU with initial number of balls $\geq Y_0$;
this concludes the proof.
\endproof

\section{Proofs of the main results}   \label{section_proofs}

Here, we present the proofs of the results described in Section~\ref{section_the_model}.
Subsection~\ref{subsection_proof_LLN} is dedicated to the proof of Theorem~\ref{thm:convergence_almost_sure} (LLN)
and the related preliminary results.
Then, in subsection~\ref{subsection_proof_CLT} we report the proof of Theorem~\ref{thm:convergence_in_distribution} (CLT)
together with Theorem~\ref{thm:tauj}, a new result needed to compute that proof.
In the last subsections, the proofs of the remaining results of Section~\ref{section_the_model} are gathered.

\subsection{Proof of the LLN}   \label{subsection_proof_LLN}

We start by reporting some preliminary results needed in the proof of the LLN.
Initially, we show that the number of balls sampled from the urn $N_{1,n}$, $N_{2,n}$ and the total number of balls in the urn $Y_n$, increase to infinity almost surely. To do that, we first need to show a lower bound for the increments of the process $Y_n$, which is given by the following:

\begin{lem}\label{lem:lower_bound_Yn}
For any $i\geq1$, we have that
\begin{equation*}
\bm{E}\left[Y_i-Y_{i-1}|\mathcal{F}_{i-1}\right]\ \geq\ a\cdot\left(\frac{\min\{y_{1,0};y_{2,0}\}}{y_0+\left(i-1\right)b}\right).
\end{equation*}
\end{lem}

\proof
First, note that
$$Y_i-Y_{i-1}\ =\ X_{i} D_{1,i} W_{1,i-1}+\left(1-X_i\right) D_{2,i}W_{2,i-1}.$$
Since $X_{i}$ and $D_{1,i}$ are conditionally independent with respect to $\mathcal{F}_{i-1}$, and $W_{1,i-1}$ is $\mathcal{F}_{i-1}$-measurable, we have that
\[\begin{aligned}
\bm{E}\left[Y_i-Y_{i-1}|\mathcal{F}_{i-1}\right]\ &&=&\
\left(m_1Z_{i-1}W_{1,i-1}+m_2\left(1-Z_{i-1}\right)W_{2,i-1}\right)\\
&&\geq&\ a\cdot\left(Z_{i-1}W_{1,i-1}+\left(1-Z_{i-1}\right)W_{2,i-1}\right),
\end{aligned}\]
where the last inequality is because $m_1,m_2\geq a$.
We recall that the variables $W_{1,i-1}$ and $W_{2,i-1}$ can only take the values 0 and 1,
and by construction we have that $W_{1,i-1}+W_{2,i-1}\geq 1$ for any $i\geq1$;
then, we can give a further lower bound
\begin{equation}\label{eq:conditional_increments_Y_i}
\bm{E}\left[Y_i-Y_{i-1}|\mathcal{F}_{i-1}\right]\ \geq\ a\cdot\left(\min\left\{Z_{i-1};1-Z_{i-1}\right\}\right).
\end{equation}
Finally, the result follows by noting that
$$\min\left\{Z_{i-1};1-Z_{i-1}\right\}\ =\ \frac{\min\left\{Y_{1,i-1};Y_{2,i-1}\right\}}{Y_{i-1}}\ \geq\
\frac{\min\{y_{1,0};y_{2,0}\}}{y_0+\left(i-1\right)b}.$$
\endproof

Here, we present the lemma on the divergence of the sequences $Y_n$, $N_{1,n}$ and $N_{2,n}$.
This result is obtained by using the conditional Borel-Cantelli lemma.

\begin{lem}\label{lem:urn_not_stop}
Consider the urn model presented in Section~\ref{section_the_model}. Then,
\begin{itemize}
\item[(a)] $Y_n\stackrel{a.s.}{\rightarrow}\infty$;
\item[(b)] $\min\{N_{1,n};N_{2,n}\}\stackrel{a.s.}{\rightarrow}\infty$.
\end{itemize}
\end{lem}

\proof
\noindent We begin with the proof of part (a).
First, notice that $Y_n=\sum_{i=1}^n\left(Y_i-Y_{i-1}\right)+y_0$.
Then, by Theorem 1 in~\cite{Chen.78}, it is sufficient to show that
$$\left\{\ \omega\in\Omega\ :\ \sum_{i=1}^{\infty}\left[Y_i-Y_{i-1}|\mathcal{F}_{i-1}\right]=\infty\ \right\}$$
occurs with probability one.
To this end, we will now use the lower bound of Lemma~\ref{lem:lower_bound_Yn}, so obtaining
$$\sum_{i=1}^{n}\bm{E}\left[Y_i-Y_{i-1}|\mathcal{F}_{i-1}\right]\
\geq\ a\left(\sum_{i=1}^{n}\frac{\min\{y_{1,0};y_{2,0}\}}{y_0+\left(i-1\right)b}\right)\ \stackrel{a.s.}{\rightarrow}\ \infty.$$
Hence, we have that $Y_{n}\stackrel{a.s.}{\rightarrow}\infty$.

\noindent We now report the proof of part (b).
We will show that $N_{1,n}\stackrel{a.s.}{\rightarrow}\infty$, since the proof for $N_{2,n}$ is analogous.
Since $N_{1,n}=\sum_{i=1}^n X_i$, by Theorem 1 in~\cite{Chen.78}, it is sufficient to show that
$$\left\{\ \omega\in\Omega\ :\ \sum_{i=1}^{\infty}\bm{P}\left(X_i|\mathcal{F}_{i-1}\right)=\infty\ \right\}$$
occurs with probability one.
Then, we obtain
$$\sum_{i=1}^{n}\bm{P}\left(X_i|\mathcal{F}_{i-1}\right)\ =\ \sum_{i=1}^{n}Z_i\ \geq\ \sum_{i=1}^{n}\frac{y_{1,0}}{y_0+\left(i-1\right)b}\ \stackrel{a.s.}{\rightarrow}\ \infty.$$
Hence, we have that $N_{1,n}\stackrel{a.s.}{\rightarrow}\infty$.
\endproof

The following lemma corresponds to Theorem 2.1 of~\cite{Aletti.et.al.13}, and it is needed in the proof of Theorem~\ref{thm:convergence_almost_sure}.
This result provides multiple equivalent ways to show the almost sure convergence of a real-valued process.
We consider a general real-valued process $\{Z_n;n\geq 0\}$ and two real numbers $d$ (down) and $u$ (up), with $d<u$.
The result requires two sequences of times $t_j(d,u)$ and $\tau_j(d,u)$ defined as follows:
for each $j\geq0$, $t_j(d,u)$ represents the time of the first up-cross of $u$ after $\tau_{j-1}(d,u)$, and
$\tau_j(d,u)$ represents the time of the first down-cross of $d$ after $t_j$.
Note that $t_j(d,u)$ and $\tau_j(d,u)$ are stopping times,
since the events $\{t_j(d,u)=k\}$ and $\{\tau_j(d,u)=k\}$ depend on $\{Z_n;n\leq k\}$,
which are measurable with respect to $\mathcal{F}_k$.
We omit the proof since it is reported in Theorem 2.1 of~\cite{Aletti.et.al.13}, using the same notation.

\begin{lem}\label{lem:Zn_converge_ifonlyif}
Let $\{Z_n;n\geq 0\}$ be a real-valued process in $\left[0,1\right]$.
Let $\tau_{-1}(d,u)=-1$ and define for every $j\geq0$ two stopping times
\begin{equation}\label{eq:def_tau}
\begin{aligned}
t_j(d,u) & =
\begin{cases}
\inf\{n>\tau_{j-1}(d,u):Z_n>u\} & \text{if }\{n>\tau_{j}(d,u):Z_n>u\} \neq\emptyset;
\\
+ \infty & \text{otherwise}.
\end{cases}
\\
\tau_j(d,u) & =
\begin{cases}
\inf\{n>t_{j}(d,u) :Z_n<d\}\ \  & \text{if }\{n>t_{j-1}(d,u) :Z_n<d\} \neq\emptyset;
\\
+ \infty & \text{otherwise}.
\end{cases}
\end{aligned}
\end{equation}
Then, the following three events are a.s. equivalent
\begin{itemize}
\item[(a)] $Z_n$ converges a.s.;
\item[(b)] for any $0<d<u<1$,
$$\lim_{j\rightarrow\infty} \bm{P}\left( t_{j}(d,u) < \infty \right) =0;$$
\item[(c)] for any $0<d<u<1$,
$$\sum_{j\geq1} \bm{P}\left( t_{j+1}(d,u) = \infty | t_j(d,u) < \infty\right) =\infty;$$
\end{itemize}
using the convention that $\bm{P}\left(t_{j+1}(d,u) = \infty  | t_j(d,u) < \infty\right) = 1$ when\\
$\bm{P}\left( t_j(d,u) = \infty\right)=1$.
\end{lem}

The following lemma provides lower bounds for the total number of balls in the urn at the times of up-crossings $Y_{t_{j}}$.
The lemma gets used in the proof of Theorem~\ref{thm:convergence_almost_sure},
where conditioning to a fixed number of up-crossing ensures to have at least a number of balls $Y_n$ determined by the lower bounds of this lemma.
This result has been taken by Lemma 2.1 of~\cite{Aletti.et.al.13}.
We omit the proof since adaptive thresholds does not play any role during up-crossings
and the proof reported in Lemma 2.1 of~\cite{Aletti.et.al.13} carries over to our model, with $D_n$ replaced by $Y_n$.

\begin{lem}\label{lem:Y_geometric_increasing}
For any $0<d<u<1$, we have that
$$Y_{t_{j}(d,u)}\geq \left(\frac{u\left(1-d\right)}{d\left(1-u\right)}\right)Y_{t_{j-1}(d,u)} \geq ... \geq \left(\frac{u\left(1-d\right)}{d\left(1-u\right)}\right)^{j}Y_{t_{0}(d,u)}.$$
\end{lem}

The following lemma provides a uniform bound for the generalized P\'{olya} urn with same reinforcement means,
which is needed in the proof of Theorem~\ref{thm:convergence_almost_sure}.
This result has been taken from Lemma 3.2 of~\cite{Aletti.et.al.13}.
The proof is omitted since it is reported in~\cite{Aletti.et.al.13}.
\begin{lem}\label{lem:same_mean}
Consider a generalized P\'{o}lya urn with $m_1 = m_2$. If $Y_0\geq 2b$,
then
$$\bm{P}\left( \sup_{n\geq1}|Z_n-Z_0| \geq h \right) \leq
\frac{b}{Y_0}\left(\frac{4}{h^2}+\frac{2}{h}\right)$$
for every $h>0$.
\end{lem}

Here, we provide the proof of Theorem~\ref{thm:convergence_almost_sure}.

\proof[Proof of Theorem~\ref{thm:convergence_almost_sure}]

Wlog assume $m_1>m_2$, which implies $\hat{\rho}_n=\hat{\rho}_{1,n}$ and $\rho=\rho_1$.
We divide the proof in two steps:
\begin{description}
\item[(a)] $\bm{P}\left(\ \limsup_{n\rightarrow\infty}Z_n =\rho_1\ \right)\ =\ 1,$
\item[(b)] $\bm{P}\left(\ \lim_{n\rightarrow\infty}Z_n \textit{ exists}\ \right)\ =\ 1.$
\end{description}

\noindent \textit{Proof of part} $\bm{(a)}$:\\
We begin by proving that $\bm{P}(\limsup_{n\rightarrow\infty}Z_n \leq\rho_1)=1$.
To this end, we show that there cannot exist $\epsilon>0$ and $\rho^{\prime} >\rho_1$ such that
\begin{equation}\label{eq:hp_assurda_part_a}
\bm{P}\left(\limsup_{n\rightarrow\infty}Z_n>\rho_1^{\prime}\right)\ \geq\ \epsilon\ >\ 0.
\end{equation}
We prove this by contradiction using a comparison argument with a RRU model.
The proof involves last exit time arguments.
Now, suppose~\eqref{eq:hp_assurda_part_a} holds and let $A_1:=\{\limsup_{n\rightarrow\infty}Z_n>\rho_1^{\prime}\}$.
Let
$$R_1:=\left\{\ k\geq0\ :\ \hat{\rho}_{1,k}\geq \frac{\rho_1^{\prime}+\rho_1}{2}\ \right\},$$
and denote the last time the process $\{\hat{\rho}_{1,n};n\geq1\}$ is above $\left(\rho_1^{\prime}+\rho_1\right)/2$ by
\[\begin{aligned}
t_{\frac{\rho_1^{\prime}+\rho_1}{2}} & =
\begin{cases}
\sup\{R_1\} & \text{if }R_1 \neq\emptyset;\\
0 & \text{otherwise}.
\end{cases}
\end{aligned}\]
Since $\hat{\rho}_{1,n}\stackrel{a.s.}{\rightarrow}\rho_1$, then we have that $\bm{P}\left(t_{\frac{\rho_1^{\prime}+\rho_1}{2}}<\infty\right)=1$.
Hence, there exists $n_{\epsilon}\in\mathbb{N}$ such that
\begin{equation}\label{eq:epsilon_mezzi_1}
\bm{P}\left(t_{\frac{\rho_1^{\prime}+\rho_1}{2}}>n_{\epsilon}\right)\ \leq\ \frac{\epsilon}{2}.
\end{equation}
Setting $B_1:=\left\{t_{\frac{\rho_1^{\prime}+\rho_1}{2}}>n_{\epsilon}\right\}$ and using~\eqref{eq:epsilon_mezzi_1}, it follows that
$$\epsilon\ \leq\ \bm{P}\left(A_1\right)\ \leq\ \epsilon/2\ +\ \bm{P}\left(A_1\cap B_1^c\right).$$
Now, we show that $\bm{P}\left(A_1\cap B_1^c\right)=0$.
Setting
$$C_1\ =\ \left\{\ \omega\in\Omega\ :\ \liminf_{n\rightarrow\infty} Z_n<\frac{\rho_1^{\prime}+\rho_1}{2}\ \right\},$$
we decompose $\bm{P}\left(A_1\cap B_1^c\right)$ as follows:
$$\bm{P}\left(A_1\cap B_1^c\right)\ \leq\ \bm{P}\left(E_1\right)\ +\ \bm{P}\left(E_2\right),$$
where $E_1=A_1\cap B_1^c\cap C_1$ and $E_2=A_1\cap B_1^c\cap C_1^c.$

Consider the term $\bm{P}\left(E_2\right)$. Note that on the set $C_1^c$, we have\\ $\left\{\liminf_{n\rightarrow\infty} Z_n\geq\frac{\rho_1^{\prime}+\rho_1}{2} \right\}$ and  on the set $B_1^c$ we have $\{\hat{\rho}_{1,n}\leq \frac{\rho_1^{\prime}+\rho_1}{2}\}$ for any $n\geq n_{\epsilon}$.
Hence, since $B_1^c\cap C_1^c\supset E_2$, on the set $E_2$ we have that $W_{1,n}=\ind_{\{Z_n\leq \hat{\rho}_{1,n}\}}\stackrel{a.s.}{\rightarrow} 0$.
Then, letting $\tau_W:=\sup\{k\geq1:W_{1,k}=1\}$ we have $\bm{P}(E_2\cap\{\tau_W<\infty\})=\bm{P}(E_2)$ and,
on the set $E_2$, for any $n\geq \tau_W$ the ARRU model can be written as follows:
\[\left\{
\begin{array}{l}
Y_{1,n+1}=Y_{1,\tau_W}\\
\\
Y_{2,n+1}=Y_{2,\tau_W}+\sum_{i=\tau_W}^{n+1}\left(1-X_i\right)D_{2,i},
\end{array}
\right.\]
where $W_{1,i-1}=0$ for any $i\geq\tau_W$, and $W_{2,i-1}=1$ because $W_{2,i-1}+W_{2,i-1}\geq1$ by construction.
Now, consider an RRU model $\{Z^{R}_{i};i\geq1\}$ with initial composition $(Y^{R}_{1,0},Y^{R}_{2,0})=(Y_{1,\tau_W},Y_{2,\tau_W})$ a.s.;
the reinforcements are defined as $D^{R}_{1,i}=0$ and $D^{R}_{2,i}=D_{2,\tau_W+i}$ for any $i\geq1$ a.s.;
the drawing process is modeled by $X^R_{i+1}:=\ind_{\{U^R_{i}<Z^R_{i}\}}$ and $U^{R}_{i}=U_{\tau_W+i}$ a.s.,
where $\{U_n;n\geq1\}$ is the sequence such that $X_{n+1}=\ind_{\{U_n<Z_n\}}$ for any $n\geq1$.
Formally, this RRU model can be described for any $n\geq1$ as follows:
\[\left\{
\begin{array}{l}
Y^R_{1,n+1}=Y^R_{1,0}=Y_{1,\tau_W}\\
\\
Y^R_{2,n+1}=Y^R_{2,0}+\sum_{i=0}^{n+1}\left(1-X^R_i\right)D^R_{2,i}=Y_{2,\tau_W}+\sum_{i=\tau_W}^{n+\tau_W+1}\left(1-X_{i}\right)D_{2,i}.
\end{array}
\right.\]
Hence, on the set $E_2$ we have that
$$(Y_{1,n},Y_{2,n})=(Y^{R}_{1,n-\tau_W},Y^{R}_{2,n-\tau_W})\ \ \ a.s.,$$
for any $n\geq\tau_W$.
Since from~\cite{Muliere.et.al.06} $\bm{P}(\limsup_{n\rightarrow\infty} Z^{R}_n=0)=1$,
on the set $E_2$ we have that $\{\limsup_{n\rightarrow\infty} Z_n=0\}$.
This is incompatible with the set $A_1$ which includes $E_2$. Hence $\bm{P}\left(E_2\right)=0$.

We now turn to the proof that $\bm{P}\left(E_1\right)=0$.
To this end, let
$$\tau_{\epsilon}\ :=\ \inf\left\{k\geq n_{\epsilon}\ :\ \left\{Z_k<\frac{\rho_1^{\prime}+\rho_1}{2}\right\}\cap
\left\{Y_k>\frac{b}{(\rho_1^{\prime}-\rho_1)/2}\right\}\ \right\}$$
and note that, since by Lemma~\ref{lem:urn_not_stop} $Y_n\stackrel{a.s.}{\rightarrow}\infty$,
$\bm{P}(C_1\cap\{\tau_{\epsilon}<\infty\})=\bm{P}(C_1)$.
Moreover, on the set $B_1^c$ we have that $\{\hat{\rho}_{1,n}\leq \frac{\rho_1^{\prime}+\rho_1}{2}\}$ for any $n\geq n_{\epsilon}$.
We now show by induction that on the set $B_1^c\cap C_1$ we have $\{Z_{n}<\rho_1^{\prime}\ \forall n\geq \tau_{\epsilon}\}$.
By definition we have $Z_{\tau_{\epsilon}}<\frac{\rho_1^{\prime}+\rho_1}{2}$, and by Lemma~\ref{lem:Y_increments} this implies $Z_{\tau_{\epsilon}+1}<\rho_1^{\prime}$;
now, consider an arbitrary $n> \tau_{\epsilon}$; if $Z_n<\frac{\rho_1^{\prime}+\rho_1}{2}$, then by Lemma~\ref{lem:Y_increments} we have $Z_{n+1}<\rho_1^{\prime}$;
if $\frac{\rho_1^{\prime}+\rho_1}{2}<Z_n<\rho_1^{\prime}$ we have $W_{1,n}=0$ and so $Z_{n+1}\leq Z_n<\rho_1^{\prime}$.
Hence, since $B_1^c\cap C_1\subset E_1$, on the set $E_1$ we have $\{Z_{n}<\rho_1^{\prime}\ \forall n\geq \tau_{\epsilon}\}$.
This is incompatible with the set $A_1$ which also includes $E_1$. Hence $\bm{P}\left(E_1\right)=0$.

Combining all together we have $\epsilon\ \leq\ \epsilon/2 + \bm{P}\left(E_1\right) + \bm{P}\left(E_2\right)\ =\ \epsilon/2$,
which is impossible.
Thus, we conclude that $\bm{P}(A_1^c)=\bm{P}(\limsup_{n\rightarrow\infty}Z_n \leq\rho_1)=1$.\\

\noindent We now prove that $\bm{P}(\limsup_{n\rightarrow\infty}Z_n \geq\rho_1)=1$.
To this end, we now show that there cannot exist $\epsilon>0$ and $\rho^{\prime} <\rho_1$ such that
\begin{equation}\label{eq:hp_assurda_part_b}
\bm{P}\left(\limsup_{n\rightarrow\infty}Z_n<\rho_1^{\prime}\right)\ \geq\ \epsilon\ >\ 0.
\end{equation}
We prove this by contradiction, using a comparison argument with a RRU model.
Now suppose~\eqref{eq:hp_assurda_part_b} holds and let $A_2:=\{\limsup_{n\rightarrow\infty}Z_n<\rho_1^{\prime}\}$.

Let
$$R_2:=\left\{\ k\geq0\ :\ \hat{\rho}_{1,k}< \frac{\rho_1^{\prime}+\rho_1}{2}\ \right\},$$
and define the last time the process $\{\hat{\rho}_{1,n};n\geq1\}$ is less than $\left(\rho_1^{\prime}+\rho_1\right)/2$ by
\[\begin{aligned}
\tau_{\frac{\rho_1^{\prime}+\rho_1}{2}} & =
\begin{cases}
\sup\{R_2\} & \text{if }R_2 \neq\emptyset;\\
0 & \text{otherwise}.
\end{cases}
\end{aligned}\]
Since $\hat{\rho}_{1,n}\stackrel{a.s.}{\rightarrow}\rho_1$, then we have that $\bm{P}\left(\tau_{\frac{\rho_1^{\prime}+\rho_1}{2}}<\infty\right)=1$.
Hence, there exists $n_{\epsilon}\in\mathbb{N}$ such that
\begin{equation}\label{eq:epsilon_mezzi_2}
\bm{P}\left(\tau_{\frac{\rho_1^{\prime}+\rho_1}{2}}>n_{\epsilon}\right)\ \leq\ \frac{\epsilon}{2}.
\end{equation}
Setting $B_2:=\left\{\tau_{\frac{\rho_1^{\prime}+\rho_1}{2}}>n_{\epsilon}\right\}$ and using~\eqref{eq:epsilon_mezzi_2}, it follows that
$$\epsilon\ \leq\ \bm{P}\left(A_2\right)\ \leq\ \epsilon/2\ +\ \bm{P}\left(A_2\cap B_2^c\right).$$
Let $E_3:=A_2\cap B_2^c$. We now show that $\bm{P}\left(E_3\right)$=0.
On the set $A_2$, we have $\left\{\liminf_{n\rightarrow\infty} Z_n\leq\rho_1^{\prime} \right\}$ and
on the set $B_2^c$,
we have $\{\hat{\rho}_{1,n}\geq \frac{\rho_1^{\prime}+\rho_1}{2}\}$ for any $n\geq n_{\epsilon}$.
Hence, on the set $E_3$ we have that $W_{1,n}=\ind_{\{Z_n\leq \hat{\rho}_{1,n}\}}\stackrel{a.s.}{\rightarrow} 1$.
Then, letting $\tau_W:=\sup\{k\geq1:W_{1,n}=0\}$ we have $\bm{P}(E_3\cap\{\tau_W<\infty\})=\bm{P}(E_3)$.
Now, analogously to the proof of $\bm{P}\left(E_2\right)=0$, we can use comparison arguments with the RRU model
to show that on the set $E_3$ we have $\{\limsup_{n\rightarrow\infty} Z_n=1\}$.
This is incompatible with the set $A_2$, which also includes $E_3$. Hence $\bm{P}\left(E_3\right)=0$.

Combining all together we have $\epsilon\ \leq\ \epsilon/2 + \bm{P}\left(E_3\right)\ =\ \epsilon/2$,
which is impossible.
Thus, we conclude that the event $A_2^c=\{\limsup_{n\rightarrow\infty}Z_n \geq\rho_1\}$ occurs with probability one.\\

\noindent \textit{Proof of part} $\bm{(b)}$:\\
In part (a), we have shown that $\bm{P}\left(\limsup_{n\rightarrow\infty}Z_n= \rho_1\right)=1$.
Therefore, if the process $\{Z_n;n\geq1\}$ converges almost surely, then its limit has to be equal to $\rho_1$.
First, let $d$, $u$, $\gamma$ and $\rho_1^{\prime}$ ($d<u<\gamma<\rho_1^{\prime}<\rho_1$) be four constants in $\left(0,1\right)$.
Let $\{\tau_j(d,u);j\geq1\}$ and $\{t_j(d,u);j\geq1\}$ be the sequences of random variables defined in~\eqref{eq:def_tau}.
Since $d$ and $u$ are fixed in this proof, we sometimes denote $\tau_j(d,u)$ by $\tau_j$ and $t_j(d,u)$ by $t_j$.
It is easy to see that $\tau_n$ and $t_n$ are stopping times with respect to $\left\{\mathcal{F}_n;n\geq1\right\}$.

Recall that, by Lemma~\ref{lem:Zn_converge_ifonlyif}, we have that for every $0<d<u<1$
\[\begin{aligned}
Z_n\ \text{converges a.s.}\ &&\Leftrightarrow&\ \bm{P}\left( t_{n}(d,u) < \infty \right)\rightarrow0,\\
&&\Leftrightarrow&\ \sum_{n=1}^{\infty} \bm{P}\left( t_{n+1}(d,u) = \infty | t_n(d,u) < \infty\right) =\infty.
\end{aligned}\]
Now, to prove that $Z_n$ converges a.s., it is sufficient to show that
$$\bm{P}\left( t_n(d,u) < \infty \right)\rightarrow0,$$
for all $0<d<u<1$.
Suppose $Z_n$ does not converges a.s.. This implies that $\bm{P}\left( t_n < \infty \right)\downarrow \phi_1>0$,
since $\bm{P}\left( t_n < \infty \right)$ is a non-increasing sequence.
We will show that for large $j$ there exists a constant $\phi<1$ dependent on $\phi_1$, such that
\begin{equation}\label{eq:goal_a_s_convergence}
\bm{P}\left( t_{j+1} < \infty | t_j<\infty \right) \leq \phi.
\end{equation}
This result implies that $\sum_{n} \bm{P}\left( t_{n+1} = \infty | t_n < \infty\right) =\infty$,
establishing by Lemma~\ref{lem:Zn_converge_ifonlyif} that $\bm{P}\left( t_{n} < \infty \right)\rightarrow0$,
which is a contradiction.\\

\noindent Consider the term $\bm{P}\left( t_{i+1} < \infty | t_i < \infty\right)$.
First, let us denote by $\tau_{\rho_1^{\prime}}$ the last time the process $\hat{\rho}_{1,n}$ is below $\rho_1^{\prime}$, i.e.
\begin{equation*}
\begin{aligned}
\tau_{\rho_1^{\prime}} & =
\begin{cases}
\sup\{n\geq 1 :\hat{\rho}_{1,n}\leq\rho_1^{\prime}\} &
\text{if } \{n\geq 1 :\hat{\rho}_{1,n}\leq\rho_1^{\prime}\}  \neq\emptyset;
\\
0 & \text{otherwise}.
\end{cases}
\end{aligned}
\end{equation*}
Since $\hat{\rho}_{1,n}\stackrel{a.s.}{\rightarrow}\rho_1$, we have that $\bm{P}\left(\tau_{\rho_1^{\prime}}<\infty\right)=1$.
Hence, for any $\epsilon\in\left(0,\frac{1}{2}\right)$ there exists $n_{\epsilon}\in\mathbb{N}$ such that
\begin{equation}\label{eq:epsilon_convergence}
\frac{1}{\phi_1}\bm{P}\left(\tau_{\rho_1^{\prime}}> n_{\epsilon}\right)\ \leq\ \epsilon.
\end{equation}
By denoting $\bm{P_{i}}\left(\cdot\right)=\bm{P}\left(\cdot|t_i<\infty\right)$ and
using $t_i\leq\tau_i\leq t_{i+1}$ we obtain
$$\bm{P}\left(t_{i+1}<\infty|t_i<\infty\right)\ \leq\ \bm{P_{i}}\left(\tau_i<\infty\right).$$
Hence
\begin{equation}\label{eq:decomposition_as}
\bm{P_{i}}\left(\tau_i<\infty\right)\leq\bm{P_{i}}\left(\{\tau_i<\infty\}\cap\{\tau_{\rho_1^{\prime}}\leq n_{\epsilon}\}\right)\ +\
\bm{P_{i}}\left(\tau_{\rho_1^{\prime}}> n_{\epsilon}\right).
\end{equation}
We start with the second term in~\eqref{eq:decomposition_as}. Note that
$$\bm{P_{i}}\left(\tau_{\rho_1^{\prime}}> n_{\epsilon}\right)\ \leq\ \frac{\bm{P}\left(\tau_{\rho_1^{\prime}}> n_{\epsilon}\right)}{\bm{P}\left(t_i<\infty\right)}\ \leq\
\frac{\bm{P}\left(\tau_{\rho_1^{\prime}}> n_{\epsilon}\right)}{\phi_1}\ \leq\ \epsilon,$$
where the last inequality follows from~\eqref{eq:epsilon_convergence}.

Now, consider the first term in~\eqref{eq:decomposition_as}.
Since the probability is conditioned to the set $\{t_i<\infty\}$,
in what follows we will consider the urn process at times $n$ after the stopping time $t_i$.
Since we want to show~\eqref{eq:goal_a_s_convergence} for large $i$,
we can choose an integer $i\geq n_{\epsilon}$ and
$$i> \log_{\frac{u\left(1-d\right)}{d\left(1-u\right)}}\left(\frac{b}{Y_0\left(\gamma-u\right)}\right),$$
so that
\begin{itemize}
\item[(i)] $t_i\geq i\geq n_{\epsilon}$ a.s.;
\item[(ii)] from Lemma~\ref{lem:Y_geometric_increasing}, we have that $Y_{\tau_i} > b/\left(\gamma-u\right)$ a.s.
\end{itemize}
These two properties imply respectively that, on the set $\{n\geq t_i\}$
\begin{itemize}
\item[(i)] $\hat{\rho}_{1,n}\geq\rho_1^{\prime}$, since from $\{\tau_{\rho_1^{\prime}}\leq n_{\epsilon}\}$ we have that $n\geq \tau_{\rho_1^{\prime}}$;
\item[(ii)] $Z_{t_i}\in\left(u,\gamma\right)$, since $Z_{t_i-1}\leq u$ and $Z_{t_i}> u$ and
from Lemma~\ref{lem:Y_increments} we have that $|Z_{n}-Z_{n-1}|<\left(\gamma-u\right)$.
\end{itemize}

Now, let us define two sequences of stopping times $\{t^\ast_n;n\geq1\}$ and $\{\tau^\ast_n;n\geq1\}$,
where $t^\ast_n$ represents the first time after $\tau^\ast_{n-1}$ the process
$Z_{t_i+n}$ up-crosses $\rho_1^{\prime}$,
while $\tau^\ast_n$ represents the first time after $t^\ast_{n}$ the process
$Z_{t_i+n}$ down-crosses $\gamma$.
Formally, let $\tau^\ast_{0}=0$ and define for every $j\geq 1$ two stopping times
\begin{equation}\label{eq:def_tauast2}
\begin{aligned}
t^\ast_j & =
\begin{cases}
\inf\{n> \tau^\ast_{j-1}:Z_{\tau_i+n}>\rho_1^{\prime}\} &
\text{if } \{n>\tau^\ast_{j}:Z_{\tau_i+n}>\rho_1^{\prime}\} \neq\emptyset;
\\
+ \infty & \text{otherwise}.
\end{cases}
\\
\tau^\ast_j & =
\begin{cases}
\inf\{n>t^\ast_{j} :Z_{\tau_i+n}\leq\gamma\} &\ \ \ \
\text{if } \{n>t^\ast_{j-1} :Z_{\tau_i+n}\leq\gamma\} \neq\emptyset;
\\
+ \infty &\ \ \ \ \text{otherwise}.
\end{cases}
\end{aligned}
\end{equation}
Note that, since $Z_{t_i+\tau^\ast_j-1} \geq \gamma$ and $Z_{t_i+\tau^\ast_j} < \gamma$,
from (ii) we have that $Z_{t_i+\tau^\ast_j}\in\left(u,\gamma\right)$.\\

\noindent For any $j\geq0$, let $\{\tilde{Z}^j_n;n\geq1\}$ be a RRU model defined as follows:
\begin{itemize}
\item[(1)] $\left(\tilde{Y}^j_{1,0},\tilde{Y}^j_{2,0}\right) = \left(Y_{1,t_i+\tau^\ast_j},Y_{1,t_i+\tau^\ast_j}\frac{u+d}{2-u-d}\right)$ a.s.,
which implies that $\tilde{Z}^j_0=\frac{u+d}{2}$;
\item[(2)] the drawing process is modeled by
$\tilde{X}^j_{n+1}=\ind_{\{\tilde{U}^j_{n+1}<\tilde{Z}^j_n\}}$, where $\tilde{U}^j_{n+1}=U_{t_i+\tau^\ast_j+n+1}$ a.s.
and $U_n$ is such that $X_{n}=\ind_{\{U_{n}<Z_{n-1}\}}$;
\item[(3)] the reinforcements are defined as
$\tilde{D}^j_{2,n+1}=D_{2,t_i+\tau^\ast_j+n+1}+\left(m_1-m_2\right)$,
$\tilde{D}^j_{1,n+1}=D_{1,t_i+\tau^\ast_j+n+1}$ a.s.;
this means $\bm{E}[\tilde{D}^j_{1,n}]=\bm{E}[\tilde{D}^j_{2,n}]$ for any $n\geq1$;
\item[(4)] the urn process evolves as a RRU model, i.e. for any $n\geq0$
\[\left\{
\begin{array}{l}
\tilde{Y}^j_{1,n+1}=\tilde{Y}^j_{1,n}+\tilde{X}^j_{n+1}\tilde{D}^j_{1,n+1},\\
\tilde{Y}^j_{2,n+1}=\tilde{Y}^j_{2,n}+\left(1-\tilde{X}^j_{n+1}\right)\tilde{D}^j_{2,n+1},\\
\tilde{Y}^j_{n+1}=\tilde{Y}^j_{1,n+1}+\tilde{Y}^j_{2,n+1},\\
\tilde{Z}^j_{n+1}=\frac{\tilde{Y}^j_{1,n+1}}{\tilde{Y}^j_{n+1}}.
\end{array}
\right.\]
\end{itemize}
We will compare the process $\{Z_{t_i+n};n\geq1\}$ with the ARRU process $\{Z_{t_i+n};n\geq1\}$.
Note that at time $n$, we have defined only the processes $\tilde{Z}^j$ such that $\tau^\ast_j<n$.\\

\noindent We will prove, by induction, that on the set $\{\tau_{\rho_1^{\prime}}\leq n_{\epsilon}\}$, for any $j\in\mathbb{N}$
and for any $n\leq {t^\ast_{j+1}}-{\tau^\ast_j}$
\begin{equation}\label{hp:ind2}
\tilde{Z}^j_n < Z_{t_i+\tau^\ast_j+n}, \qquad
\tilde{Y}^j_{2,n} \geq Y_{2,t_i+\tau^\ast_j+n},\qquad
\tilde{Y}^j_{1,n} < Y_{1,t_i+\tau^\ast_j+n}.
\end{equation}
In other words, we will show, provided that $t_i>\tau_{\rho_1^{\prime}}$, that for each $j\geq 1$ the process $\tilde{Z}^j_n$ is always
dominated by the original process $Z_{t_i+\tau^\ast_j+n}$, as long as $Z_{t_i+\tau^\ast_j+n}$ is dominated by $\rho_1^{\prime}$
(i.e. for $n\leq {t^\ast_{j+1}}-{\tau^\ast_j}$).
By construction we have that
\[\tilde{Z}^j_0 = \frac{d+u}{2} < u < Z_{t_i+\tau^\ast_j}, \qquad
\tilde{Y}^j_{1,0} = Y_{1,t_i+\tau^\ast_j}\]
which immediately implies $\tilde{Y}^j_{2,0} > Y_{2,t_i+\tau^\ast_j}$.
To this end, we assume~\eqref{hp:ind2} by induction hypothesis.
First, we will show that $\tilde{Y}^j_{2,n+1} > Y_{2,t_i+\tau^\ast_j+n+1}$.
Since from~\eqref{hp:ind2} $\tilde{Z}^j_n < Z_{t_i+\tau^\ast_j+n}$ for $n\leq {t^\ast_{j+1}}-{\tau^\ast_j}$, by construction we obtain that
$$\tilde{X}^j_{n+1} = \ind_{\{\tilde{U}^j_n<\tilde{Z}^j_n\}} \leq
\ind_{\{U_{t_i+\tau^\ast_j+n}<Z_{t_i+\tau^\ast_j+n}\}} =X_{t_i+\tau^\ast_j+n+1}.$$
As a consequence, since $W_n\leq1$ for any $n\geq1$, we have that
\[\begin{aligned}
\left(Y_{2,t_i+\tau^\ast_j+n+1}-Y_{2,t_i+\tau^\ast_j+n}\right)\ &&=&\ \left(1-X_{t_i+\tau^\ast_j+n+1}\right)D_{2,t_i+\tau^\ast_j+n+1}
W_{2,t_i+\tau^\ast_j+n}\\
&&\leq&\ (1-\tilde{X}^j_{n+1})\tilde{D}^j_{2,n+1}\\
&&=&\ \left(\tilde{Y}^j_{2,n+1}-\tilde{Y}^j_{2,n}\right),
\end{aligned}\]
which using hypothesis~\eqref{hp:ind2} implies $\tilde{Y}^j_{2,n+1} > Y_{2,t_i+\tau^\ast_j+n+1}$.
Similarly, we now show that $\tilde{Y}^j_{1,n+1} \leq Y_{1,t_i+\tau^\ast_j+n+1}$. We have
$$\left(Y_{1,t_i+\tau^\ast_j+n+1}-Y_{1,t_i+\tau^\ast_j+n}\right)\ =\
X_{t_i+\tau^\ast_j+n+1} D_{1,t_i+\tau^\ast_j+n+1}
W_{1,t_i+\tau^\ast_j+n}.$$
From (i) we have that, as long as $Z$ remains below $\rho_1^{\prime}$, $Z$ is also above the process $\hat{\rho}_{1,n}$.
Since we consider the behavior of $Z_{t_i+\tau^\ast_j+n}$ when it is below $\rho_1^{\prime}$, i.e. $n\leq {\tau^\ast_{j+1}}-{t^\ast_j}$,
we have that $W_{1,t_i+\tau^\ast_j+n}=1$. Thus,
$$\left(Y_{1,t_i+\tau^\ast_j+n+1}-Y_{1,t_i+\tau^\ast_j+n}\right)\ \geq\ \tilde{X}^j_{n+1}\tilde{D}^j_{1,n+1}\ =\
\left(\tilde{Y}^j_{1,n+1}-\tilde{Y}^j_{1,n}\right),$$
which using hypothesis~\eqref{hp:ind2} implies $\tilde{Y}^j_{1,n+1} \leq Y_{1,t_i+\tau^\ast_j+n+1}$.
Thus, we have shown that, on the set $\{\tau_{\rho_1^{\prime}}\leq n_{\epsilon}\}$, for any $n\leq {t^\ast_{j+1}}-{\tau^\ast_j}$,
$\tilde{Z}^j_{n+1} < Z_{t_i+\tau^\ast_j+{n+1}}$, $\tilde{Y}^j_{1,n+1} \leq Y_{1,t_i+\tau^\ast_j+{n+1}}$ and
$\tilde{Y}^j_{2,n+1} > Y_{2,t_i+\tau^\ast_j+{n+1}}$ hold.\\

Now, for any $j\geq1$, let $T_j$ be the stopping time for $\tilde{Z}^j_{n}$
to exit from $\left(d,u\right)$, i.e.:
\[T_j =\left\{
\begin{aligned}
\inf \{R_3\}\ &\ \text{if } R_3 \neq\emptyset;\\
+ \infty\ &\ \text{otherwise},
\end{aligned}
\right.\]
where $R_3:=\{n\geq1: \tilde{Z}^j_{n}\leq d\text{ or }\tilde{Z}^j_{n}\geq u\}$.
Note that, on the set $\{\tau_{\rho_1^{\prime}}\leq n_\epsilon\}$,
\[\begin{aligned}
\left\{\tau_i<\infty\right\}\ =\ \left\{\inf_{n\geq1}\left\{Z_{t_i+n}\right\}<d\right\}\ &&\subset&\
\left\{\cup_{j: \tau^\ast_{j}\leq n}\left\{\inf_{n\geq1}\left\{\tilde{Z}^j_{n-\tau^\ast_{j}}\right\}<d\right\}\right\}\\
&&\subset&\ \left\{\cup_{j=0}^{\infty}\left\{T_j<\infty\right\}\right\}.
\end{aligned}\]
Hence,
\[\begin{aligned}
\bm{P_{i}}\left(\{\tau_i<\infty\}\cap\{\tau_{\rho_1^{\prime}}\leq n_{\epsilon}\}\right)\ &&\leq&\
\bm{P_{i}}\left(\left\{\cup_{j=0}^{\infty}\left\{T_j<\infty\right\}\right\}\cap\{\tau_{\rho_1^{\prime}}\leq n_{\epsilon}\}\right)\\
&&\leq&\ \sum_{j=0}^{\infty}\bm{P_{i}}\left(\left\{T_j<\infty\right\}\cap\{\tau_{\rho_1^{\prime}}\leq n_{\epsilon}\}\right).
\end{aligned}\]
Consider a single term of the series; by setting $h=\frac{u-d}{2}$ we get
\[\begin{aligned}
\bm{P_{i}}\left(\left\{T_j<\infty\right\}\cap\{\tau_{\rho_1^{\prime}}\leq n_{\epsilon}\}\right)\ &\leq&&\
\bm{P_{i}}\left(\left\{\sup_{n\geq1}| \tilde{Z}_{n}^j -
\tilde{Z}_{0}^j| \geq  h\right\}\cap\{\tau_{\rho_1^{\prime}}\leq n_{\epsilon}\}\right)\\
&\leq&&\ \bm{P_{i}}\left(\sup_{n\geq1}| \tilde{Z}^j_{n} - \tilde{Z}^j_{0}| \geq  h\right).
\end{aligned}\]
Note that
$\{\tilde{Z}^j_{n};n\geq1\}$
is the proportion of red balls in a RRU model with same reinforcement means.
Then, using Lemma~\ref{lem:same_mean} we obtain
\[\begin{aligned}
\bm{P_{i}}\left(\sup_{n\geq1}| \tilde{Z}^j_{n} - \tilde{Z}^j_{0}| \geq  h\right)\ &=&&\
\bm{E_i}\left[ \bm{P}\left(\left.\left\{\sup_{n\geq1}| \tilde{Z}^j_{n} - \tilde{Z}^j_{0}| \geq  h\right\}\right|\mathcal{F}_{\tau_i+t^\ast_{j}}\right) \right]\\
&\leq&&\
\bm{E_i}\left[\frac{b}{{Y}_{t^\ast_{j}}}\right]\left(\frac{4}{h^2}+\frac{2}{h}\right).
\end{aligned}\]
where $\bm{E_i}\left[\cdot\right]=\bm{E}\left[\cdot|t_i<\infty\right]$.
Moreover, using Lemma~\ref{lem:Y_geometric_increasing}, the right hand side can be expressed as
$$\bm{E_i}\left[\frac{b}{Y_{t_i}}\right]\left(\frac{\rho_1^{\prime}\left(1-\gamma\right)}{\gamma\left(1-\rho_1^{\prime}\right)}\right)^{j}
\left(\frac{4}{h^2}+\frac{2}{h}\right).$$
Since from Lemma~\ref{lem:urn_not_stop} $Y_n$ converges a.s. to infinity, and since $\tau_i\rightarrow\infty$ a.s. because $\tau_i\geq i$,
we have that $\bm{E_i}\left[Y_{t_i}^{-1}\right]$ tends to zero as $i$ increases.
As a consequence, we can choose an integer $i$ large enough such that
$$\bm{E_i}\left[\frac{b}{Y_{t_i}}\right]\left(\frac{4}{h^2}+\frac{2}{h}\right)\left(\frac{1-\rho_1^{\prime}}{1-\rho_1^{\prime}/\gamma}\right)\ <\ \frac{1}{2},$$
which setting $\phi=1/2+\epsilon$ implies~\eqref{eq:goal_a_s_convergence}, i.e.
$$\bm{P}\left(t_{i+1}<\infty|t_i<\infty\right)\ \leq\ \phi\ <\ 1.$$
This concludes the proof.
\endproof

\proof[Proof of Corollary~\ref{cor:convergence_almost_sure_NR}]
This corollary has been proved in Proposition 2.1 of~\cite{Ghiglietti.et.al.12} for the MRRU.
That proof is only based on the fact that the urn proportion $Z_n$ converges a.s. to a value within the interval $(0,1)$,
while the reinforcement rules do not play any role.
Hence, the proof used in~\cite{Ghiglietti.et.al.12} can be applied to the ARRU,
since $Z_n\stackrel{a.s.}{\rightarrow}\rho\in(0,1)$ for ARRU using Theorem~\ref{thm:convergence_almost_sure}.
\endproof

\subsection{Proof of the central limit theorem}   \label{subsection_proof_CLT}


Before the proof of Theorem~\ref{thm:convergence_in_distribution},
we recall that $\{\tau_j; j\geq1\}$ is the sequence defined in~\eqref{def:tauj} as follows:
\[\begin{aligned}
\tau_j & :=
\begin{cases}
\inf\left\{\ k\geq1\ :\ \tilde{T}_{j,k}\in \left[-b,0\right]\ \right\}& \text{if }
\left\{\ k\geq1\ :\ \tilde{T}_{j,k}\in \left[-b,0\right]\ \right\}\neq \emptyset;
\\
\infty & \text{otherwise}.
\end{cases}
\end{aligned}\]
Fix $\nu\in(0,1/2)$ and, for any $j\geq1$, let $r_j:=q^{j\frac{1+\nu}{2}}$ and $\mathcal{R}_j:=\{\tau_j>r_j\}$.
The following theorem is critical to the proof of Theorem~\ref{thm:convergence_in_distribution}.
\begin{thm}\label{thm:tauj}
Let $\tilde{\rho}_{1,n}$ and $\tilde{\rho}_{2,n}$ be as in~\eqref{def:adaptive_thresholds}.
Then, under assumption~\eqref{ass:different_mean} and~\eqref{eq:rho_exponentially_fast},
we have that
\begin{equation}\label{eq:tauj_almsot_sure}
\bm{P}\left(\mathcal{R}_j,i.o.\right)=0.
\end{equation}
\end{thm}
We delay the proof of this theorem to Subsection~\ref{subsubsection_proof_tau_j}.

\proof[Proof of Theorem~\ref{thm:convergence_in_distribution}]
Wlog assume $m_1>m_2$, which implies $\rho=\rho_1$.
To prove the main result, we establish
\begin{description}
\item[(a)] $\sqrt{n}\left(\frac{N_{1,n}}{n}-\frac{\sum_{i=1}^nZ_{i-1}}{n}\right)\ \stackrel{d}{\rightarrow}\ \mathcal{N}\left(0,\rho_1\left(1-\rho_1\right)\right)$, and
\item[(b)] $\sqrt{n}\left(\frac{\sum_{i=1}^nZ_{i-1}}{n}-\frac{\sum_{i=1}^{n}\tilde{\rho}_{1,i-1}}{n}\right)\ \stackrel{a.s.}{\rightarrow}\ 0$.
\end{description}
Finally, result~\eqref{eq:X_eta_convergence_distribution} is obtained by using Slutsky's Theorem to combine (a) and (b) together.\\

\noindent \textit{Proof of part} $\bm{(a)}$:
Let us define a random variable $J_{ni}:=\frac{1}{\sqrt{n}}\left(X_i-\bm{E}\left[X_i|\mathcal{F}_{i-1}\right]\right)$,
for any $n,i\in\mathbb{N}$ with $i\leq n$.
Then, for each $n\in\mathbb{N}$, the sequence\\
$\left\{S_{nj}=\sum_{i=1}^jJ_{ni};\ 1\leq j\leq n\right\}$ is a martingale.
Now we apply the Martingale CLT (MCLT).
First note that $J_{ni}^2\leq 1/n$ for any $n\in\mathbb{N}$ and $|J_{ni}|<\epsilon$ for any $n\geq\epsilon^{-2}$; thus
$$\sum_{i=1}^n\ \bm{E}\left[\ J_{ni}^2\ind_{\{|J_{ni}|>\epsilon\}}\ |\ \mathcal{F}_{i-1}\ \right]\ \leq\ \sum_{i=1}^{\left[\epsilon^{-2}\right]+1}1/n\ =\ \frac{\left[\epsilon^{-2}\right]+1}{n}\ \rightarrow\ 0.$$
Also,
\[\begin{aligned}
\bm{E}\left[\ J_{ni}^2\ |\mathcal{F}_{i-1}\ \right]\ &&=&\ \frac{1}{n}\cdot\bm{E}\left[\ \left(X_{ni}-\bm{E}\left[X_{ni}|\mathcal{F}_{i-1}\right]\right)^2\ |\ \mathcal{F}_{i-1}\ \right]\\
&&=&\ \frac{Z_{i-1}\left(1-Z_{i-1}\right)}{n};
\end{aligned}\]
since $\hat{\rho}_{1,n}\stackrel{a.s.}{\rightarrow}\rho_1$,
from Theorem~\ref{thm:convergence_almost_sure} we get $Z_n\stackrel{a.s.}{\rightarrow}\rho_1$, which implies
\begin{equation*}
\sum_{i=1}^n\ \bm{E}\left[\ J_{ni}^2\ |\ \mathcal{F}_{i-1}\ \right]\ =\ \frac{\sum_{i=1}^nZ_{i-1}\left(1-Z_{i-1}\right)}{n}\ \stackrel{a.s.}{\rightarrow}\ \rho_1\left(1-\rho_1\right).
\end{equation*}
From MCLT~\cite{Hu.et.al.06}, it follows that
\[\begin{aligned}
\frac{1}{\sqrt{n}}\cdot\sum_{i=1}^n\ \left(X_i-\bm{E}\left[X_i|\mathcal{F}_{i-1}\right]\right)\ &&=&\ \sqrt{n}\left(\frac{\sum_{i=1}^nX_i}{n}-\frac{\sum_{i=1}^nZ_{i-1}}{n}\right)\\
&&\stackrel{d}{\rightarrow}&\ \mathcal{N}\left(0,\rho_1\left(1-\rho_1\right)\right).
\end{aligned}\]
We now turn to the proof of part $\bm{(b)}$.
We first express
\[\begin{aligned}
\sqrt{n}\left(\frac{\sum_{i=1}^nZ_{i-1}}{n}-\frac{\sum_{i=1}^{n}\tilde{\rho}_{1,i-1}}{n}\right)\ &&=&\ \frac{1}{\sqrt{n}}\sum_{i=0}^{n-1}\left(Z_{i}-\tilde{\rho}_{1,i}\right)\\
&&=&\ B_{1n}\ +\ B_{2n},
\end{aligned}\]
where
\[\begin{aligned}
B_{1n}\ &&:=&\ \frac{1}{\sqrt{n}}\sum_{i=0}^{\left[q^{k_n}\right]}\left(Z_{i}-\tilde{\rho}_{1,i}\right),\\
B_{2n}\ &&:=&\ \frac{1}{\sqrt{n}}\sum_{i=\left[q^{k_n}\right]+1}^{n-1}\left(Z_{i}-\tilde{\rho}_{1,i}\right),
\end{aligned}\]
and we recall $k_n$ is defined in~\eqref{def:k_n} as $k_n:=[log_{q}(n)]$.
We begin with $B_{1n}$. Note that
$$\sum_{i=0}^{\left[q^{k_n}\right]}\left(Z_{i}-\tilde{\rho}_{1,i}\right)\
=\ \sum_{j=1}^{k_n-1}\sum_{i=1}^{d_j}\left(Z_{q^{j}+i}-\hat{\rho}_{1,q^j}\right)\ =
\ \sum_{j=1}^{k_n-1}\sum_{i=1}^{d_j}\left(-\Delta_{j,i}\right),$$
where we recall that $d_j=q^{j+1}-q^j$ and $\Delta_{j,i}=\hat{\rho}_{1,q^j}-Z_{q^{j}+i}$ for any $j\geq1$ and $1\leq i\leq d_j$.
Hence
$$|B_{1n}|\ =\ \frac{1}{\sqrt{n}}\cdot \left|\sum_{j=1}^{k_n-1}\sum_{i=1}^{d_j}\left(-\Delta_{j,i}\right)\right|
\leq\ \frac{1}{\sqrt{n}}\cdot \sum_{j=1}^{k_n-1}\left(\frac{\sum_{i=1}^{d_j}|\Delta_{j,i}|}{\sqrt{d_j}}\right)
\sqrt{d_j};$$
similarly
$$|B_{2n}|\ \leq\ \frac{1}{\sqrt{n}}\cdot \left(\frac{\sum_{i=1}^{d_{k_n}}|\Delta_{k_n,i}|}{\sqrt{d_{k_n}}}\right)\sqrt{d_{k_n}}.$$
Now, for any $j\geq 1$ define
\begin{equation}\label{bj}
b_j\ :=\ \frac{\sum_{i=1}^{d_j}|\Delta_{k_n,i}|}{\sqrt{d_j}},
\end{equation}
it follows that
$$|B_{1n}|+|B_{2n}|\ \leq\ \frac{1}{\sqrt{n}}\cdot \sum_{j=1}^{k_n}b_j\sqrt{d_j}.$$
Now, we have
\[\begin{aligned}
|B_{1n}|+|B_{2n}|\ &&\leq&\ \frac{1}{\sqrt{n}}\cdot \sum_{j=1}^{k_n/2-1}b_j\sqrt{d_j}\ +\ \frac{1}{\sqrt{n}}\cdot \sum_{j=k_n/2}^{k_n}b_j\sqrt{d_j}\\
&&\leq&\ \left(\frac{\sup_{i\geq1}\{b_i\}}{\sqrt[4]{n}}\right)\cdot H_{1n}\ +\ \left(\sup_{i\geq k_n/2}\{b_i\}\right)\cdot H_{2n}.
\end{aligned}\]
where
$$H_{1n}:=\frac{1}{\sqrt[4]{n}}\sum_{j=1}^{k_n/2-1}\sqrt{d_j},\ \ \ \ \ \ H_{2n}:=\frac{1}{\sqrt{n}}\sum_{j=1}^{k_n}\sqrt{d_j}.$$
Using $d_j=(q-1)q^j$ we express
\[\begin{aligned}
H_{1n}&&=&\frac{\sqrt{q-1}}{\sqrt[4]{n}}\cdot \sum_{j=1}^{k_n/2-1}\left(\sqrt{q}\right)^j=
\left(\frac{\sqrt{q}^{k_n/2}-1}{\sqrt[4]{n}}\right)\cdot\left(\frac{\sqrt{q-1}}{\sqrt{q}-1}\right),\\
H_{2n}&&=&\frac{\sqrt{q-1}}{\sqrt{n}}\cdot \sum_{j=1}^{k_n}\left(\sqrt{q}\right)^j=
\left(\frac{\sqrt{q}^{k_n+1}-\sqrt{q}^{k_n/2}}{\sqrt{n}}\right)\cdot\left(\frac{\sqrt{q-1}}{\sqrt{q}-1}\right).
\end{aligned}\]
Since $n\geq q^{k_n}$, it follows that $H_{1n}\leq C$ and $H_{2n}\leq \sqrt{q}C$,
where
$C=\left(\frac{\sqrt{q-1}}{\sqrt{q}-1}\right)$.
Thus,
$$|B_{1n}|+|B_{2n}|\ \leq\ \left(\frac{\sup_{i\geq1}\{b_i\}}{\sqrt[4]{n}}\right)\cdot C\ +\ \left(\sup_{i\geq k_n/2}\{b_i\}\right)\cdot \sqrt{q}C.$$
To conclude the proof we will show that $b_j\stackrel{a.s.}{\rightarrow}0$.

\noindent First, fix an arbitrary constant $\nu\in\left(0,1/2\right)$ and let $r_j:=q^{j\frac{1+\nu}{2}}$ for any $j\geq1$; then, write
\[\begin{aligned}
b_j\ =\ \frac{1}{\sqrt{d_j}}\sum_{i=1}^{d_j}|\Delta_{j,i}|\
&&=&\ \left(\frac{1}{\sqrt{d_j}}\sum_{i=1}^{r_j}|\Delta_{j,i}|\right)\ +\
\left(\frac{1}{\sqrt{d_j}}\sum_{i=r_j+1}^{d_j}|\Delta_{j,i}|\right)\\
&&=&\ F_{1j}\ +\ F_{2j},
\end{aligned}\]
Let us consider term $F_{1j}$, we have that
$$F_{1j}\
=\ \frac{r_j}{\sqrt{d_j}}\cdot
\left(\frac{1}{r_j}\sum_{i=1}^{r_j}|\Delta_{j,i}|\right)
=\frac{\left[q^{j\frac{\nu}{2}}\right]}{\sqrt{q-1}}\cdot
\left(\frac{1}{r_j}\sum_{i=1}^{r_j}|\Delta_{j,i}|\right),$$
since $d_j=(q-1)q^j$ and $r_j/\sqrt{q^j}=q^{j\frac{\nu}{2}}$.
Now, for any $i=1,..,r_j$ we note that
$$|\Delta_{j,i}|\ \leq\ |Z_{q^{j}+i}-Z_{q^j}|\ +\ |\Delta_{j-1,d_{j-1}}|\ +\ |\hat{\rho}_{1,q^{j-1}}-\hat{\rho}_{1,q^j}|;$$
hence, we have
$$F_{1j}\ \leq\ E_{1j}\ +\ E_{2j}\ +\ E_{3j},$$
where
\[\begin{aligned}
E_{1j}\ &&:=&\ \frac{\left[q^{j\frac{\nu}{2}}\right]}{\sqrt{q-1}}\cdot\left(\
\frac{1}{r_j}\sum_{i=1}^{r_j}|Z_{q^{j}+i}-Z_{q^j}|\right),\\
E_{2j}\ &&:=&\ \frac{\left[q^{j\frac{\nu}{2}}\right]}{\sqrt{q-1}}\cdot|\Delta_{j-1,d_{j-1}}|,\\
E_{3j}\ &&:=&\ \frac{\left[q^{j\frac{\nu}{2}}\right]}{\sqrt{q-1}}\cdot|\hat{\rho}_{1,q^{j-1}}-\hat{\rho}_{1,q^j}|.
\end{aligned}\]
Let us consider the term $E_{1j}$.
Since from Lemma~\ref{lem:Y_increments} we have $|Z_k-Z_{k-1}|\leq b/Y_{k-1}$, we have that
$$E_{1j}\ \leq\ \frac{\left[q^{j\frac{\nu}{2}}\right]}{\sqrt{q-1}}\cdot\frac{br_j}{Y_{q^j}}\ =\
\left(\frac{b}{\sqrt{q-1}}\right)\cdot\left(\frac{q^{j\left(\frac{1}{2}+\nu\right)}}{Y_{q^j}}\right).$$
Then, by using Markov's inequality we obtain
$$\sum_{j=1}^{\infty}\bm{P}\left(\frac{q^{j\left(\frac{1}{2}+\nu\right)}}{Y_{q^j}}>\epsilon\right)\ \leq\
\frac{1}{\epsilon}\sum_{j=1}^{\infty}\bm{E}\left[\frac{q^j}{Y_{q^j}}\right]q^{-j\left(\frac{1}{2}-\nu\right)}\ \leq\
\frac{C}{\epsilon}\sum_{j=1}^{\infty}q^{-j\left(\frac{1}{2}-\nu\right)}\ <\ \infty,$$
where $C=\sup_{k\geq1}\{\bm{E}[k/Y_{k}]\}$ is finite from Theorem~\ref{thm:Y_bounded_and_diverge}.
Thus, from the Borel-Cantelli lemma it follows that $E_{1j}\stackrel{a.s.}{\longrightarrow}\ 0$.\\

Now, consider the term $E_{2j}$. We have
\[\begin{aligned}
\bm{P}\left(\lim_{k\rightarrow\infty}\cup_{j\geq k}\{E_{2j}>\epsilon\}\right)\ &&\leq&\ \bm{P}\left(\lim_{k\rightarrow\infty}\cup_{j\geq k}\mathcal{R}_j\right)\\
&&+&\ \bm{P}\left(\lim_{k\rightarrow\infty}\cup_{j\geq k}\left\{\frac{\left[q^{(j+1)\frac{\nu}{2}}\right]}{\sqrt{q-1}}\cdot|\Delta_{j,d_j}|\ >\ \epsilon \right\}\cap\mathcal{R}_j^c\right).
\end{aligned}\]
where the term $\bm{P}\left(\lim_{k\rightarrow\infty}\cup_{j\geq k}\mathcal{R}_j\right)=0$ from Theorem~\ref{thm:tauj}.
Then, by using Markov's inequality we obtain
$$\sum_{j=1}^{\infty}\bm{P}\left(\left\{\frac{\left[q^{(j+1)\frac{\nu}{2}}\right]}{\sqrt{q-1}}\cdot|\Delta_{j,d_j}|\ >\ \epsilon \right\}\cap\mathcal{R}_j^c\right)\ \leq\ M,$$
where
$$M\ :=\ \frac{1}{\epsilon}\sum_{j=1}^{\infty}\bm{E}\left[\frac{\left[q^{(j+1)\frac{\nu}{2}}\right]}{\sqrt{q-1}}\cdot|\Delta_{j,d_j}|
\ind_{\mathcal{R}_j^C}\right].$$
Now, for any $j\geq1$ let us introduce the set $\mathcal{Q}_j:=\{\tau_j>d_j\}$.
Using $\mathcal{R}_j^C\subseteq\mathcal{Q}_j^C$ from $r_j\leq d_j$, and
by multiplying and dividing by $q^{j+1}$, we have that
\[\begin{aligned}
M\ &&=&\ \frac{1}{\epsilon\sqrt{q-1}}\sum_{j=1}^{\infty}\bm{E}\left[q^{j+1}|\Delta_{j,d_j}|
\ind_{\mathcal{R}_j^C}\right]\cdot q^{-\left(j+1\right)\left(1-\frac{\nu}{2}\right)}\\
&&\leq&\ \frac{1}{\epsilon\sqrt{q-1}}\sum_{j=1}^{\infty}\bm{E}\left[q^{j+1}|\Delta_{j,d_j}|
\ind_{\mathcal{Q}_j^C}\right]\cdot q^{-\left(j+1\right)\left(1-\frac{\nu}{2}\right)}\\
&&\leq&\ \frac{1}{\epsilon\sqrt{q-1}}\left(\sup_{k\geq1}\left\{\bm{E}\left[q^{k+1}|\Delta_{k,d_k}|\ind_{\mathcal{Q}_k^C}\right]\right\}
\right)\sum_{j=1}^{\infty}q^{-\left(j+1\right)\left(1-\frac{\nu}{2}\right)}\\
&&<&\ \infty,
\end{aligned}\]
using Theorem~\ref{thm:how_to_use_branching_process} and the result follows from the Borel-Cantelli lemma.\\

Let us consider the term $E_{3j}$. For any $\epsilon>0$, by using Markov's inequality we have
$$\bm{P}\left(E_{3j}>\epsilon\right)\ \leq\
\frac{1}{\epsilon\sqrt{q-1}}\bm{E}\left[q^{j\frac{\nu}{2}}\cdot|\hat{\rho}_{1,q^{j}}-\hat{\rho}_{1,q^{j-1}}|\right].$$
The right-hand side (RHS) of the above expression can be rewritten as
$$\frac{q^{-j\left(\frac{1-\nu}{2}\right)}}{\epsilon\sqrt{q-1}}
\bm{E}\left[q^{\frac{j}{2}}\cdot|\hat{\rho}_{1,q^{j}}-\hat{\rho}_{1,q^{j-1}}|\right],$$
Now, by decomposing the last expectation into
$$\bm{E}\left[q^{\frac{j}{2}}\cdot|\hat{\rho}_{1,q^{j}}-\hat{\rho}_{1,q^{j-1}}|\right]\ =\
\bm{E}\left[q^{\frac{j}{2}}\cdot|\rho_1-\hat{\rho}_{1,q^{j-1}}|\right]
+\bm{E}\left[q^{\frac{j}{2}}\cdot|\rho_1-\hat{\rho}_{1,q^{j}}|\right]$$
we can see that
$$\sum_{j=1}^{\infty}\bm{P}\left(E_{3j}>\epsilon\right)\ \leq\ \left(\frac{2\sup_{k\geq1}\left\{\bm{E}\left[q^{\frac{k}{2}}\cdot|\rho_1-\hat{\rho}_{1,q^{k}}|\right]\right\}}
{\epsilon\sqrt{q-1}}\right)\sum_{j=1}^{\infty}q^{-j\left(\frac{1-\nu}{2}\right)},$$
which is finite because of~\eqref{eq:rho_exponentially_fast}.
Hence, by another application the Borel-Cantelli lemma, $E_{3j}\stackrel{a.s.}{\longrightarrow}\ 0$;
then, we have $F_{1j}\stackrel{a.s.}{\longrightarrow}\ 0$.\\

\noindent Finally, let us consider term $F_{2j}$.
First, we multiply and divide by $\left(d_j-r_j\right)q^{-\frac{j}{2}}$ to obtain $F_{2j}=c_jF_{3j}$,
where
$$c_j=\frac{d_j-r_j}{q^{\frac{j}{2}}\sqrt{d_j}},\ \ \ \ \ \
F_{3j}=\frac{1}{d_j-r_j}\sum_{i=r_j+1}^{d_j}q^{\frac{j}{2}}|\Delta_{j,i}|.$$
Since $c_j\rightarrow\sqrt{q-1}$, let us focus on $F_{3j}$.
Since $\bm{P}\left(\mathcal{R}_j,i.o.\right)=0$ (Theorem~\ref{thm:tauj}),
it is sufficient to show that $F_{3j}\ind_{R_j^c}\stackrel{a.s.}{\longrightarrow}0$.
For any $\epsilon>0$, by Markov's inequality it follows that
$$\bm{P}\left(\left\{F_{3j}>\epsilon\right\}\cap \mathcal{R}_j^C\right)\ \leq\
\frac{1}{\epsilon}\left(\frac{1}{d_j-r_j}
\sum_{i=r_j+1}^{d_j}\bm{E}\left[q^{\frac{j}{2}}|\Delta_{j,i}| \ind_{\mathcal{R}_j^C}\right]\right),$$
Now, since
$$\sum_{i=r_j+1}^{d_j}\bm{E}\left[q^{\frac{j}{2}}|\Delta_{j,i}|\ind_{\mathcal{R}_j^C}\right]\ \leq\
(d_j-r_j)\left(\max_{i=r_j+1,..,d_j}\left\{\bm{E}\left[q^{\frac{j}{2}}|\Delta_{j,i}|
\ind_{\mathcal{R}_j^C}\right]\right\}\right),$$
we have that
\[\begin{aligned}
\bm{P}\left(\left\{F_{3j}>\epsilon\right\}\cap \mathcal{R}_j^C\right)\
&&\leq&\ \frac{1}{\epsilon}\left(\max_{i=r_j+1,..,d_j}\left\{\bm{E}\left[q^{\frac{j}{2}}|\Delta_{j,i}|
\ind_{\mathcal{R}_j^C}\right]\right\}\right)\\
&&=&\ \frac{1}{\epsilon}
\left(\max_{i=r_j+1,..,d_j}\left\{\bm{E}\left[q^j|\Delta_{j,i}|
\ind_{\mathcal{R}_j^C}\right]\right\}\right)q^{-\frac{j}{2}}\\
&&\leq&\ \frac{1}{\epsilon}\left(\sup_{k\geq1}\left\{\max_{i=\left[r_k\right]+1,..,d_k}
\left\{\bm{E}\left[q^k|\Delta_{k,i}|
\ind_{\mathcal{R}^C_k}\right]\right\}\right\}\right)q^{-\frac{j}{2}}\\
&&\leq&\ Cq^{-\frac{j}{2}},
\end{aligned}\]
where the last inequality follows from Theorem~\ref{thm:how_to_use_branching_process}.
Now, summing over $j$ we have that
$$\sum_{j=1}^{n}\bm{P}\left(\left\{F_{3j}>\epsilon\right\}\cap \mathcal{R}_j^C\right)\ \leq\ C\sum_{j=1}^{n}q^{-\frac{j}{2}}\ <\ \infty.$$
Now, using the Borel-Cantelli lemma we get that $F_{2j}\stackrel{a.s.}{\longrightarrow}\ 0$, which concludes the proof.
\endproof

\subsubsection{Proof of Theorem~\ref{thm:tauj}}   \label{subsubsection_proof_tau_j}

\proof
Wlog assume $m_1>m_2$, which implies $\hat{\rho}_n=\hat{\rho}_{1,n}$ and $\rho=\rho_1$.
To prove~\eqref{eq:tauj_almsot_sure} we need to study the sequence of sets $\{\mathcal{R}_j;j\geq1\}$.
On the set $\mathcal{R}_j$, the urn proportion does not cross the thresholds at times $q^j,..,q^j+r_j$.
Hence, $\mathcal{R}_j$ will be included in $\mathcal{A}_{j}\cup\mathcal{B}_{j}$,
where $\mathcal{A}_{j}$ and $\mathcal{B}_{j}$ represent the events in which the urn proportion is always above and below, respectively,
the thresholds at times $q^j,..,q^j+r_j$.
To show that $\mathcal{A}_{j}$ and $\mathcal{B}_{j}$ cannot occur i.o., we need to appropriately express them
by using the following scaling processes:
\begin{itemize}
\item[(a)] $\tilde{T}_{j,k}=Y_{q^j+k}\Delta_{j,k}=Y_{q^j+k}\left(\hat{\rho}_{1,q^j}-Z_{q^j+k}\right)$,
defined for any $j\geq1$ and any $k=1,..,d_j$. This process models the closeness among the urn proportion and the adaptive threshold.
\item[(b)] $T_n=Y_n\left(\rho_1-Z_n\right)$, defined for any $n\geq1$. This process models the closeness among the urn proportion and the limit of the threshold's sequence.
\item[(c)]  $T^{(\rho_1)}_{j,k}:=Y_{q^j+k}\left(\rho_1-\hat{\rho}_{1,q^j}\right)$, defined for any $j\geq1$ and $k=1,..,d_j$.
 This process models the closeness between the adaptive threshold and its limit.
\end{itemize}
Let us now define formally the sets $\mathcal{A}_{j}$ and $\mathcal{B}_{j}$.
First, note that if the urn proportion crosses the threshold at time $(q^j+k)$, then
$\tilde{T}_{q^{j}+k}\cdot\tilde{T}_{q^{j}+k-1}<0$, since only one among
$\tilde{T}_{q^{j}+k}$ and $\tilde{T}_{q^{j}+k-1}$ is within the interval $[-b,0]$.
Thus, from the definition of $\tau_j$ in~\eqref{def:tauj}, we have that
$$\left\{\Delta_{j,k-1}\cdot\Delta_{j,k}<0\right\}\subseteq
\left\{\tau_j\leq k\right\}.$$
This implies that
\begin{equation*}
\begin{aligned}
\mathcal{R}_j\ &&\subset&\ \left\{\cap_{k=1}^{r_j}\left\{\Delta_{j,k-1}\cdot\Delta_{j,k}>0\right\}\right\}\\
&&=&\ \left\{\ \cap_{k=1}^{r_j}\left\{\Delta_{j,k}<0\right\}\ \right\}
\cup\ \left\{\ \cap_{k=1}^{r_j}\left\{\Delta_{j,k}>0\right\}\right\}.
\end{aligned}
\end{equation*}
Since $Y_{q^j+k}\Delta_{j,k}=T_{q^{j}+k}-T^{(\rho_1)}_{j,k}$, we can write
$$\mathcal{R}_j\ \subseteq\ \mathcal{A}_{j}\ \cup\ \mathcal{B}_{j},$$
where
$$\mathcal{A}_j\ :=\ \cap_{k=1}^{r_j}\mathcal{D}_{j,k},\ \ \ \ \ \  \mathcal{B}_j\ :=\ \cap_{k=1}^{r_j}\mathcal{D}_{j,k}^C,$$
and $\mathcal{D}_{j,k}:=\left\{T_{q^{j}+k}<T^{(\rho_1)}_{j,k}\right\}$ for $k=1,..,r_j$.

The idea to prove that these events cannot occur infinitely often is the following:
consider $\mathcal{A}_j$ (for instance) and rewrite the set $\mathcal{D}_{j,r_j}$ as follows:
\begin{equation}\label{eq:second_decomposition_tau_j}
\mathcal{D}_{j,r_j}=\left\{T_{q^j+r_j}<T^{(\rho_1)}_{j,r_j}\right\}
=\left\{\sum_{i=1}^{r_j}\left(T_{q^{j}+i}-T_{q^{j}+i-1}\right)<T^{(\rho_1)}_{j,r_j}-T_{q^j}\right\},
\end{equation}
where the last inequality follows using telescopic series.
In the set $\mathcal{D}_{j,r_j}$ we have a sum of bounded random variables, i.e. $\left(T_{q^{j}+i}-T_{q^{j}+i-1}\right)$,
whose means are strictly positive on $\mathcal{A}_j$, because $\mathcal{A}_j$ in included in $\cap_{k=1}^{r_j-1}\mathcal{D}_{j,k}$;
hence, provided that the difference $(T^{(\rho_1)}_{j,r_j}-T_{q^j})$ increases with $j$ slower than $r_j$,
we could prove that the set cannot occur infinitely often.
Roughly speaking, it means that, if the adaptive threshold $\hat{\rho}_{1,q^j}$ is not far enough from the urn proportion $Z_{q^j}$, then
the average increments of the urn proportion make very likely that $Z_{q^j+k}$ crosses $\hat{\rho}_{1,q^j}$ before $q^j+r_j$.
Similar arguments apply for $\mathcal{B}_j$.
More formally, fix $\epsilon>0$ and define the set $\mathcal{C}_{j}$ as follows:
$$\mathcal{C}_{j}\ :=\ \left\{|T^{(\rho_1)}_{j,r_j}-T_{q^j}|>\epsilon j^2q^{\frac{j}{2}}\right\},$$
so that $\mathcal{C}_{j}^C$ is the set where the difference $|T^{(\rho_1)}_{j,r_j}-T_{q^j}|$ increases with $j$ slower than $r_j$.
Hence, it follows that
$$\mathcal{R}_j\ \subseteq\ \left\{\mathcal{A}_{j}-\mathcal{C}_{j}\right\}\ \cup\
\left\{\mathcal{B}_{j}-\mathcal{C}_{j}\right\}\ \cup\ \mathcal{C}_{j},$$
and the result~\eqref{eq:tauj_almsot_sure} is obtained by showing that
$$\bm{P}\left(\mathcal{A}_{j}-\mathcal{C}_{j},i.o.\right)=
\bm{P}\left(\mathcal{B}_{j}-\mathcal{C}_{j},i.o.\right)=\bm{P}\left(\mathcal{C}_{j},i.o.\right)=0.$$

We will now begin with the proof of $\bm{P}\left(\mathcal{A}_{j}-\mathcal{C}_{j},i.o.\right)=
\bm{P}\left(\mathcal{B}_{j}-\mathcal{C}_{j},i.o.\right)=0$.
From~\eqref{eq:second_decomposition_tau_j} we note that, on the set $\mathcal{C}_{j}^C$,
$$\mathcal{D}_{j,r_j}\ \subseteq\ \left\{\sum_{i=1}^{r_j}\left(T_{q^{j}+i}-T_{q^{j}+i-1}\right)<\epsilon j^2q^{\frac{j}{2}}\right\}\ =\ \mathcal{E}_{j},$$
$$\mathcal{D}_{j,r_j}^C\ \subseteq\ \left\{\sum_{i=1}^{r_j}\left(T_{q^{j}+i}-T_{q^{j}+i-1}\right)>-\epsilon j^2q^{\frac{j}{2}}\right\}\ =\ \mathcal{F}_{j}.$$
As a consequence, we have
$$\mathcal{A}_{j}-\mathcal{C}_{j}\ \subseteq\
\left\{\ \cap_{k=1}^{r_j-1}\mathcal{D}_{j,k}\cap
\mathcal{E}_{j}\ \right\},$$
$$\mathcal{B}_{j}-\mathcal{C}_{j}\ \subseteq\ \left\{\ \cap_{k=1}^{r_j-1}\mathcal{D}_{j,k}^C\cap
\mathcal{F}_{j}\ \right\}.$$
Now, consider the increments $\left(T_{q^{j}+i}-T_{q^{j}+i-1}\right)$ for $i=1,..,r_j$ contained in the
sets $\mathcal{E}_{j}$ and $\mathcal{F}_{j}$ above;
recall that
\[\begin{aligned}
(T_{q^{j}+i}-T_{q^{j}+i-1})&&=&\rho_1 \left(1-X_{q^{j}+i}\right) D_{2,q^{j}+i}W_{2,q^{j}+i-1}\\
&&-&\left(1-\rho_1\right)X_{q^{j}+i} D_{1,q^{j}+i} W_{1,q^{j}+i-1}.
\end{aligned}\]
Fix an arbitrarily small $\epsilon_1>0$ and introduce two collections of i.i.d. random variables $\left(A_1,..,A_{r_j}\right)$ and $\left(B_1,..,B_{r_j}\right)$ defined as follows:
\[\begin{aligned}
A_i\ &&:=&\ \rho_1 \left(1-\ind_{\{U_{q^{j}+i}<\rho_1+\epsilon_1\}}\right) D_{2,q^{j}+i},\\
B_i\ &&:=&\ \rho_1 \left(1-\ind_{\{U_{q^{j}+i}<\rho_1-\epsilon_1\}}\right) D_{2,q^{j}+i}\ -
\ \left(1-\rho_1\right)\ind_{\{U_{q^{j}+i}<\rho_1-\epsilon\}} D_{1,q^{j}+i},
\end{aligned}\]
where $\left(U_{q^j+1},..,U_{q^{j+1}}\right)$ are the i.i.d. $\left(0,1\right)$ uniform random variables such that $X_{q^{j}+i}:=\ind_{\{U_{q^{j}+i}<Z_{q^{j}+i-1}\}}$.

First note that, by construction, on the set $\mathcal{A}_{j}$ we have $\cap_{k=1}^{r_j}\left\{Z_{q^{j}+k}>\tilde{\rho}_{1,q^{j}+k}\right\}$,
and hence $\mathcal{A}_{j}\subset\cap_{k=1}^{r_j}\{W_{1,q^{j}+k}=0\}$.
Thus, since using~\eqref{eq:rho_exponentially_fast} we have $Z_n\stackrel{a.s.}{\rightarrow}\rho_1$ by Theorem~\ref{thm:convergence_almost_sure},
on the set $\mathcal{A}_{j}$ we have that
$$\left\{ (T_{q^{j}+i}-T_{q^{j}+i-1})\geq A_i \right\},\ \ \ i\in1,..,r_j,$$
occurs with probability 1 as $n\rightarrow\infty$.
Similarly, by construction, on the set $\mathcal{B}_{j}$ we have $\cap_{k=1}^{r_j}\left\{Z_{q^{j}+k}<\tilde{\rho}_{1,q^{j}+k}\right\}$,
and hence $\mathcal{B}_{j}\subset\cap_{k=1}^{r_j}\{W_{1,q^{j}+k}=1\}$.
Thus, using $Z_n\stackrel{a.s.}{\rightarrow}\rho_1$,
on the set $\mathcal{A}_{j}$ we have that the event
$$\left\{ (T_{q^{j}+i}-T_{q^{j}+i-1})\leq B_i \right\},\ \ \ i\in1,..,r_j,$$
occurs with probability 1 as $n\rightarrow\infty$.
As a consequence, for large $j$ we have that
\[\begin{aligned}
\bm{P}\left(\mathcal{A}_{j}-\mathcal{C}_{j},i.o.\right)\ &&\leq&\
\bm{P}\left(\sum_{i=1}^{r_j}A_i<\epsilon j^2q^{\frac{j}{2}},i.o.\right)\\
\bm{P}\left(\mathcal{B}_{j}-\mathcal{C}_{j},i.o.\right)\ &&\leq&\
\bm{P}\left(\sum_{i=1}^{r_j}B_i>-\epsilon j^2q^{\frac{j}{2}},i.o.\right).
\end{aligned}\]
Set
$$P_{Aj}:=\bm{P}\left(\ \sum_{i=1}^{r_j}A_i<\epsilon j^2q^{\frac{j}{2}}\ \right),\ \ \ \texttt{and}\ \ \
P_{Bj}:=\bm{P}\left(\ \sum_{i=1}^{r_j}B_i>-\epsilon j^2q^{\frac{j}{2}}\ \right).$$
We will now use Chernoff's upper bounds on the i.i.d. bounded random variables $A_i$ and $B_i$ (see~\eqref{eq:bernoulli_chernoff_result}).
First notice that
\begin{itemize}
\item[(1)] $\bm{E}\left[A_i\right]=\rho_1\left(1-\rho_1-\epsilon\right)m_2>0$,
\item[(2)] $\bm{E}\left[B_i\right]=\rho_1\left(1-\rho_1+\epsilon\right)m_2-\left(1-\rho_1\right)\left(\rho_1-\epsilon\right) m_1<0$,
\item[(3)] $|A_i|,|B_i|<b$ a.s. for any $i\geq1$.
\end{itemize}
Note that $P_{Aj}$ can be written as $\bm{P}\left( S_{j}\leq c_j\cdot\bm{E}[S_{j}]\right)$,
where $S_j=\sum_{i=1}^{r_j}(A_i/b)$ and
$$c_j\ =\ \frac{\epsilon j^2q^{\frac{j}{2}}}{r_j\bm{E}\left[A_1\right]/b};$$
since $c_j\rightarrow0$, we can define an integer $j_0$ such that $c_j<c_0$ for any $j\geq j_0$, so that
$$\bm{P}\left( S_j\leq c_j\cdot\bm{E}[S_j]\right)\ \leq\ \bm{P}\left( S_j\leq c_0\cdot\bm{E}[S_j]\right).$$
Hence, by using~\eqref{eq:bernoulli_chernoff_result}, for any $j\geq j_0$ we have that
$$P_{Aj}\ \leq\ \exp\left(-\frac{(1-c_0)^2}{2}\cdot\bm{E}[S_j]\right),$$
which converges to zero exponentially fast since
$$\bm{E}[S_j]\ =\ r_j\frac{\bm{E}\left[A_1\right]}{b}\ \sim\ q^{j\frac{1+\nu}{2}}.$$
We can repeat the same arguments for $P_{Bj}$, with the i.i.d. random variables $(-B_i+b)/2b\in(0,1)$ for $i=1,..,r_j$;
in this case, $c_j$ tends to a constant $c<1$, so that the proof follows with $c_0\in(c,1)$.
Thus,
$$\sum_{j=1}^{\infty} (P_{Aj}+P_{Bj})\ <\ \infty,$$
yielding
$$\bm{P}\left(\mathcal{A}_{j}-\mathcal{C}_{j},i.o.\right)\ =\ \bm{P}\left(\mathcal{B}_{j}-\mathcal{C}_{j},i.o.\right)\ =\ 0.$$

We will now show that $\bm{P}\left(\mathcal{C}_{j},i.o.\right)=0$.
Note that since $|T^{(\rho_1)}_{j,r_j}|\leq|T^{(\rho_1)}_{j,d_j}|$ and
$$T_{q^j}=Y_{q^j}(\rho_1-\hat{\rho}_{1,q^{j-1}})\ +\ Y_{q^j}(\hat{\rho}_{1,q^{j-1}}-Z_{q^j})=T^{(\rho_1)}_{j-1,d_{j-1}}+
\tilde{T}_{j-1,d_{j-1}},$$
it follows that
$$|T^{(\rho_1)}_{j,r_j}-T_{q^j}|\ \leq\ |T^{(\rho_1)}_{j,d_j}|\ +\ |T^{(\rho_1)}_{j-1,d_{j-1}}|\ +\ |\tilde{T}_{j-1,d_{j-1}}|,$$
which implies that
$$\{\mathcal{C}_{j},i.o.\}\subset\left\{|T^{(\rho_1)}_{j,d_j}|>\frac{\epsilon}{3} j^2q^{\frac{j}{2}},i.o.\right\}
\cup\left\{|\tilde{T}_{j,d_{j}}|>\frac{\epsilon}{3} j^2q^{\frac{j}{2}},i.o.\right\}.$$
Now, since $Y_{n}\leq Y_0+bn$, it follows that
$$\{\mathcal{C}_{j},i.o.\}\subset \{\mathcal{G}_{1j},i.o.\}\cup\{\mathcal{G}_{2j},i.o.\},$$
where
\[\begin{aligned}
\mathcal{G}_{1j}\ &&:=&\ \left\{\left(\frac{Y_0}{bq^{j+1}}+1\right)q^{\frac{j}{2}}|\rho_1-\hat{\rho}_{1,q^{j}}|>j^2\frac{\epsilon}{3qb}\right\},\\
\mathcal{G}_{2j}\ &&:=&\ \left\{\left(\frac{Y_0}{bq^{j+1}}+1\right)q^{\frac{j}{2}}|Z_{q^{j+1}}-\hat{\rho}_{1,q^{j}}|>j^2\frac{\epsilon}{3bq} \right\}.
\end{aligned}\]

We will now show that $\bm{P}(\mathcal{G}_{1j},i.o.)=0$. By using the Markov's inequality we have
\[\begin{aligned}
\sum_{j=1}^\infty \bm{P}\left(\mathcal{G}_{1j}\right)\ &&\leq&\
\frac{3qb}{\epsilon}\sum_{j=1}^\infty \left(\frac{Y_0}{bq^{j+1}}+1\right)\frac{\bm{E}\left[q^{\frac{j}{2}}|\rho_1-\hat{\rho}_{1,q^{j}}|\right]}{j^2}\\
&&=&\ \frac{3qb}{\epsilon}\left(\frac{Y_0}{bq}+1\right)C
\sum_{j=1}^\infty \frac{1}{j^2}\ <\ \infty,
\end{aligned}\]
where
$$C:=\sup_{k\geq 1}\bm{E}\left\{\left[q^{\frac{k}{2}}|\rho_1-\hat{\rho}_{1,q^{k}}|\right]\right\}<\infty$$
from~\eqref{eq:rho_exponentially_fast}.
Hence, using the Borel-Cantelli lemma, it follows that\\
$\bm{P}\left(\mathcal{G}_{1j},i.o.\right)=0$.

Now, consider $\mathcal{G}_{2j}$.
Let $\mathcal{H}_j:=\{j^{-2}q^{\frac{j}{2}}\cdot|\Delta_{j,d_j}|\ >\ \epsilon\}$ and since
$$\bm{P}\left(\mathcal{G}_{2j},i.o.\right)=\bm{P}\left(\mathcal{H}_j,i.o.\right)$$
we now focus on $\mathcal{H}_j$.
First, for each $j\geq1$, we recall that $\mathcal{Q}_j=\{\tau_j>d_j\}$ and we decompose $\mathcal{H}_j$ as follows:
$$\mathcal{H}_j\ \subseteq\ \mathcal{Q}_j\ \cup\ \left\{\mathcal{H}_j\cap\mathcal{Q}_j^C\right\},$$
which leads to
$$\bm{P}\left(\mathcal{H}_j,i.o.\right)\ \leq\
\bm{P}\left(\mathcal{Q}_j,i.o.\right)\ +\ \bm{P}\left(\mathcal{H}_j\cap\mathcal{Q}_j^C,i.o.\right).$$
First, consider $\bm{P}\left(\mathcal{H}_j\cap\mathcal{Q}_j^C,i.o.\right)$.
By using Markov's inequality we have
\[\begin{aligned}
\sum_{j=1}^{\infty}\bm{P}\left(\mathcal{H}_j\cap\mathcal{Q}_j^C\right)\ &&\leq&\ \sum_{j=1}^{\infty}\bm{E}\left[q^{j}\cdot|\Delta_{j,d_j}|\ind_{\mathcal{Q}_j^C}\right]
\frac{q^{-\frac{j}{2}}}{\epsilon j^2}\\
&&\leq&\ \frac{\left(\sup_{k\geq1}\left\{\bm{E}\left[q^{k}\cdot|\Delta_{k,d_k}|\ind_{\mathcal{Q}_k^C}\right]\right\}\right)}
{\epsilon}\sum_{j=1}^{\infty}\frac{q^{-\frac{j}{2}}}{j^2},
\end{aligned}\]
which is finite from Theorem~\ref{thm:how_to_use_branching_process}.
Hence, again from the Borel-Cantelli lemma we have that
$$\bm{P}\left(\mathcal{H}_j\cap\mathcal{Q}_j^C,i.o.\right)=0.$$

We will now show that $\bm{P}\left(\mathcal{Q}_j,i.o.\right)=0$.
To this end, we can follow the same arguments used in the first part of this proof, except that here we define
$$\mathcal{C}_{j}\ :=\ \left\{|T^{(\rho_1)}_{j,d_j}-T_{q^j}|>\epsilon q^j\right\}.$$
In this case, to show $\bm{P}\left(\mathcal{C}_{j},i.o.\right)=0$
we have to prove that the following two events cannot occur infinitely often
\begin{itemize}
\item[(i)] $\mathcal{G}_{3j}:=\left\{\left(\frac{Y_0}{bq^{j+1}}+1\right)|\rho_1-\hat{\rho}_{1,q^{j}}|>\frac{\epsilon}{2qb}\right\},$
\item[(ii)] $\mathcal{G}_{4j}:=\left\{\left(\frac{Y_0}{bq^{j}}+1\right)|\rho_1-Z_{q^{j}}|>\frac{\epsilon}{2b} \right\}.$
\end{itemize}
Result (i) is implied by~\eqref{eq:rho_exponentially_fast}, while (ii) follows from Theorem~\ref{thm:convergence_almost_sure}.
Hence, we have that
$$\bm{P}\left(\mathcal{C}_{j},i.o.\right)=0.$$
Then, similarly to the first part of the proof, we deal with the sets $\mathcal{A}_{j}-\mathcal{C}_{j}$ and $\mathcal{B}_{j}-\mathcal{C}_{j}$ by applying
Chernoff's upper bound to the probabilities
$$P_{Aj}=\bm{P}\left(\ \sum_{i=1}^{d_j}A_i<\epsilon q^{j}\ \right)\ \ \ and
\ \ \ P_{Bj}=\bm{P}\left(\ \sum_{i=1}^{d_j}B_i>-\epsilon q^{j}\ \right),$$
which implies $\sum_{j=1}^{\infty} P_{Aj}<\infty$ and $\sum_{j=1}^{\infty}P_{Bj}<\infty$.
Hence, from the Borel-Cantelli lemma we get
$$\bm{P}\left(\mathcal{A}_{j}-\mathcal{C}_{j},i.o.\right)\ =\ \bm{P}\left(\mathcal{B}_{j}-\mathcal{C}_{j}\right)\ =\ 0.$$
which implies $\bm{P}\left(\mathcal{Q}_j,i.o.\right)=0$.
This concludes the proof.
\endproof

\begin{rem}
The result of Theorem~\ref{thm:tauj} continues to hold if~\eqref{eq:rho_exponentially_fast} is not satisfied,
but~\eqref{eq:rho_convergence_almost_sure} and condition (c1) hold.
Moreover, since in the proof we use Theorem~\ref{thm:Y_bounded_and_diverge},
if~\eqref{eq:rho_exponentially_fast} does not hold condition (c2) must be assumed (see Remark~\ref{rem:Y_bounded_and_diverge}).
\end{rem}

\subsection{Proof of Proposition~\ref{prop:rhomean_closed_rho}}   \label{subsection_proof_proposition_closed_rho}

\proof
Wlog assume $m_1>m_2$, which implies $\hat{\rho}_n=\hat{\rho}_{1,n}$ and $\rho=\rho_1$.
First, we have
\begin{equation}\label{eq:bias_first_eq}
\bm{E}\left[n|\bar{\rho}_{1,n}-\rho_1|^2\right]\ =\
\frac{1}{n}\bm{E}\left[\left|\sum_{i=0}^{n-1}\left(\tilde{\rho}_{1,i}-\rho_1\right)\right|^2\right],
\end{equation}
and note that
\[\begin{aligned}\sum_{i=0}^{n-1}\left(\tilde{\rho}_{1,i}-\rho_1\right)\ &&=&\
\sum_{j=0}^{k_n}\sum_{i=0}^{d_j}
\left(\tilde{\rho}_{1,q^j+i}-\rho_1\right)\ind_{\{q^{k_n}+i\leq n\}}\\
&&=&\ \sum_{j=0}^{k_n-1}d_j\left(\hat{\rho}_{1,q^j}-\rho_1\right)+\left(n-q^{k_n}\right)\left(\hat{\rho}_{1,q^{k_n}}-
\rho_1\right),
\end{aligned}\]
where we recall $k_n$ is defined in~\eqref{def:k_n} as $k_n:=[log_{q}(n)]$.
Since $d_j=(q-1)q^j$, the LHS of~\eqref{eq:bias_first_eq} is equal to
$$ \frac{\left(q-1\right)^2}{n}\bm{E}\left[\left|\sum_{j=0}^{k_n-1}\left(\sqrt{q}\right)^j\cdot\left(\sqrt{q}^j\left(\hat{\rho}_{1,q^j}-\rho_1\right)\right)+
\left(\frac{n-q^{k_n}}{q-1}\right)\left(\hat{\rho}_{1,q^{k_n}}-\rho_1\right)\right|^2\right],$$
and, defining $c_j:=\sqrt{q}^j|\hat{\rho}_{1,q^j}-\rho_1|$, we can rewrite the last expression as follows:
$$\frac{\left(q-1\right)^2}{n}\bm{E}\left[\left(\sum_{j=0}^{k_n-1}\left(\sqrt{q}\right)^j\cdot c_j+
\left[\frac{n-q^{k_n}}{\sqrt{q}^{k_n}\left(q-1\right)}\right]c_{k_n}
\right)^2\right];$$
Now, using Cauchy Schwartz inequality and using $\left(\frac{n-q^{k_n}}{\sqrt{q}^{k_n}\left(q-1\right)}\right)\leq \sqrt{q}^{k_n}$,
the above expectation is less than or equal to
$$K_n\ :=\ \sum_{j_1=0}^{k_n}\sum_{j_2=0}^{k_n}\left(\sqrt{q}\right)^{j_1}\left(\sqrt{q}\right)^{j_2}\cdot
 \sqrt{\bm{E}\left[c_{j_1}^2\right]\bm{E}\left[c_{j_2}^2\right]}.$$
Now, by the symmetry in $K_n$, we can use the following decomposition
\[\begin{aligned}
\sum_{j_1=0}^{k_n}\sum_{j_2=0}^{k_n}\left(\cdot\right)\ &&=&\
\sum_{j_1=0}^{\sqrt{k_n}}\sum_{j_2=0}^{\sqrt{k_n}}\left(\cdot\right)\ +\
2\sum_{j_1=0}^{\sqrt{k_n}}\sum_{j_2=\sqrt{k_n}}^{k_n}\left(\cdot\right)\ +\
\sum_{j_1=\sqrt{k_n}}^{k_n}\sum_{j_2=\sqrt{k_n}}^{k_n}\left(\cdot\right)\\
&&\leq&\ 2\sum_{j_1=0}^{k_n}\sum_{j_2=0}^{\sqrt{k_n}}\left(\cdot\right)\ +\
\sum_{j_1=\sqrt{k_n}}^{k_n}\sum_{j_2=\sqrt{k_n}}^{k_n}\left(\cdot\right),\\
\end{aligned}\]
we obtain
\[\begin{aligned}
K_n\ &&\leq&\ \sup_{j\geq 1}\left\{\bm{E}\left[c_{j}^2\right]\right\}\cdot2
\sum_{j_1=0}^{k_n}\sum_{j_2=0}^{\sqrt{k_n}}\left(\sqrt{q}\right)^{j_1}\left(\sqrt{q}\right)^{j_2}\\
&&+&\  \max_{\sqrt{k_n}\leq j\leq k_n}\left\{\bm{E}\left[c_{j}^2\right]\right\}\cdot
\sum_{j_1=\sqrt{k_n}}^{k_n}\sum_{j_2=\sqrt{k_n}}^{k_n}\left(\sqrt{q}\right)^{j_1}\left(\sqrt{q}\right)^{j_2}\\
&&=&\ K_{1n}\ +\ K_{2n}.
\end{aligned}\]
Now, consider $K_{1n}$; we have that
$$K_{1n}\leq\sup_{j\geq 1}\left\{\bm{E}\left[c_{j}^2\right]\right\}\cdot2
\left(\frac{\left(\sqrt{q}\right)^{\sqrt{k_n}+1}-1}{\sqrt{q}-1}\right)
\left(\frac{\left(\sqrt{q}\right)^{k_n+1}-1}{\sqrt{q}-1}\right),$$
and by multiplying for $(q-1)^2/n$ we obtain
$$2\left(\frac{q-1}{\sqrt{q}-1}\right)^2\sup_{j\geq 1}\left\{\bm{E}\left[c_{j}^2\right]\right\}\cdot
\left(\frac{\left(\sqrt{q}\right)^{\sqrt{k_n}+1}-1}{\sqrt{n}}\right)
\left(\frac{\left(\sqrt{q}\right)^{k_n+1}-1}{\sqrt{n}}\right).$$
Using~\eqref{eq:rho_exponentially_fast} we have that $\sup_{j\geq 1}\left\{\bm{E}\left[c_{j}^2\right]\right\}$ is finite.
Moreover, since $n\leq q^{k_n+1}$ by definition of $k_n$, we have that
$$\left(\frac{\left(\sqrt{q}\right)^{k_n+1}-1}{\sqrt{n}}\right)\leq\sqrt{q},\ \ \ \left(\frac{\left(\sqrt{q}\right)^{\sqrt{k_n}+1}-1}{\sqrt{n}}\right)\ \rightarrow\ 0.$$

Similarly, we can consider $K_{2n}$ and write
$$K_{2n}\leq\max_{\sqrt{k_n}\leq j\leq k_n}\left\{\bm{E}\left[c_{j}^2\right]\right\}\cdot
\left(\frac{\left(\sqrt{q}\right)^{k_n+1}-1}{\sqrt{q}-1}\right)^2.$$
Then, by multiplying for $(q-1)^2/n$ we obtain
$$\left(\frac{q-1}{\sqrt{q}-1}\right)^2\max_{\sqrt{k_n}\leq j\leq k_n}\left\{\bm{E}\left[c_{j}^2\right]\right\}\cdot
\left(\frac{\left(\sqrt{q}\right)^{k_n+1}-1}{\sqrt{n}}\right)^2,$$
and from~\eqref{eq:rho_exponentially_fast} and $n\leq q^{k_n+1}$ we have $\max_{\sqrt{k_n}\leq j\leq k_n}\left\{\bm{E}\left[c_{j}^2\right]\right\}$ is finite and
$$\left(\frac{\left(\sqrt{q}\right)^{k_n+1}-1}{\sqrt{n}}\right)^2\ \leq\ q.$$

Then, combining all together we obtain
\[\begin{aligned}
\limsup_{n\rightarrow\infty}\bm{E}\left[n|\bar{\rho}_{1,n}-\rho_1|^2\right]\ &&\leq&\ \limsup_{n\rightarrow\infty}\frac{(q-1)^2}{n}K_{2n}\\
&&\leq&\
q\left(1+\sqrt{q}\right)^2\cdot\limsup_{n\rightarrow\infty}\bm{E}\left[n|\hat{\rho}_{1,n}-\rho_1|^2\right],
\end{aligned}\]
which is finite because of condition~\eqref{eq:rho_exponentially_fast}.
\endproof

\subsection{Proof of Corollary~\ref{cor:MRRU_CLT}}   \label{subsection_proof_corollary_MRRU_CLT}

\proof
To prove this result, we apply Theorem~\ref{thm:convergence_in_distribution} to the urn model with fixed thresholds, i.e.
$\tilde{\rho}_{1,n}=\rho_1$ and $\tilde{\rho}_{2,n}=\rho_2$ for all $n\geq0$, since in this case $\bar{\rho}_n=\rho$ for all $n\geq0$.
\endproof

\subsection{Remarks on the CLT for $Z_n$}   \label{subsection_proof_theorem_no_CLT_Z}

In this subsection, we discuss the second-order behavior of the proportion $Z_n$ of balls in the urn in the ARRU model.
Specifically, we establish CLT of $Z_{n_j}$ for some specific subsequences $\{n_j;j\geq1\}$ and
we highlight the challenges to the proof of a full CLT.
To this end, let us assume that
\begin{equation}\label{eq:rho_converge_normal}
\sqrt{n}(\hat{\rho}_{n}-\rho)\ \stackrel{d}{\rightarrow}\ Z,
\end{equation}
where $Z$ is a Gaussian random variable with zero mean and variance $\sigma^2>0$.
It is worth noticing that~\eqref{eq:rho_converge_normal} is usually verified in applications,
since $\hat{\rho}_{n}$ is typically a continuous function of maximum likelihood estimators, e.g. see~\eqref{eq:f_clinical_trials}.

Wlog, assume $m_1>m_2$ and consider the sequence $\left\{\sqrt{n}(\rho_1-Z_n);n\geq1\right\}$.

Now, along the subsequence $\{q^j;j\geq1\}$
\begin{equation}\label{eq:first_subsequence}
\sqrt{q^j}(\rho_1-Z_{q^j})\ =\ \sqrt{q^j}(\rho_1-\hat{\rho}_{1,q^{j-1}})\ +\ \sqrt{q^j}(\hat{\rho}_{1,q^{j-1}}-Z_{q^j}).
\end{equation}
Using~\eqref{eq:rho_converge_normal}, it follows that
$$\sqrt{q^j}(\rho_1-\hat{\rho}_{1,q^{j-1}})\ =\ \sqrt{q}\cdot \left[\sqrt{q^{j-1}}(\rho_1-\hat{\rho}_{1,q^{j-1}})\right]\
\stackrel{d}{\rightarrow}\ \sqrt{q}\cdot Z.$$
As for the second term in~\eqref{eq:first_subsequence}, it can be expressed as
\[\begin{aligned}\sqrt{q^j}(\hat{\rho}_{1,q^{j-1}}-Z_{q^j})\ &&=&\ \sqrt{q^j}\Delta_{j-1,d_{j-1}}\\
&&=&\ \sqrt{q^j}\Delta_{j-1,d_{j-1}}\ind_{\mathcal{R}_{j-1}}\ +\ \sqrt{q^j}\Delta_{j-1,d_{j-1}}\ind_{\mathcal{R}_{j-1}^c},
\end{aligned}\]
where we recall $\mathcal{R}_j:=\{\tau_j>r_j\}$ and $r_j:=q^{j(1+\nu)/2}$, with $\nu\in(0,1/2)$.
The first term converges to zero a.s. from Theorem~\ref{thm:tauj}, while using $r_{j-1}\leq d_{j-1}$ from Theorem~\ref{thm:how_to_use_branching_process} we have
$$\sqrt{q^j}\bm{E}\left[|\Delta_{j-1,d_{j-1}}|\ind_{\mathcal{R}_{j-1}^c}\right]\ \leq\
\sqrt{q^j}\bm{E}\left[|\Delta_{j-1,d_{j-1}}|\ind_{\{\tau_{j-1}\leq d_{j-1}\}}\right]\
\rightarrow\ 0.$$
Thus,
$$\sqrt{q^j}(\rho_1-Z_{q^j})\stackrel{d}{\rightarrow}\sqrt{q}\cdot Z.$$

Now, consider the subsequence $\{q^j+r_j;j\geq1\}$.
Now,
\begin{equation}\label{eq:second_subsequence}
\sqrt{q^j}(\rho_1-Z_{q^j+r_j})\ =\ \sqrt{q^j}(\rho_1-\hat{\rho}_{1,q^{j}})\ +\ \sqrt{q^j}(\hat{\rho}_{1,q^{j}}-Z_{q^j+r_j}).
\end{equation}
As before, from~\eqref{eq:rho_converge_normal} we have that
$$\sqrt{q^j}(\rho_1-\hat{\rho}_{1,q^{j}})\ \stackrel{d}{\rightarrow}\ Z.$$
Once again expressing the second term in~\eqref{eq:first_subsequence} below,
\[\begin{aligned}\sqrt{q^j}(\hat{\rho}_{1,q^{j}}-Z_{q^j+r_j})\ &&=&\ \sqrt{q^j}\Delta_{j,r_j}\\
&&=&\ \sqrt{q^j}\Delta_{j,r_j}\ind_{\mathcal{R}_j}\ +\ \sqrt{q^j}\Delta_{j,r_j}\ind_{\mathcal{R}_j^c},
\end{aligned}\]
one can show the first term in the RHS converges to zero a.s. from Theorem~\ref{thm:tauj}.
The second term in the RHS tends to zero in $L^1$
from Theorem~\ref{thm:how_to_use_branching_process}.
Thus,
$$\sqrt{q^j+r_j}(\rho_1-Z_{q^j+r_j})\stackrel{d}{\rightarrow}Z.$$

A crucial result to obtain CLT for $\{Z_{n_j};j\geq1\}$, with $n_j=q^j$ and $n_j=q^j+r_j$, is Theorem~\ref{thm:tauj}, which establishes $\bm{P}(\tau_j>n_j,i.o.)=0$.
From these results, it follows that the asymptotic distribution of $Z_{n_j}$ only involves times $t\in\cup_{j}(q^j+\tau_j,q^{j+1})$;
for these times, Theorem~\ref{thm:how_to_use_branching_process} establishes a uniform bound for $n|Z_n-\tilde{\rho}_n|$.
However, at times $t\in \cup_{j}(q^j,q^j+\tau_j)$ it seems difficult to obtain the detailed behavior of $(Z_n-\tilde{\rho}_n)$.
This gap needs to be handled for a CLT for $\{Z_n;n\geq1\}$.
This is beyond the scope of the current paper.

\section{Simulation studies}   \label{section_simulation}

In this section, we describe some simulation studies that illustrate the theoretical results presented in Section~\ref{section_the_model} in the
context of clinical trials.
We recall from subsection~\ref{subsection_choices_f} that, in the context of clinical trials,
the random variables $\xi_{1,n}$ and $\xi_{2,n}$ are interpreted as potential responses to competing treatments $\mathcal{T}_1$ and $\mathcal{T}_2$,
whose distributions $\mu_1$ and $\mu_2$ depend on parameters $\bm{\theta_1}$ and $\bm{\theta_2}$ respectively.
Let $\bm{\theta}=(\bm{\theta}_1,\bm{\theta}_2)$.
Now, letting $f_{1}$ and $f_{2}$ are two continuous functions, we recall that $\rho_1=f_1(\bm{\theta})$ and $\rho_2=f_2(\bm{\theta})$.
Moreover, the adaptive thresholds $\hat{\rho}_{1,n}$ and $\hat{\rho}_{2,n}$ are defined as follows:
$$\hat{\rho}_{1,n}:=f_{1}\left(\bm{\hat{\theta}}_{1,n}\right)\ \ \ \ \ \textit{and}\ \ \ \ \ \hat{\rho}_{2,n}:=f_{2}\left(\bm{\hat{\theta}}_{2,n}\right), \ \ \forall n\geq1,$$
where $\bm{\hat{\theta}}_{1,n}$ and $\bm{\hat{\theta}}_{2,n}$ are the adaptive estimators of $\bm{\theta_1}$ and $\bm{\theta_2}$
after the first $n$ allocations.

The main goal of this section is to illustrate the asymptotic behavior of the allocation proportion $N_{1,n}/n$
and of the parameter estimators $\bm{\hat{\theta}_n}$.
Simulations are performed with $N=10^5$ independent urn processes, each which evolve following the model described in Section~\ref{section_the_model} with
adaptive thresholds $\tilde{\rho}_{1,n}$ and $\tilde{\rho}_{2,n}$ that change at exponential times $\{q^j;j\geq1\}$, with $q=1.25$, [see~\eqref{def:adaptive_thresholds}].
For all the $N$ urn processes we used initial composition $\left(y_{1,0},y_{2,0}\right)=\left(2,2\right)$ and sample size $n=200$.
The functions $f_1$ and $f_2$ are chosen as in~\eqref{eq:linear_combination} with $p=0.75$.
We analyze both Bernoulli and Gaussian responses.\\

\subsection{Bernoulli responses}
We assume responses to treatments $\mathcal{T}_1$ and $\mathcal{T}_2$
are from Bernoulli distributions with parameters $p_1$ and $p_2$, respectively.
In this case, $\bm{\theta}=\left(p_1,p_2\right)$ is the vector of unknown parameters.
We examine two target allocations:
\begin{itemize}
\item[(a)] $\eta(\bm{\theta})=\left(1-p_1\right)/\left(2-p_1-p_2\right)$, proposed by~\cite{Wei.et.al.78};
\item[(b)] $\eta(\bm{\theta})=\sqrt{p_1}/\left(\sqrt{p_1}+\sqrt{p_2}\right)$, proposed by~\cite{Rosenberger.et.al.01}.
\end{itemize}
Hence, from~\eqref{eq:linear_combination} with $p=0.75$, we have
$$\rho_1= 0.25\cdot 1 +0.75\cdot \eta(p_1,p_2), \qquad \text{and} \qquad \rho_2= 0.25\cdot 0 +0.75\cdot \eta(p_1,p_2).$$
In Table~\ref{tab:Bernoulli} we report the simulation results on the mean and standard error of the allocation proportion $N_{1,n}/n$
and of the estimators $\hat{p}_{1,n}$ and $\hat{p}_{2,n}$, defined as
$$\hat{p}_{1,n}=\frac{\sum_{i=1}^n X_i \xi_{1,i}}{N_{1,n}}, \qquad \text{and} \qquad \hat{p}_{2,n}=\frac{\sum_{i=1}^n \left(1-X_i\right) \xi_{2,i}}{N_{2,n}}.$$
Hence,
$$\hat{\rho}_{1,n}= 0.25\cdot 1 +0.75\cdot \eta(\hat{p}_{1,n},\hat{p}_{2,n}),$$
and
$$\hat{\rho}_{2,n}= 0.25\cdot 0 +0.75\cdot \eta(\hat{p}_{1,n},\hat{p}_{2,n}).$$

\begin{table}[t]
\begin{center}
\caption{\small Simulation of $N_{1,n}/n$ and $\bm{\hat{\theta}}_n$ are given for different designs, with mean square errors given in parenthesis.
 The target allocation is $\rho_1= (1-p)\cdot 1 +p\cdot \eta(p_1,p_2)$ with $p=0.75$.
Simulation used $N=10^5$ ARRU processes with $n=200$ and changes at times $\{q^j;j\geq1\}$ with $q=1.25$.
Initial composition $\left(y_{1,0},y_{2,0}\right)=(2,2)$.}
\begin{tabular}{|cc|cccc|}
  \hline
  $p_1$ & $p_2$ & $\rho_1$ & $N_{1,n}/n$ & $\hat{p}_{1,n}$ & $\hat{p}_{2,n}$ \\\hline
   & & \multicolumn{4}{c|}{(a) $\eta=(1-p_1)/(2-p_1-p_2)$}\\\hline
  0.9 & 0.7 & 0.44 & 0.44(0.07) & 0.89(0.03) & 0.7(0.04)\\
  0.9 & 0.5 & 0.38 & 0.41(0.06) & 0.89(0.03) & 0.50(0.05)\\
  0.9 & 0.3 & 0.34 & 0.40(0.07) & 0.89(0.03) & 0.30(0.04)\\
  0.9 & 0.1 & 0.33 & 0.43(0.12) & 0.89(0.03) & 0.11(0.03)\\
  0.7 & 0.5 & 0.53 & 0.50(0.07) & 0.70(0.05) & 0.50(0.05)\\
  0.7 & 0.3 & 0.48 & 0.48(0.05) & 0.70(0.05) & 0.30(0.04)\\
  0.7 & 0.1 & 0.44 & 0.48(0.06) & 0.70(0.05) & 0.11(0.03)\\
  0.5 & 0.3 & 0.56 & 0.53(0.06) & 0.50(0.05) & 0.30(0.05)\\
  0.5 & 0.1 & 0.52 & 0.53(0.04) & 0.50(0.05) & 0.11(0.03)\\
  0.3 & 0.1 & 0.58 & 0.56(0.05) & 0.30(0.04) & 0.11(0.03)\\\hline
   & & \multicolumn{4}{ c| }{(b) $\eta=\sqrt{p_1}/(\sqrt{p_1}+\sqrt{p_2})$}\\\hline
  0.9 & 0.7 & 0.65 & 0.57(0.11) & 0.89(0.03) & 0.69(0.05)\\
  0.9 & 0.5 & 0.68 & 0.63(0.08) & 0.89(0.03) & 0.50(0.06)\\
  0.9 & 0.3 & 0.73 & 0.69(0.06) & 0.89(0.03) & 0.30(0.06)\\
  0.9 & 0.1 & 0.81 & 0.76(0.07) & 0.89(0.02) & 0.11(0.04)\\
  0.7 & 0.5 & 0.66 & 0.58(0.11) & 0.69(0.04) & 0.50(0.06)\\
  0.7 & 0.3 & 0.70 & 0.66(0.07) & 0.70(0.04) & 0.30(0.06)\\
  0.7 & 0.1 & 0.79 & 0.74(0.07) & 0.70(0.04) & 0.12(0.04)\\
  0.5 & 0.3 & 0.67 & 0.60(0.10) & 0.50(0.05) & 0.30(0.05)\\
  0.5 & 0.1 & 0.77 & 0.70(0.08) & 0.50(0.04) & 0.11(0.04)\\
  0.3 & 0.1 & 0.73 & 0.64(0.11) & 0.30(0.04) & 0.11(0.03)\\\hline
\end{tabular}
\end{center}
\label{tab:Bernoulli}
\end{table}

\medskip

\subsection{Gaussian responses}
We now assume responses to treatments $\mathcal{T}_1$ and $\mathcal{T}_2$ are from a Gaussian distribution with parameters $\left(m_1,\sigma_1^2\right)$ and $\left(m_2,\sigma_2^2\right)$, respectively.
In this case, $\bm{\theta}=\left(m_1,\sigma_1^2,m_2,\sigma_2^2\right)$ is the vector of unknown parameters.
We examine two target allocation:
\begin{itemize}
\item[(c)] $\eta(\bm{\theta})=\sigma_1/\left(\sigma_1+\sigma_2\right)$, used in~\cite{Hu.et.al.06};
\item[(d)] $\eta(\bm{\theta})=\sigma_1\sqrt{m_2}/\left(\sigma_1\sqrt{m_2}+\sigma_2\sqrt{m_1}\right)$, proposed by~\cite{Zhang.et.al.06}.
\end{itemize}
Hence, from~\eqref{eq:linear_combination} with $p=0.75$, we have
$$\rho_1= 0.25\cdot 1 +0.75\cdot \eta(\bm{\theta}), \qquad \text{and} \qquad \rho_2= 0.25\cdot 0 +0.75\cdot \eta(\bm{\theta}).$$
In Table~\ref{tab:normal} we report the simulation results on the mean and standard error of the allocation proportion $N_{1,n}/n$
and the parameter estimators $\hat{\sigma^2}_{1,n}$ and $\hat{\sigma^2}_{2,n}$, defined as
$$\hat{\sigma^2}_{1,n}=\frac{\sum_{i=1}^n X_i\left( \xi_{1,i}-\hat{m}_{1,n}\right)^2}{N_{1,n}} \qquad \text{and}$$
$$\hat{\sigma^2}_{2,n}=\frac{\sum_{i=1}^n \left(1-X_i\right)\left( \xi_{2,i}-\hat{m}_{2,n}\right)^2}{N_{2,n}},$$
where $\hat{m}_{1,n}=\sum_{i=1}^n X_i\xi_{1,i} / N_{1,n}$ and
$\hat{m}_{2,n}=\sum_{i=1}^n (1-X_i)\xi_{2,i} / N_{2,n}$.
Hence,
$$\hat{\rho}_{1,n}= 0.25\cdot 1 +0.75\cdot \eta(\bm{\hat{\theta}}_n),$$
and
$$\hat{\rho}_{2,n}= 0.25\cdot 0 +0.75\cdot \eta(\bm{\hat{\theta}}_n).$$

\begin{table}[t]
\begin{center}
\caption{\small Simulations of $N_{1,n}/n$ and $\bm{\hat{\theta}}_n$ are given for different designs, with mean square errors given in parenthesis.
 The target allocation is $\rho_1= (1-p)\cdot 1 +p\cdot \eta(\bm{\theta})$ with $p=0.75$.
Simulation used $N=10^5$ ARRU processes with $n=200$ and changes at times $\{q^j;j\geq1\}$ with $q=1.25$.
Initial composition $\left(y_{1,0},y_{2,0}\right)=(2,2)$.}
\begin{tabular}{|cccc| cccc|}
  \hline
  $m_1$ & $m_2$ & $\sigma_1^2$ & $\sigma_2^2$ & $\rho_1$ & $N_{1,n}/n$ & $\hat{\sigma^2}_{1,n}$ & $\hat{\sigma^2}_{2,n}$ \\\hline
   \multicolumn{8}{|c|}{(c) $\eta=\sigma_1/\left(\sigma_1+\sigma_2\right)$}\\\hline
  10 & 5 & 1 & 1 & 0.63 & 0.61(0.05) & 1.01(0.13) & 1.01(0.16)\\
   8 & 5 & 1 & 1 & 0.63 & 0.59(0.07) & 1.01(0.13) & 1.01(0.16)\\
   6 & 5 & 1 & 1 & 0.63 & 0.55(0.12) & 1.01(0.14) & 1.01(0.15)\\
  10 & 5 & 4 & 1 & 0.75 & 0.73(0.06) & 4.00(0.47) & 1.01(0.20)\\
   8 & 5 & 4 & 1 & 0.75 & 0.71(0.07) & 4.00(0.48) & 1.01(0.19)\\
   6 & 5 & 4 & 1 & 0.75 & 0.66(0.13) & 4.03(0.50) & 1.01(0.18)\\
  10 & 5 & 1 & 4 & 0.50 & 0.49(0.05) & 1.01(0.14) & 4.00(0.57)\\
   8 & 5 & 1 & 4 & 0.50 & 0.48(0.07) & 1.01(0.15) & 4.03(0.56)\\
   6 & 5 & 1 & 4 & 0.50 & 0.43(0.11) & 1.01(0.16) & 4.03(0.54)\\ \hline
   \multicolumn{8}{|c|}{(d) $\eta=\sigma_1\sqrt{m_2}/(\sigma_1\sqrt{m_2}+\sigma_2\sqrt{m_1})$}\\\hline
  10 & 5 & 1 & 1 & 0.56 & 0.55(0.05) & 1.01(0.14) & 1.01(0.15)\\
   8 & 5 & 1 & 1 & 0.58 & 0.55(0.07) & 1.01(0.14) & 1.01(0.15)\\
   6 & 5 & 1 & 1 & 0.61 & 0.53(0.12) & 1.01(0.14) & 1.01(0.15)\\
  10 & 5 & 4 & 1 & 0.69 & 0.67(0.06) & 4.03(0.49) & 1.01(0.18)\\
   8 & 5 & 4 & 1 & 0.71 & 0.67(0.07) & 4.03(0.49) & 1.01(0.18)\\
   6 & 5 & 4 & 1 & 0.73 & 0.65(0.13) & 4.03(0.51) & 1.01(0.18)\\
  10 & 5 & 1 & 4 & 0.45 & 0.44(0.05) & 1.01(0.15) & 4.03(0.54)\\
   8 & 5 & 1 & 4 & 0.46 & 0.44(0.07) & 1.01(0.16) & 4.03(0.54)\\
   6 & 5 & 1 & 4 & 0.48 & 0.42(0.11) & 1.01(0.16) & 4.03(0.53)\\ \hline
\end{tabular}
\end{center}
\label{tab:normal}
\end{table}

The results show that our methods target the true parameters effectively.
In real clinical trials, further calibration may be performed to reduce small bias.

\section{Extensions to multi-color urn models}   \label{section_discussion}

It is important to note that all the results presented in this paper can be extended to the case of $K>2$ colors,
when $\exists j\in\{1,..,K\}$ such that $m_j>m_k$ for any $k\neq j$.
In the context of clinical trials, the functions $f_j$ should be interpreted as the target allocations for $N_{j,n}/n$
when $\mathcal{T}_j$ is the superior treatment, and the variables $W_{j,n}$ should be all defined as $\ind_{\{Z_{n\leq\hat{\rho}_{j,n}}\}}$.

\section{Acknowledgments}
The authors thank Prof. A.M. Paganoni of Politecnico di Milano for stimulating discussions and ideas on the urn models considered in the paper.
The authors thank also Prof. G. Aletti of Universit\'{a} degli Studi di Milano for supporting this research and several useful discussions.
Part of Andrea Ghiglietti's work was carried out while he was a doctoral student, Department of Mathematics, Politecnico di Milano.
This research was started while Andrea Ghiglietti was a doctoral student visiting the Department of Statistics, George Mason University.
He thanks the Department of Statistics for its hospitality.
Part of Prof. Rosenberger's research was conducted while he was a Visiting Scholar in the
Department of Mathematics, University of Southern California. He thanks the Department for its hospitality.
Part of Anand Vidyashankar's work was carried out when he was visiting the Department of Mathematics, Universit\`{a} degli Studi di Milano.
He thanks the Department for its hospitality.

\end{document}